\documentclass[11pt,a4paper,fleqn]{article}
\usepackage{amsfonts,amsmath}
\usepackage{latexsym}
\usepackage{amssymb}
\usepackage{euscript}
\usepackage{graphicx}
\usepackage{a4wide}

\newtheorem{prop}{Proposition}[section]

\newtheorem{lemma}[prop]{Lemma}

\newtheorem{rem}[prop]{Remark}
\newtheorem{theo}[prop]{Theorem}

\newcommand{\dis}{\displaystyle }
\newcommand{\noi}{\noindent}

\newenvironment{proof}[1]{\begin{trivlist}\item {\it
\bf Proof.}\quad} {\qed\end{trivlist}}
\newenvironment{prooff}[1]{\begin{trivlist}\item {\it
\bf Proof}\quad} {\qed\end{trivlist}}

\newcommand{\qed}{\nopagebreak\hspace*{\fill}
{\vrule width6pt height6ptdepth0pt}\par}

\begin{document}


\title{\bf LIMITING LAWS ASSOCIATED WITH  BROWNIAN MOTION PERTURBED
BY ITS MAXIMUM, MINIMUM AND LOCAL TIME II}
\author{{\small{\text{\bf Bernard ROYNETTE}$^{(1)}$}},\
{\small {\text{\bf  Pierre VALLOIS}$^{(1)}$}} {\small and} {\small
{\text{\bf Marc YOR }$^{(2),(3)}$}}}



\maketitle {\small

\noindent (1)\,\, Universit\'e Henri Poincar\'e, Institut de
Math\'ematiques Elie Cartan, B.P. 239, F-54506 Vand\oe uvre-l\`es-Nancy Cedex\\

\noi (2)\,\, Laboratoire de Probabilit\'{e}s et
 Mod\`{e}les Al\'{e}atoires, Universit\'{e}s Paris VI et VII -  4, Place Jussieu
 - Case 188 -
 F-75252 Paris Cedex 05.\\

 \noi (3) Institut Universitaire de France.\\



\vskip 40 pt \noi {\bf Abstract.} Let $P _0$ denote the Wiener
measure defined on the canonical space $\big( \Omega={\cal
C}(\mathbb{R} _+, \mathbb{R}),$ $ (X_t)_{t \geq 0},\; ({\cal
F}_t)_{t \ge 0}\big)$, and
 $(S_t)$ (resp. $(I_t)$), be the one sided-maximum (resp.
minimum), $(L^0_t)$ the local time at $0$, and  $(D_t)$ the number
of down-crossings from $b$ to $a$ (with $b>a$).
%
Let $f : \mathbb{R} \times \mathbb{R}^d \longrightarrow ]0, +
\infty[$ be a Borel function, and $(A_t)$ be a process chosen
within the set : $\big\{(S_t);\ (S_t,t);\ (L^0_t);\
(S_t,I_t,L^0_t);\ (D_t)\big\}$, which consists of 5 elements. We
prove a penalization result : under some suitable assumptions on
$f$, there exists a positive $\big( ({\cal
F}_t),P_0\big)$-martingale $(M_t^{f})$, starting at $1$, such that
:
\begin{equation}\label{ab1}
    \lim_{t\rightarrow \infty}\frac{E_0\big[1_{\Gamma
_s}f(X_t,A_t)\big]}{E_0[f(X_t,A_t)]}=Q_0^{f}(\Gamma _s):=E_0\big[
1_{\Gamma _s} M_s^{f}\big],\quad \forall \Gamma _s \in {\cal F}_s
\ \mbox{and }\ s\geq 0.
\end{equation}
We determine the law of $(X_t)$ under the  p.m. $Q_0^{f}$ defined
on $\big(\Omega, {\cal F}_\infty\big)$  by (\ref{ab1}). For the
$1^{st}, 3^{rd}$ and $5^{th}$ elements of the set, we prove first
that $Q_0^{f} (A_{\infty} < \infty)=1$, and more generally
$Q_0^{f} (0<g < \infty)=1$ where $g={\rm sup} \{s>0, \,
A_s=A_{\infty}\}$ (with the convention $\sup \emptyset=0$).
Secondly, we  split the
 trajectory of $(X_t)$ in two parts : $(X_t)_{0\leq t \leq g}$ and
$(X_{t+g})_{ t \geq 0}$, and we describe their laws  under
$Q_0^{f}$, conditionally on $A_{\infty}$. For the $2^{nd}$ and
$4^{th}$ elements,  a similar result holds replacing $A_\infty$ by
resp. $S_\infty$,  $S_\infty \vee I_\infty$.


\vskip 40pt


\noi {\bf Key words and phrases} : penalization, enlargement of
filtration, maximum, minimum, local time, down-crossings.
\smallskip



\noi {\bf AMS 2000 subject classifications} :  60 B 10, 60 G 17,
60 G 40, 60 G 44, 60 J 25, 60 J 35, 60 J 55, 60 J 60, 60 J 65.

\section{Introduction}\label{int}

\setcounter{equation}{0}

 \noi {\bf 1.1} Let $(P _x)_{x \in
\mathbb{R}}$ be the family of Wiener measures defined on the
canonical space $\big(\Omega = {\cal C}(\mathbb{R} _+,
\mathbb{R}),$ $ (X_t)_{t \geq 0},\; ({\cal F}_t)_{t \ge 0}\big)$.
To   an $({\cal F}_t)$-adapted, non negative process $( F_t)_{t
\ge 0}$ such that $0 < E_x (F_t) < \infty$, for any $t \ge 0, \; x
\in \mathbb{R}$, we associate   the probability measure
$Q_{x,t}^{F}$ defined on $(\Omega, \; {\cal F}_t)$ as follows :
\begin{equation}\label{int1}
    Q_{x,t}^{F} (\Gamma_t) = \frac{1 }{ E_x [F_t]} \; E_x[1_{\Gamma_t} F_t],
\quad \Gamma_t \in {\cal F}_t .
\end{equation}
\noi A priori, the family $(Q_{x,t}^{F} \;;\; t \ge 0)$ is not
consistent : $Q_{x,t}^{F} (\Gamma_s)$ may be different from
$Q_{x,s}^{F} (\Gamma_s)$ for $\Gamma_s \in {\cal F}_s$ and $s <
t$; in fact, it is easy to see that $(Q_{x,t}^{F} ; t \ge 0)$ is
consistent, if and only if $\widetilde{F}_t = F_t / E_x [F_t]\;;\;
t \ge 0$ is a $P_{x}$-martingale. When this condition holds, we
write $Q_x^{\widetilde{F}}$ instead of $Q_{x,t}^{F}$.
$Q_{x}^{\widetilde{F}}$ is well defined since $Q_x^{\widetilde{F}}
(\Gamma_s) = E_x [1_{\Gamma_s} \widetilde{F}_t], \;\Gamma_s \in
{\cal F}_s$ and $s \le t$.

\noi In a previous study (\cite{RoyValYor}, \cite{RVY2}), we have
considered $F_t = \exp\Big\{ - \dis \frac{1}{ 2} \int_0^t
V(X_s)ds\Big\}$, where $V : \mathbb{R} \longrightarrow
\mathbb{R}_+$ is a Borel function. Our basic result was the
following : under some suitable assumptions on $V$, for any given
$s\geq 0$ and $\Gamma_s$ in ${\cal F}_s, \; Q_{x,t}^{F}
(\Gamma_s)$ converges as $t \to \infty$, to $Q_x^{\widetilde{F}}
(\Gamma_s)$ where $\widetilde{F}$ is the $(P_x)$-martingale :
\begin{equation}\label{int2}
\widetilde{F}_t = \frac{\varphi_V (X_t)}{\varphi_V (X_0)} F_t =
\frac{\varphi_V (X_t)}{\varphi_V (X_0)} \; \exp \Big\{- \frac{1 }{
2} \int_0^t V(X_s) ds\Big\},
\end{equation}

\noi and $\varphi_V$ is a "good" positive solution of the
Sturm-Liouville equation $\varphi'' =V \varphi$. The weak
convergence of $Q_{x,t}^{F}$ to $Q_{x}^{\widetilde{F}}, \; t \to
\infty$, is a direct consequence of the two following facts :
\begin{equation}\label{int3}
\big((X_t, \; F_t)\;;\;t \ge 0\big) \; \mbox{is a Markov process},
\end{equation}
\begin{equation}\label{int4}
    E_x [F_t] = E_x \Big[\exp\Big\{ - \frac{1}{2}
 \int_0^t V (X_s)ds \Big\}\Big] _{\stackrel{  \sim}{t \to \infty}}
t^{-k} \varphi_V(x), \; \mbox{for some}\; k \ge 0.
\end{equation}

\smallskip

\noi {\bf 1.2} The goal of this paper is to deal with a more
general setting by  considering $F_t = f(X_t, A_t)$, where $f :
\mathbb{R} \times \mathbb{R}^d \longrightarrow ]0, + \infty[$ is a
Borel function, $(A_t \;;\; t \ge 0)$ is $({\cal F}_t)$ adapted
and  $\mathbb{R}^d$- valued. We suppose moreover :
\begin{equation}\label{int5}
\big(Y_{t} {\stackrel {{\rm (def)}}{=}}(X_t, A_t) \;;\; t \ge
0\big)\; \mbox{is a } \ \Big((P_x)_{x \in \mathbb{R}} \;;\; ({\cal
F}_t)_{t \ge 0}\big)\mbox{-Markov process}.
\end{equation}
\noi Let $(\Lambda_t)_{t \ge
0}=\Big(\Lambda_t(y,dy')=\Lambda_t(x,a;dx'da') ; t \ge 0\Big)$ be
its semigroup (we denote $y=(x,a)$ and $y'=(x',a')$).

\noi To recover the setting of \cite{RoyValYor}, \cite{RVY2}
recalled above, it suffices to choose $d=1$, $f(x,a)=e^{- a/2}$
and $A_t= \dis \int_0^t V (X_s)ds$. Let $(Q_{x,t}^{F};\  x \in
\mathbb{R}, t \ge 0)$ be the family of p.m. ($\equiv$ probability
measures) associated with $(F_t)_{t \ge 0}$ :
\begin{equation}\label{int6}
    Q_{x,t}^{F} (\Gamma_t) = \frac{1 }{ E_x[F_t]} E_x[1_{\Gamma_t} F_t]
=\frac{1 }{ E_x\big[f(X_t, A_t)\big]} E_x\big[1_{\Gamma_t}
f(X_t,A_t)\big] , \; \Gamma_t \in {\cal F}_t, \; t \ge 0.
\end{equation}

\noi We now present a "meta-theorem", i.e. a statement which will
hold in great generality, so much so that our remaining study
shall consist in verifying that the hypotheses of Theorem
\ref{tint1} hold in various cases.

\begin{theo}  \label{tint1} Let $y_0=(x_0,a_0)$. We suppose $P_{x_{0}}(A_0=a_0)=1$ and
\begin{equation}\label{int7}
   M_s(y_0, f ;y):=  \lim_{t \to \infty}
\frac{\Lambda_{t-s} (f) (y) }{ \Lambda_t (f) (y_0)}
\end{equation}
\noi exists, for any  $s \ge 0 $ and $y=(x,a) \in \mathbb{R}
\times \mathbb{R}^d$,
\begin{equation}\label{int8}
    \frac{\Lambda_{t-s} (f)(y)}{\Lambda_t (f) (y_0)} \le C(s,y)
\;;\; \forall t >s,
\end{equation}
\noi where $E_{x_{0}} \big[C(s, Y_s)\big] < \infty$.

\noi  Then :
\begin{enumerate}
    \item $( M_s:=M_s(y_0, f ;Y_s);s\geq 0)$ is a non-negative
$P_{x_{0}}$-martingale, and $M_0=1, P_{x_0}$-a.s.

    \item $Q_{x_{0},t}^F$  converges weakly to $Q_{x_{0}}^M$, $t \to \infty$
(i.e. $ \displaystyle \lim_{t \to \infty} Q_{x_{0},t}^F (\Gamma_s)
= Q_{x_{0}}^M (\Gamma_s), \ \forall \Gamma_s \in {\cal F}_s, \;
\forall s
> 0)$, where :
\begin{equation}\label{int9}
    Q_{x_{0}}^M (\Gamma_s) =  E_{x_{0}} [1_{\Gamma_s} M_s]
, \; \Gamma_s \in {\cal F}_s \;;\; s \ge 0.
\end{equation}
\end{enumerate}
\end{theo}

\noi By definition, the p.m. $Q_{x_0,t}^F$  is absolutely
continuous on $(\Omega ,{\cal F}_t)$, with respect to the Wiener
measure. The Radon-Nikodym density $F_t/E[F_t]$ may be interpreted
as a weight or a penalization as it is done in statistical
mechanics. In the sequel, Theorem \ref{tint1} will be refered to
as a penalization principle.

\begin{prooff} \ {\bf of Theorem \ref{tint1}} Let $s > 0$ and $\Gamma_s \in {\cal F}_s$
fixed. Using the definition of $Q_{x,t}^{F} (\Gamma_s)$ and the
Markov property we have :
$$Q_{x_{0},t}^{F} (\Gamma_s)=E_{x_{0}} \Big[1_{\Gamma_s} \frac{\Lambda_{t-s} (f)(X_s,A_s)
}{\Lambda_t (f) (x_0,A_0)}\Big] .$$

\noi Property (\ref{int7}) and inequality (\ref{int8}) allow us to
apply the dominated convergence theorem :
$$ \lim_{t \to \infty}Q_{x_{0},t}^F (\Gamma_s) = Q_{x_{0}}^M
(\Gamma_s), $$
\noi where $M_s=M_s(y_0, f ;Y_s)$.

\noi It is clear that $M_0(y_0, f ;Y_0)=1$; consequently
$P_{x_0}(M_0=1)=1$.

\noi Let $0 \le s < s'$ and $\Gamma_s \in {\cal F}_s$; since
$\Gamma_s \in {\cal F}_{s'}$ we have :  $ E_{x_{0}} [1_{\Gamma_s}
M_s] = E_{x_{0}} [1_{\Gamma _s}M_{s'}]$. This means that $(M_t)$
is a $P_{x_{0}}$-martingale.
\end{prooff}

\noi {\bf 1.3} In this paper, we investigate  four cases of
examples involving respectively for $(A_t)$ :
\begin{itemize}
    \item the unilateral maximum (resp.
minimum) $S_t$ (resp. $I_t)$ : $ \displaystyle S_t =\max_{0 \le u
\le t} X_u$ (resp. $\displaystyle I_t = - \min_{0 \le u \le t}
X_u)$. We also consider, in the same case study the
two-dimensional process $(S_t,t)$.
    \item $(L_t^0; t \ge 0)$ the local time at $0$ of $(X_t)_{t \ge
0}$.

\item The triplet $((S_t,I_t,L^0_t); t\geq 0)$.

    \item $(D_t; t \ge 0)$ the number of down-crossings of $X$ from level
$b$ to level $a$.
\end{itemize}

\noi We observe that, in all cases, the function $M_s(y_0,f;y)$
may be written as :
$$ \frac{M
(f;y)}{M (f;y_0)}e^{\alpha s} ,$$
for some function $M$ and some $\alpha \in \mathbb{R}$; in fact
$\alpha =0$, except for the case 1, b), as shown below.

 \noi Since
$s\rightarrow M_s(y_0,f;Y_s)$ is  a $P_{x_0}$-martingale, it is
clear that $s\rightarrow M(f;Y_s)e^{\alpha s}$ is also a
$P_{x_0}$-martingale.

 \noi The results are summarized in the following Table
:

\begin{center}
\begin{tabular}{| c | c | c | c |c|}
\hline Cases & $A_t$ & $F _t$ & $M(f;Y_t)e^{\alpha t}$&Theorem \\
\hline
& & & \\
1 & $ a) \qquad S_t$ & $\varphi (S_t)$ & $M_t^\varphi$& \ref{tosm1} \\
\cline{2-5}
& & & \\
 & $ b) \  (S_t ,t)$ & $\displaystyle \varphi (S_t)e^{\lambda (S_t
 -X_t)}$& $M_t^{\lambda ,\varphi}$ & \ref{tosm2}\\
\hline  & & & \\
2& $L^0_t$ & $h^+(L^0_t)1_{\{X_t>0\}}+h^-(L^0_t)1_{\{X_t<0\}}$ &
$M_t^{h^+,h^-}$ & \ref{tloc1}\\

\hline
& & & \\
 3 & $(S_t,I_t,L^0_t)$ & $A_\nu (S_t ,I_t,L^0 _t)$
&$ M_t
^\nu $ & \ref{tmil1}\\

\hline
& & & \\
4 & $D_t$ & $\Delta G(D_t)$ & $ M_t^{\downarrow , G}$ & \ref{tdo1}\\

\hline
\end{tabular}

\end{center}

\smallskip

\noi These four cases will be  treated in sections
\ref{osm}-\ref{do} respectively.

\noi Let us describe the martingale $M(f;Y_t)e^{\alpha t}$ for
each case.

\smallskip

\noi {\bf Case 1.} a) We have :
\begin{equation}\label{int10}
    M_t^{\varphi}=(S_t - X_t)\varphi (S_t) +1- \Phi (S_t),
\end{equation}
where $\varphi : \mathbb{R} \mapsto \mathbb{R}_+$ is bounded,
$\dis \int_{\mathbb{R}} \varphi (u)du=1$ and $\Phi(x)=\dis \int_{-
\infty}^x \varphi (u)du$.

b) More generally, for $\lambda >0$ :
\begin{equation}\label{int11}
M_t^{\lambda, \varphi} =\Big\{ (1-\Phi(S_t)) \cosh \big( \lambda
(S_t-X_t)\big) + \varphi (S_t)\frac { \sinh\big(\lambda
(S_t-X_t)\big) }{ \lambda}\Big\} e^{-\lambda ^2 t/2},
\end{equation}
\noi where $ \psi : \mathbb{R} \mapsto \mathbb{R}_+, \xi_\lambda
(x) = e^{\lambda x} 1_{\{x<0\}}, 1-\Phi =\psi*\xi_\lambda$ and
$\varphi=\Phi'$.

\noi Both families $(M_t^{\varphi}, M_t^{\lambda, \varphi})$ were
intensively used in  \cite{AY1} to solve Skorokhod's problem for
Brownian motion. The class $(M_t^{\varphi})$ (resp.
$(M_t^{\lambda, \varphi})$ is defined in detail in Proposition
\ref{posm1} (resp. Proposition \ref{posm2}) and the special case
of the related meta-theorem in Theorem \ref{tosm1} (resp. Theorem
\ref{tosm2}).

\smallskip

\noi {\bf Case 2.} Let $h^+,h^-:\mathbb{R}_+ \longrightarrow
\mathbb{R}_+$ be two Borel, bounded functions, $H(l):=\dis \frac{1
}{ 2} \int_0^l \big(h^+(u)+h^-(u)\big)du$. We suppose $H(+
\infty)=1$. $(M_t^{h^{+},h^{-}})$ is the martingale :
%
\begin{equation}\label{int12}
    M_t^{h^{+},h^{-}} = X_t^+ h^+ (L_t^0)+X_t^- h^- (L_t^0)+1-
H(L_t^0).
\end{equation}
These martingales appear in\cite{JY}.

 \noi The limit theorem
associated with $A_t=L_t^0$ is stated, only with $x_0=0$, in
Theorem \ref{tloc1}.

\smallskip

\noi {\bf Case 3.} Let $\nu$ be a p.m. on $[\alpha, \infty[\times
[\alpha , \infty[$ for some $\alpha > 0$ and :
$$A_{\nu} (s,i,l) = \int_{\mathbb{R}_{+}^{2}} e^{\frac{1}{2} \big(
{1 \over a}+{1 \over b}\big)l} 1_{\{a \ge s, b \ge i\}} \nu
(da,db) \;;\; s , i,l \ge 0.$$
\noi $(M_t^{\nu})$ is the martingale :
\begin{equation}\label{int13}
    M_t^{\nu} = \int_{\mathbb{R}_{+}^{2}}\Big(1-\frac{X_{t}^{+}} {
a}\Big) \Big(1-\frac{X_{t}^{-}} { b} \Big) \exp \Big\{\frac{1 }{
2} \Big(\frac{1}{ a}+\frac{1 }{ b}\Big) L_t^0\Big\} \; 1_{\{S_t
\le a, I_t \le b\}} \nu(da,db).
\end{equation}
\noi Some properties  of the family $(M_t^{\nu})$ are given in
Proposition \ref{pmil1}, and the penalization principle is stated
in Theorem \ref{tmil1}.

\smallskip

\noi {\bf Case 4.} Let $D_t$ be the number of down-crossings from
$b$ to $a$, achieved by $(X_t)$ up to time $t$, and
$\big(G(n)\big)_{n \ge 0}$ be a decreasing sequence of positive
numbers such that $G(0)=1$ and $ \dis \lim_{n \to \infty} G
(n)=0$. Then
\begin{equation}\label{int14}
 M_{t}^{\downarrow , G}= \sum_{n \ge 0} \Biggr\{
1_{[\sigma_{2n},\sigma_{2n+1}[}(t)\Biggr( \frac{G(n)}{ 2} \Big(
1+\frac{b-X_t }{ b-a}\Big) +\frac {G(n+1) }{ 2} \frac{X_t - a }{
b-a}\Biggr)
\end{equation}
$$ + 1_{[\sigma_{2n+1},\sigma_{2n+2}[}(t)\Biggr( \frac{G(n+1)}{ 2} \Big(
1+\frac{b-X_t }{ b-a}\Big) +\frac {G(n) }{ 2} \frac{X_t - a }{
b-a}\Biggr)\Biggr\}.$$
\noi where $\sigma_0=0$ and $(\sigma_n)_{n \ge 1}$ is defined
inductively as follows : $\sigma_1 = {\rm inf}\{t \ge 0 \;;\; X_t
>b\}, \; \sigma_2={\rm inf} \{t \ge 0, \; X_t < a\},\;
\sigma_{2n+1}=\sigma_1 \circ \theta_{\sigma_{2n}}, \;
\sigma_{2n+2} = \sigma_{2} \circ \theta_{\sigma_{2n+1}}$, where
$(\theta_u)_{u \geq 0}$ denotes the family of shift operators on
the canonical space.

\smallskip

\noi The corresponding case of the  meta-theorem is stated in
 Theorem \ref{tdo1}.


\noi {\bf 1.4} Theorem \ref{tint1} leads naturally to ask for a
description of the law of $(X_t)_{t \ge 0}$ under $Q_{x_{0}}^M$.
Since $(M_t)_{t \ge 0}$ is a strictly positive
$(P_{x_{0}})$-martingale, it may be written as an exponential
martingale :
$$ M_t={\cal E}(J)_t = \exp \Big\{
 \int_0^t J_s dX_s -\frac {1}{ 2} \int_0^t J_s^2 ds\Big\},$$
 \noi
for some adapted process $(J_t)_{t \ge 0}$.

\noi Girsanov's Theorem implies that $\dis (\beta_t ; t \ge
0)=\Big(X_t -  \int_0^t J_s \, ds ; t \ge 0 \Big)$ is a
$Q_{x_{0}}^M$-Brownian motion.

\noi Suppose that $(x,a) \longrightarrow M (f ;x,a)$ is of class
$C^{2,1}$ and $(A_t)$ has bounded variation. Then by Ito's formula
we obtain :
$$J_t = \alpha +\frac{\partial }{ \partial x} \Big( \log\big( M(
f; X_t, A_t)\big) \Big).$$

\noi Consequently $(X_t)$ solves the following stochastic
differential equation :
\begin{equation}\label{int15}
    X_t = \beta_t + \int_0^t \Big(\alpha +\frac{\frac{\partial M }{ \partial x}
}{ M} (f;X_s,A_s)\Big)ds.
\end{equation}

\noi We may recover (\ref{int15}) in a different manner. Suppose
for simplicity that $\alpha =0$ . It is clear that $\big(M(f;
Y_s); s \ge 0\big)$ is a non-negative $P_{x_{0}}$-martingale. In
other words, $y \to M(f;y)$ is a non-negative   harmonic function
with respect to $(Y_t)_{t \ge 0}$, under $P_{x_{0}}$. In
particular $(Y_t)_{t \ge 0}$ is a $Q_{x_{0}}^M$-Markov process
with semigroup :
$$\Lambda_t^M (g) (y) = \frac{\Lambda_t \big(g M (f;
\cdot)\big) (y) }{ M (f ;y)} .$$
 \noi If $(x,a) \mapsto  M (f ; x,a)$ is of class $C^{2,1}$, then
the generator associated with the semigroup $(\Lambda_t^M)$ is
$$ \frac{1 }{2} \frac{\partial^2 }{ \partial x^2} + \frac{\frac{\partial M }{
\partial x} }{ M} (f;x,a) \frac{\partial }{ \partial x}
+ \frac{\frac{\partial M }{
\partial a} }{ M} (f;x,a) \frac{\partial }{ \partial a}.$$
 This gives a new proof of (\ref{int15}).

\noi However in our four classes of examples we observe that the
drift term in (\ref{int15}) is explicit but complicated and
therefore does not allow to identify directly the law of $(X_t)$
under $Q_{x_0}^M$. Suppose that $(A_t)_{t \ge 0}$ is a
one-dimensional process, $Q_{x_{0}} (A_0=0)=1$ and $t \to A_t$ is
continuous and non-decreasing. The r.v. $A_\infty$, which may be
infinite, plays a central role in our approach. We claim that we
may compute the distribution function of $A_{\infty}$ under
$Q_{x_{0}}$. Let $(A_t^{-1})_{t \ge 0}$ be the right inverse of
$s\mapsto A_s$, i.e. : $A_t^{-1}={\rm inf} \{s \ge 0, \, A_s
> t\}$.  We have :
$$Q_{x_{0}} (A_t> \alpha)  = Q_{x_{0}} (A_{\alpha}^{-1}
<t)=E_{x_{0}} [1_{\{A_{\alpha}^{-1} <t\}}M_t]
  = E_{x_{0}}
[1_{\{A_{\alpha}^{-1} <t\}} M_{A_{\alpha}^{-1}}].
$$
\noi Note that the last equality  follows from  the optional
stopping theorem.

\noi Since $M_{A_{\alpha}^{-1}} \ge 0$, taking $t \to \infty$, in
the previous expression, we get :
$$
Q_{x_{0}} (A_{\infty} > \alpha)= E_{x_{0}}
\big[1_{\{A_{\alpha}^{-1}<\infty\}}M_{A_{\alpha}^{-1}}\big]=
E_{x_{0}}\big[1_{\{A_{\alpha}^{-1}< \infty\}} M(x_0,0, f ;
X_{A_\alpha^{-1}}, \alpha)\big].$$
\noi In cases 1 (with $A=S$ only), 2 and 4, we prove that
$Q_{x_{0}} (A_{\infty} < \infty)=1$, and more generally $Q_{x_{0}}
(0<g < \infty)=1$ where $g={\rm sup} \{s>0, \, A_s=A_{\infty}\}$
(with the convention $\sup \emptyset=0$). To describe the law of
$(X_t)$ under $Q_{x_0}$, it is convenient to split the whole
trajectory in two parts : $(X_t)_{0\leq t \leq g}$ and
$(X_{t+g})_{ t \geq 0}$. We observe that the random time $g$ is
not a $({\cal F}_t)$-stopping time but it is a last exit time.
Using the technique of enlargement of filtrations, we are able to
describe the law of $(X_t)_{0\leq t \leq g}$ and $(X_{t+g})_{ t
\geq 0}$ under $Q_{x_{0}}^M$, conditionally to $A_{\infty}$ (see
Theorems \ref{tQ2}, \ref{tQ3} and  \ref{tQ5}).

\noi As for case 3 \big(i.e. $A_t=(S_t, I_t, L_t^0)\big)$, we
prove that $Q_0(X_{\infty}^* < \infty)=1$, with $\dis
X_{\infty}^*=
 \sup_{t \ge 0} |X_t|=S_\infty \,\vee\, I_\infty$. However
$Q_0(L_{\infty}^0 =\infty)=1$. Conditionally on $X_{\infty}^*$,
the law of $(X_t)_{t \ge 0}$ is given in Theorem \ref{tQ4}, via a
path decomposition at time $g=\sup\{ t\geq 0; |X_t|=X_\infty
^*\}$.

\noi {\bf 1.5} In Section \ref{dir}, we recover previous results,
using a direct approach based on a disintegration of $Q_0$. We
exhibit a family $\big(Q_0^{(a)};a \in \mathbb{R}^d\big)$ of p.m.
on $\big(\Omega, {\cal F}_\infty\big)$ such that :
\begin{equation}\label{int16}
    Q_0(\cdot)=\int_{\mathbb{R}^d}Q_0^{(a)}(\cdot)\mu (da),
\end{equation}
\noi where in cases  1 (resp. 2, 3), $\mu$ is the law of
$S_\infty$ (resp. $L^0_\infty$, $(S_\infty,I_\infty,L^0_\infty)$).

\noi We are able to determine the law of $(X_t)$ under
$Q_0^{(a)}$, for any $a$. Moreover this distribution does not
depend on $f$.

\noi {\bf 1.6} Finally, in Section \ref{fur}, we present several
other directions of research, which we have now begun to
investigate and will be the subject of a future publication.


\section{ Notation}\label{no}

\setcounter{equation}{0}

 \noi In the sequel of the paper we shall
use intensively the following notation and conventions.

\begin{itemize}
    \item $\big(\Omega = {\cal C} (\mathbb{R}_+, \mathbb{R}),
\, (X_t)_{t \ge 0},\, ({\cal F}_t)_{t \ge 0}\big)$ is the
canonical space, with $(X_t)_{t \ge 0}$ the coordinate maps :
$X_t(\omega)=\omega(t)$ and ${\cal F}_t=\sigma\{X_s; 0\leq s \leq
t\}\ ; t \ge 0$.

    \item $(P_x)_{x \in \mathbb{R}}$ is the family of Wiener
measures on the canonical space $\Omega$ : under $P_x, \, (X_t)_{t
\ge 0}$ is a one-dimensional Brownian motion started at $x$.

\noi Let $\lambda$ and $x$ be two real numbers. We denote by
$P_x^{(\lambda )}$  the p.m. on the canonical space, under which
$(X_t)_{t \ge 0}$ is a Brownian motion with drift $\lambda$,
started at $x$. Obviously $P_x=P_x^{(0 )}$ and $(X_t)_{t \ge 0}$
is distributed under $P_x^{(\lambda )}$  as $(X_t+\lambda t)_{t
\ge 0}$ under $P_x$.

 \noi  If $Q$ is a probability measure (p.m.) on $\Omega$,
the expectation with respect to $Q$ is denoted $E_Q$. However if
$Q=P_x$,  we shall write $E_x$ for $E_{P_x}$ for simplicity.

    \item $(\theta_t)_{t \ge 0}$ is the family of shift operators from
    $\Omega$ to $\Omega$, defined by $\theta_t (\omega)(s)=
\omega (t+s)\;;\; s\ge 0, \ t \ge 0$.

    \item For any $a \in \mathbb{R},\, T_a$ is the first
hitting time of level $a$, namely : $T_a = {\inf} \{t \ge 0, \,
X_t = a\}$. We adopt the convention $\inf\{\emptyset\} =+ \infty$.

    \item $(S_t ; t\ge 0)$  (resp. $(-I_t ; t \ge
0)$), is the one-sided maximum (resp. minimum) :
$$
S_t = \sup_{0 \le u \le t} X_u \;;\; I_t = - \inf_{0 \le u \le
t}X_u .
$$
    \item The bilateral maximum $(X_t^*)_{t \ge 0}$ is the process :
$$X_t^* = S_t \vee I_t =  \sup_{0 \le u \le t} |X_u|.$$
    \item $(L_t^x ; x \in \mathbb{R},t \ge 0)$ is the
jointly continuous family of local times associated with $(X_t)_{t
\ge 0}$. For simplicity we write $(L_t)_{t \ge 0}$ instead of
$(L_t^0)_{t \ge 0}$.

\item If $J$ is a predictable process (with respect to $({\cal
F}_t))$ such that $\dis \int_0^t J_s^2 ds < \infty\;P_x$ a.s. for
every $t$, then we denote by ${\cal E} (J)_t$ the $P_x$-local
martingale :
$$\exp \Big\{ \int_0^t J_s d X_s - \frac{1 }{ 2} \int_0^t J_s^2
ds\Big\}, \; t \ge 0.$$
\end{itemize}


\section{Penalization and associated martingales}\label{pen}

\subsection{ Case 1 : The one-sided maximum} \label{osm}


\noi In this section, we consider two families of local
martingales $(M_t^{\varphi})$ and $(M_{t}^{\lambda, \varphi})$
involving the one-sided maximum. These local martingales are
well-known : they play a prominent role in Az\'{e}ma-Yor's
solution of Skorokhod's problem  studied in \cite{AY1}. We will
show (see Theorems \ref{tosm1} and \ref{tosm2}) that a sub-class
of the previous local martingales appears naturally after a
penalization procedure.

\subsubsection{ The local martingales associated with the one-sided
maximum}\label{sosm1}

 \setcounter{equation}{0}
  \noi Let us start with local
martingales of the form $H(X_t, S_t)$, and more generally $H(X_t,
S_t,t)$, for some function $H$.

\begin{prop}\label{posm1}\begin{enumerate}
    \item  Let $\varphi : \mathbb{R} \longrightarrow ]0,+\infty[$
be a Borel function such that

\noi $\dis  \int_{- \infty}^\cdot \varphi (u)du<\infty$, and
define $\dis \Phi(s)= \int_{- \infty}^s \varphi (u)du$. Then :
\begin{equation}\label{osm1}
    M_t^{\varphi} := (S_t - X_t) \varphi (S_t)+1-\Phi (S_t),
\end{equation}
 is a $P_x$-martingale, $M_0^{\varphi} =1-\Phi(x), \ P_x$ a.s., and :
 \begin{equation}\label{osm2}
    M_t^{\varphi}=1- \Phi (x) - \dis \int_0^t \varphi (S_u)d
X_u.
\end{equation}
    \item Suppose moreover :
    \begin{equation}\label{osm3}
    \int_{\mathbb{R}} \varphi (u)du=1.
\end{equation}
\noi Then under $P_x, \, M_t^{\varphi} > 0$ and $\dis
M_t^{\varphi}=
 ( 1-\Phi (x))  {\cal E} (J^{\varphi})_t$ where :
 \begin{equation}\label{osm4}
    J_t^{\varphi} = \dis - {\varphi (S_t) \over M_t^{\varphi}}
=- {\varphi (S_t) \over (S_t - X_t) \varphi(S_t)+1-\Phi (S_t)}.
\end{equation}
\end{enumerate}

\end{prop}

\begin{rem}\label{rosm1}
\begin{enumerate}
    \item  Note that the local martingale $(M_t^{\varphi})$ is
actually a martingale since :
$$ \sup_{0 \le u \le t}
 |M_u ^{\varphi}|\le 1+ 2 ||\varphi||_{\infty} \sup_{0 \le u \le t}
|X_u|.$$
    \item  Taking  $\dis \varphi (u) = \frac{1 }{ a} 1_{[0,a]}
(u)$, with $a > 0$, we obtain : $\dis \Phi (x)= \frac{1 }{ a}
\big(x \wedge a\big)$ and $\dis M_t^{\varphi} = 1-\frac{X_{t
\wedge T_{a}}}{a}$.
\end{enumerate}
\end{rem}

\noi We now recall the definition of Kennedy's martingales, which
also played some role in the computation of the laws of the
stopping times studied in   \cite{AY1} .

\begin{prop}\label{posm2} Let $\lambda > 0$,  $\varphi : \mathbb{R} \longrightarrow \mathbb{R}$
be a locally integrable function, $\Phi$ be any primitive of
$\varphi$ ($\varphi(x)=\Phi'(x)$). Let $( M^{\lambda ,\varphi}_t)$
be the process :
\begin{equation}\label{osm5}
     M_t^{\lambda, \varphi} := \Big\{ (1-\Phi (S_t))\cosh
 \big(\lambda (S_t - X_t)\big) + \varphi (S_t)
\frac{\sinh \big(\lambda (S_t-X_t)\big) }{ \lambda} \Big\}
  e^{- \lambda^{2}t/2},
\end{equation}
\noi Then :
\begin{enumerate}
    \item   $( M^{\lambda ,\varphi}_t)$ is a $P_x$-local  martingale.
    Under $P_x,\;  M_0^{\lambda, \varphi}=1-\Phi(x)$ and :
    \begin{equation}\label{osm6}
     M_t^{\lambda, \varphi}=1-\Phi(x)-\int_0^t \Big\{-\lambda
\big(1-\Phi(S_u)\big) \sinh \big(\lambda (S_u-X_u)\big)+\varphi
(S_u ) \cosh\big(\lambda (S_u-X_u)\big)\Big\}  e^{-
\lambda^{2}u/2}dX_u.
\end{equation}

\item Let $x_0 \in\mathbb{R}$.  Then $ M_t^{\lambda, \varphi} \geq
0, \ \forall t\geq 0$ , $P_{x_0}$ a.s. if and only if there exists
a Borel function $\psi : \mathbb{R} \mapsto [0,\infty[$ and a
non-negative constant $\kappa$ such that the two following
conditions hold :
\begin{equation}\label{osm6a}
\int_{x_0}^\infty \psi(z)e^{-\lambda z}dz <\infty,
\end{equation}
\begin{equation}\label{osm6b}
   1- \Phi(y)=e^{\lambda y}\Big(\kappa +\int _y^\infty \psi(z)e^{-\lambda
    z}dz\Big), \ y \geq x_0.
\end{equation}
Assuming that (\ref{osm6a}) and (\ref{osm6b}) hold, then $
(M_t^{\lambda, \varphi})$ is a $P_{x_0}$-martingale.

    \item Let $\psi :\mathbb{R} \mapsto [0,\infty[$ be a
    Borel function satisfying :
\begin{equation}\label{osm6b1}
\int_{x}^\infty \psi(z)e^{-\lambda z}dz <\infty, \quad \forall
x\in \mathbb{R}.
\end{equation}
Let $\Phi : \mathbb{R} \mapsto \mathbb{R}$ be the function :
\begin{equation}\label{osm6b2}
\Phi(y)=1-e^{\lambda y}\int _y^\infty \psi(z)e^{-\lambda
    z}dz, \ y \in\mathbb{R}.
\end{equation}
Then :
\begin{equation}\label{osm6b3}
    \varphi(y)=\Phi'(y)=\psi(y)-\lambda e^{\lambda y}\int _y^\infty \psi(z)e^{-\lambda
    z}dz,
\end{equation}
$(M_t^{\lambda,\varphi})$ is a $P_x$-martingale and :
\begin{equation}\label{osm6b4}
  M_t^{\lambda,\varphi}= \Big\{  \psi (S_t) \frac{\sinh \big(\lambda (S_t-X_t)\big)
}{ \lambda}+ e^{\lambda X_t}\int_{S_t}^\infty \psi (z)e^{-\lambda
z} dz\Big\}
  e^{- \lambda^{2}t/2}.
\end{equation}
\noi Suppose moreover $\psi >0$, then under $P_x$ : $M_t^{\lambda,
\varphi}>0$ and $M_t^{\lambda, \varphi}= (1-\Phi(x)){\cal E}
(J^{\lambda, \varphi})_t$, where :
\begin{equation}\label{osm6b5}
    J_t^{\lambda, \varphi} = - \lambda \frac{\varphi (S_t) \cosh
\big(\lambda (S_t-X_t)\big)+\lambda \big(1- \Phi (S_t)\big) \sinh
\big(\lambda (S_t-X_t)\big) }{\lambda \big(1- \Phi (S_t)\big)
\cosh \big(\lambda (S_t-X_t)\big) + \varphi (S_t)\sinh
\big(\lambda (S_t-X_t)\big)} .
\end{equation}
\end{enumerate}
\end{prop}

\begin{prooff} \ {\bf of Proposition \ref{posm2}} 1) If $\varphi$
is  of class $C^1$, then point 1) of Proposition \ref{posm2} is a
direct consequence of It\^{o}'s formula and the fact that $dS$ is
carried by $\{S-X=0\}$. The general case follows from the monotone
class theorem.

\noi 2) Let us investigate the positivity of $ M_t^{\lambda,
\varphi}$, under $P_{x_0}$. It is clear that : $ M_t^{\lambda,
\varphi} \geq 0, \ \forall t\geq 0$ , $P_{x_0}$ a.s. is equivalent
to :
$$(1-\Phi (y))\cosh
 \big(\lambda (y - x)\big) + \varphi (y)
\frac{\sinh \big(\lambda (y-x) \big) }{ \lambda} \geq 0,
$$
for all $x$ and $y$ such that $y\geq x_0$ and $x\leq y$.

\noi Setting $u=y-x$ and rewriting $\sinh$ and $\cosh$ in terms of
exponential functions, it is easy to check that the previous
inequality is equivalent to :
$$e^{2\lambda u}\big((1-\Phi (y)) +\frac{\varphi (y)}{\lambda}\big)
+1-\Phi (y) -\frac{\varphi (y)}{\lambda} \geq 0, \quad \forall
u\geq 0, y\geq x_0.$$
\noi Since $\alpha Y_1+Y_2\geq 0, \ \forall \alpha \geq 1$ iff
$Y_1\geq 0$ and $Y_1+Y_2\geq 0$, the previous inequality is
 equivalent to :
\begin{equation}\label{osm6c}
    1-\Phi (y) +\frac{\varphi (y)}{\lambda} \geq 0, \quad \forall
y\geq x_0,
\end{equation}
\noi and
\begin{equation}\label{osm6d}
    \Phi (y)\leq 1, \quad \forall
y\geq x_0 .
\end{equation}
\noi Consequently the function $\psi$ defined by
\begin{equation}\label{osm6e}
    \psi (y)=\lambda (1-\Phi (y)) +\varphi (y), y\geq x_0,
\end{equation}
\noi takes its values in $[0,\infty[$.

\noi Recall that $\Phi'=\varphi$, hence (\ref{osm6e}) may be
interpreted as an ordinary linear  differential equation in
$\Phi$, which is easily solved :
\begin{equation}\label{osm6f}
\Phi (y)=1+e^{\lambda y}\Big( \kappa_0+\int _{x_0}^y\psi
(z)e^{-\lambda z} dz\Big), \ y \geq x_0,
\end{equation}
\noi  where $\kappa_0$ is a constant.

\noi Since (\ref{osm6c}) and (\ref{osm6f}) are equivalent, it
remains to deal with (\ref{osm6d}). Obviously this inequality is
equivalent to :
$$ \kappa_0+\int _{x_0}^y\psi
(z)e^{-\lambda z} dz \leq 0, \ \forall  \ y \geq x_0.
$$
\noi The function $\psi$ being non-negative, this last condition
is equivalent to (\ref{osm6a}) and
$$ \kappa=-\kappa_0-\int _{x_0}^\infty \psi
(z)e^{-\lambda z} dz \geq 0.$$
\noi  Relation (\ref{osm6b}) is a direct consequence of
(\ref{osm6f}).

 \noi 3) Choosing $\kappa =0$ in(\ref{osm6b}), we easily
obtain (\ref{osm6b3})-(\ref{osm6b5}).

\noi 4) Suppose that $\Phi$ is given by (\ref{osm6b}) and $\psi$
verifies (\ref{osm6a}). We know that $(M^{\lambda, \varphi}_t)$ is
a non-negative local martingale. We would like to prove that it is
in fact a  $P_x$-martingale.

 \noi a) Let $\psi _n(y)=\psi (y)1_{\{y\leq
n\}}, \varphi _n=\Phi'_n$, with $\dis \Phi _n(y)=1-e^{\lambda
y}\int_y^\infty \psi _n(z)e^{-\lambda z}dz$.

\noi Since for any $y\geq x$, we have :
$$ \int_y^\infty \psi_n(z)e^{-\lambda z}dz \leq
\int_y^\infty \psi (z)e^{-\lambda z}dz \leq \int_x^\infty \psi
(z)e^{-\lambda z}dz,$$
then
\begin{itemize}
    \item $\dis e^{\lambda X_t}\int_{S_t}^\infty \psi (z)e^{-\lambda
    z}dz$ and $\dis e^{\lambda X_t}\int_{S_t}^\infty \psi_n (z)e^{-\lambda
    z}dz$ are $P_x$-integrable r.v.'s,

    \item $\dis e^{\lambda X_t}\int_{S_t}^\infty \psi_n (z)e^{-\lambda
    z}dz$ goes to $\dis e^{\lambda X_t}\int_{S_t}^\infty \psi (z)e^{-\lambda
    z}dz$, in $L^1(\Omega,P_x)$, as $n\rightarrow\infty$.
\end{itemize}

\noi b) For $u\geq 0$, we have $\sinh (u)\leq e^u$, then
inequality (\ref{osm15a}) (which will be proved independently
later) implies that $\dis E_x[\psi
(S_t)e^{\lambda(S_t-X_t)}]<\infty$ and
$$\lim_{n\rightarrow\infty}|E_x[\psi
(S_t)e^{\lambda(S_t-X_t)}]-E_x[\psi _n(S_t)e^{\lambda(S_t-X_t)}]|
=\lim_{n\rightarrow\infty} E_x[\psi
(S_t)1_{\{S_t>n\}}e^{\lambda(S_t-X_t)}]=0.$$
\noi c) Consequently the representation (\ref{osm6b4}) implies
that $E_x[M^{\lambda, \varphi}_t]<\infty$ and $M^{\lambda,
\varphi_n}_t$ goes to $M^{\lambda, \varphi}_t$, in  $L^1(\Omega
,P_x)$, as $n\rightarrow\infty$.

 \noi It is clear that $(M^{\lambda, \varphi_n}_t)$ is a
$P_x$-martingale. As a result, $(M^{\lambda, \varphi}_t)$ is a
$P_x$-martingale.

\end{prooff}

\begin{rem}\label{rosm1b}
\noi Recall that   under $P^{(-\lambda)}_x$, $(X_t)$ is a Brownian
motion with drift $-\lambda$, started at $x$.

\noi Let $(\widetilde{M}_t)$ be the process :
\begin{equation}\label{osm6g}
    \widetilde{M}_t:=M^{\lambda, \varphi}_t e^{\big\{\lambda X_t+\lambda
^2 t/2\big\}}= \frac{\psi (S_t)}{2\lambda}\big( e^{\lambda S_t} -
e^{\{\lambda (2X_t-S_t)\}}\big) +e^{2\lambda X_t}\int_{S_t}^\infty
\psi (z)e^{-\lambda z} dz, \ t\geq 0.
\end{equation}
\noi It is clear that  $(\widetilde{M}_t)$ is a $P^{(-\lambda)}_x$
martingale.

 \noi But $\dis
P^{(-\lambda)}_x\big(\lim_{t\rightarrow
\infty}X_t=-\infty\big)=1$, therefore, under $P^{(-\lambda)}_x$  :
$$\lim_{t\rightarrow \infty}\widetilde{M}_t=\widetilde{M}_\infty :=\frac{\psi (S_\infty)}
{2\lambda} e^{\lambda S_\infty}.
$$
\noi Recall that (see e.g. Williams \cite{Wi}; but this result
also follows from Theorem \ref{tQ1b}, 1.):
\begin{equation}\label{osm6i}
    P^{(-\lambda)}_x(S_\infty >y)=P^{(-\lambda)}_0(S_\infty
+x>y)=e^{-2\lambda(y-x)}, \ y\geq x.
\end{equation}
This directly implies :
$$
E^{(-\lambda)}_x\big[\frac{\psi (S_\infty)} {2\lambda} e^{\lambda
S_\infty}\big]=e^{2\lambda x}\int_x^\infty \psi(y) e^{-\lambda
y}dy <\infty .$$
Finally, $(\widetilde{M}_t)$ is a non-negative $P^{(-\lambda)}_x$
martingale, converging a.s. to $\widetilde{M}_\infty \in
L^1(\Omega)$ and
$$ E^{(-\lambda)}_x[\widetilde{M}_0]=(1-\Phi(x))e^{\lambda x}=
e^{2\lambda x}\int_x^\infty \psi(y) e^{-\lambda y}dy
 = E^{(-\lambda)}_x[\widetilde{M}_\infty].$$

\noi As a result,
 $(\widetilde{M}_t)$ is a uniformly integrable $P^{(-\lambda)}_x$
martingale, and :
\begin{equation}\label{osm6h}
\widetilde{M}_t=E^{(-\lambda)}_x\Big[\frac{\psi (S_\infty)}
{2\lambda} e^{\lambda S_\infty}| {\cal F}_t\Big], \ t \geq 0.
\end{equation}
\noi This easily implies that $(M^{\lambda , \varphi}_t)$ is a
$P_x$-martingale. Hence the arguments developed in this remark may
be used instead of those in point 4) of the Proof of Proposition
\ref{posm2}.
\end{rem}


\begin{rem} \label{rosm2}\begin{enumerate}
    \item Suppose that $\psi$ verifies the conditions given in 3. of Proposition
     \ref{posm2} and

     \noi $\dis \int _\mathbb{R}\psi (z)dz=1$. We observe
      that $M_t^{\lambda, \varphi}$, as defined in
(\ref{osm5}), converges as $\lambda \rightarrow  0$, to
$M_t^{\varphi}$. This leads us to adopt the convention $M^{0,
\varphi}_t=M^{\varphi}_t$.
    \item If $\varphi=0$, then $M_t^{\lambda, 0} =
    \cosh \big(\lambda (S_t-X_t)\big) \, e^{- \lambda^{2}t/2}$.
\end{enumerate}
\end{rem}

\subsubsection{ Penalization involving the unilateral maximum}
\label{sosm2}

\noi As in sub-section \ref{osm1}, $\varphi : \mathbb{R} \mapsto
]0, \infty[$ is a Borel function. We suppose moreover that
(\ref{osm3}) holds. Hence $\varphi$ is actually a probability
density
 function. We denote $\dis \Phi (x)=  \int_{- \infty}^x
\varphi(y)dy$. As we shall see, the martingales involved in the
penalization result stated below, belong to the family  $(M_t^{
\varphi})$ as defined in Proposition \ref{posm1}.

\begin{theo} \label{tosm1}
Let $\varphi$ be as above.
\begin{enumerate} \item  Let $u \ge 0$ and $x \in \mathbb{R}$. For any
$\Gamma_u$ in ${\cal F}_u$, we have :
\begin{equation}\label{osm7}
    \lim_{t \to \infty} \; \frac{E_x\big[1_{\Gamma_u} \varphi
(S_t) \big] }{ E_x \big[\varphi (S_t)\big]} = \frac{1 }{1-\Phi(x)}
 E_x [1_{\Gamma_u} M_u^{\varphi}],
\end{equation}
\noi where $(M_u^{\varphi})_{u \ge 0}$ is the martingale defined
in (\ref{osm1}).

    \item Let $(Q_x^{\varphi})_{x \in \mathbb{R}}$ be the
family of probabilities on $\big(\Omega, {\cal F}_\infty\big)$ :
\begin{equation}\label{osm8}
    Q_x^{\varphi} (\Gamma_u) = \frac{1 }{ 1-\Phi(x)} \, E_x [1_{\Gamma_u}
\, M_u^{\varphi}], \quad \mbox{for any} \ u \ge 0,  \ \mbox{and} \
\Gamma_u \in {\cal F}_u.
\end{equation}
\noi Then, under $Q_x^{\varphi}$, the process $\dis
\Big(X_t-x+\int_0^t \frac{\varphi (S_u) }{ M_u^{\varphi}}du ; t
\ge 0 \Big)$ is a Brownian motion, started at $0$.

\end{enumerate}

\end{theo}

\begin{rem}\label{rosm3}
 $(M_t^{\varphi} \;;\; t \ge 0)$ is a
$P_x$-martingale but it is not uniformly integrable (u.i.) :
indeed,  $(M_t^{\varphi}; t \ge 0)$ is a positive martingale, thus
it converges a.s. to $M_{\infty}^{\varphi} \ge 0$, as $t \to
\infty$. Let $(t_n)_{n \ge 1}$ be an increasing sequence of times
such that $X_{t_n}=S_{t_n}$, e.g. :  $t_n=\inf \{u \ge 0, \; X_u
=n\}$ is convenient. Since $M_{t_n}^{\varphi}=1-\Phi (S_{t_n})$,
under $P_x$, $M_{t_n}^{\varphi}$ goes a.s. to $0$, as $n \to
\infty$. Hence, $M_{\infty}^{\varphi}=0$.
\end{rem}

\noi The proof of Theorem \ref{tosm1} requires a preliminary
result.

\begin{lemma} \label{losm1} Let $x \le a$ and
$\varphi_0 : [a,+\infty[ \mapsto \mathbb{R}_+$ such that $\dis
\int_{a} ^\infty \varphi_0 (u) du <+ \infty$. Then :
\begin{equation}\label{osm9}
    E_0 \Big[\varphi_0 \big(a \vee (x+S_u)\big)\Big]
_{\stackrel{  \sim} {u \rightarrow \infty}}
 \sqrt{\frac{2 }{ \pi u}} \Big\{(a-x) \varphi_0(a) +
\int_a^{\infty} \varphi_0(y)dy\Big\},
\end{equation}
\begin{equation}\label{osm9a}
E_0 \Big[\varphi_0 \big(a \vee (x+S_u)\big)\Big] \leq
 \sqrt{\frac{2 }{ \pi u}} \Big\{(a-x) \varphi_0(a) +
\int_a^{\infty} \varphi_0(y)dy\Big\},
\end{equation}
\begin{equation}\label{osm9b}
    E_0 \Big[\varphi_0 \big(a \vee (x+S_u)\big)\Big] \geq
    \sqrt{\frac{2 }{ \pi u}}\int_0^{\infty} \varphi_0 (y) \,
e^{-(y-x)^{2}/2} dy, \quad \mbox{ for any} \ u \geq 1.
\end{equation}
\end{lemma}

\begin{prooff} \ {\bf of Lemma \ref{losm1}}.  We have :
$$ E_0\Big[ \varphi_0 \big(a \vee (x+S_u)\big)\Big]=\varphi_0(a)
P_0(x+S_u \le a) + E_0\big[\varphi_0(x+S_u) 1_{\{x+S_u \ge a\}}
\big] $$

\noi Recall that under $P_0$,
\begin{equation}\label{osm10}
    S_u  {\stackrel {(d)}{=}} |X_u| {\stackrel {(d)}{=}} \sqrt{u} |X_1|.
\end{equation}
\noi Consequently,
$$
E_0\Big[ \varphi_0 \big(a \vee (x+S_u)\big)\Big]=\varphi_0(a)
\sqrt{\frac{2 }{ \pi }}\int_0^{(a-x)/\sqrt {u}}e^{-z^2/2}dz+
\sqrt{\frac{2 }{ \pi u }}\int_a^\infty \varphi_0 (y)
e^{-(y-x)^{2}/2u} dy.$$

\noi Then (\ref{osm9}) and (\ref{osm9a}) follow immediately.

\noi As for (\ref{osm9b}), we have :
$$E_0\Big[ \varphi \big(a \vee (x+S_u)\big)\Big] \geq
\sqrt{\frac{2 }{ \pi u }} \int_a^\infty \varphi_0 (y)
e^{-(y-x)^{2}/2u} dy
 \geq \sqrt{\frac{2 }{ \pi u }} \int_a^\infty \varphi_0 (y)
e^{-(y-x)^{2}/2} dy,
$$
if $u\geq 1$.
\end{prooff}

\begin{prooff} \ {\bf of Theorem \ref{tosm1}}

\noi  1) Let $s \ge 0, \, x \in \mathbb{R}$ and $\Gamma _s \in
{\cal F}_s$ be fixed. We consider $t>s$. Since $\dis S_t=S_s \vee
\big\{X_s +
 \sup_{0 \le u \le t-s}(X_{u+s}-X_s) \big\}$, and
under $P_x$, $(X_{u+s} - X_s ; u \ge 0)$ is a Brownian motion
started at $0$,
$$E_x \big[\varphi (S_t)|{\cal F}_s\big]= \widetilde{\varphi}
(S_s,X_s;t-s),$$

\noi where
$$\widetilde{\varphi} (a,x;r)=E_0 \Big[\varphi \big( a \vee
(x+S_r)\big)\Big].$$

\noi Lemma \ref{losm1} implies :
$$E_x\big[\varphi (S_t)|{\cal F}_s\big]  _{\stackrel {\sim}{t \rightarrow
\infty}} \sqrt{\frac{2 }{ \pi (t-s)} } M_s^{\varphi}.$$
\noi Using Lemma  \ref{losm1} with $x=a$, we obtain :
$$E_x\big[\varphi(S_t)\big]=E_0\big[\varphi(x+S_t)\big]
_{\stackrel {\sim}{t \rightarrow \infty}} \sqrt{\frac{2 }{ \pi t}}
\big(1-\Phi(x)\big).$$

\noi It is now easy to check (\ref{osm7}), using the two previous
estimates, (\ref{osm9a}), (\ref{osm9b}) and :
$$\frac {E_x \big[1_{\Gamma_s} \varphi (S_t)\big] }{ E_x \big[\varphi
(S_t)\big]} = E_x \Big[1_{\Gamma_s}
 \frac{ E_x\big[\varphi
(S_t)|{\cal F}_s \big] }{ E_x \big[\varphi (S_t) \big]}\Big] .$$

\noi  2) By Proposition \ref{posm1}, we know that $\dis \Big( {1
\over 1- \Phi(x)} M_t^{\varphi}\Big)_{t \ge 0}$ may be written as
an exponential martingale. Consequently, point 2) of Theorem
\ref{tosm1} is a direct consequence of Girsanov's theorem (cf
\cite{RevYor}, p. 311-313).

\end{prooff}


%
%
%



\noi We would like to generalize Theorem \ref{tosm1}, replacing
the normalization coefficient  $\varphi (S_t)$ by a function of
$(X_t, S_t)$. As  Theorem \ref{tosm2} below shows, a good
candidate is $\psi(S_{t}) e^{\lambda(S_{t}-X_{t})}$ where $\lambda
> 0$ and $\psi : \mathbb{R}_+ \rightarrow \mathbb{R}_+$ satisfies
some conditions.

\begin{theo}\label{tosm2}
 Let $\lambda >0$ and $\psi : \mathbb{R}\mapsto ] 0,+ \infty[$
 satisfying (\ref{osm6b1}). Let $\Phi$ be the function associated
 with $\psi$ via (\ref{osm6b2}) and $\varphi=\Phi'$.

\begin{enumerate}

    \item  Let $u>0$, $x \in \mathbb{R}$ and any $\Gamma_u$
in ${\cal F}_u$. Then :
\begin{equation}\label{osm13}
\lim_{t \rightarrow \infty} \frac{E_x \big[1_{\Gamma_u}
\psi(S_{t}) e^{\lambda (S_t-X_t)}\big] }{ E_x \big[\psi (S_t)
e^{\lambda (S_{t}- X_{t})}\big]} = \frac{1 }{1- \Phi(x)}
E_x[1_{\Gamma_u} M_u^{\lambda, \varphi}],
\end{equation}
\noi where $(M_u^{\lambda, \varphi})$ is the $P_x$-martingale
defined in Proposition \ref{posm2}.

\item Let $(Q_x^{\lambda ,\varphi})_{x \in \mathbb{R}}$ be the
family of probabilities on $\big(\Omega, {\cal F}_\infty\big)$ :
\begin{equation}\label{osm13a}
    Q_x^{\lambda ,\varphi} (\Gamma_u) = \frac{1 }{1- \Phi(x)} \, E_x [1_{\Gamma_u}
\, M_u^{\lambda ,\varphi}], \quad \mbox{for any} \ u \ge 0,  \
\mbox{and} \ \Gamma_u \in {\cal F}_u.
\end{equation}

\noi Then, under $Q_x^{\lambda ,\varphi}$, the process
$$
\Big(X_t-x + \lambda \int_0^t  \frac{\varphi (S_u) \cosh
\big(\lambda (S_u-X_u)\big)+\lambda \big(1- \Phi (S_u)\big) \sinh
\big(\lambda (S_u-X_u)\big) }{\lambda \big(1- \Phi (S_u)\big)
\cosh \big(\lambda (S_u-X_u)\big) + \varphi (S_u)\sinh
\big(\lambda (S_u-X_u)\big)} du\ ; t \ge 0 \Big)$$
\noi is a Brownian motion, started at $0$.

\item  In fact, there is the absolute continuity relationship :
\begin{equation}\label{osm15c}
Q_x^{\lambda ,\varphi}=\frac{1}{1-\Phi (x)}\Big(\frac{e^{\lambda
S_\infty}\psi(S_\infty)}{2\lambda}\Big)P_x^{(-\lambda)}.
\end{equation}
\end{enumerate}
\end{theo}

\noi Recall (cf 1) of Remark \ref{rosm2}) that $M^{0,
\varphi}=M^{\varphi}, M^{\varphi}$ being the martingale defined in
Proposition \ref{posm1}. Therefore,  Theorem  \ref{tosm1} may be
interpreted as a particular case of Theorem \ref{tosm2}, since
taking formally $\lambda =0$ in (\ref{osm13})
 we recover (\ref{osm7}). Note that 3. of Theorem
 \ref{tosm2} follows from Remark \ref{rosm1b}.

 \noi
Our proof of Theorem \ref{tosm2} is similar to that of Theorem
\ref{tosm1}. An extension of  Lemma \ref{losm1} is required to
obtain an equivalent of $E_0\big[\psi(s \vee (x+S_t)\big)
e^{\lambda \{ s \vee (x+S_t)-x-X_t)\}}\big]$, as $t \to \infty$.
This result is stated below in Lemma \ref{losm2}. We observe that
the  two asymptotic rates of growth are drastically different.
This explains why we state two separate results.

\begin{lemma} \label{losm2}
Let $s \ge x, \, s \ge 0$ and $\psi : \mathbb{R} \rightarrow
\mathbb{R}$ such that :
\begin{equation}\label{osm14}
\int_{s}^{\infty} |\psi(z)|e^{-\lambda z}dz< \infty,
\end{equation}
\begin{equation}\label{osm15}
\rho_{\lambda}(s,x):=\psi(s)\sinh \big(\lambda (s-x)\big)+ \lambda
e^{\lambda x} \int_s^{\infty} \psi (z)
 e^{- \lambda z} dz \neq 0.
\end{equation}
\noi Then :
$$E_0\Big[\psi\big(s \vee (x+S_t)\big) e^{\lambda (s_{\vee}
(x+S_t)-x-X_t)}\Big] _{\stackrel{\sim}{t\rightarrow\infty}} 2
\rho_{\lambda} (s,x) e^{\lambda^{2} t/2}.$$
\noi If $\psi \geq 0$, we have :
\begin{equation}\label{osm15a}
    E_0\Big[\psi\big(s \vee (x+S_t)\big) e^{\lambda (s_{\vee}
(x+S_t)-x-X_t)}\Big]  \leq 2 \rho_{\lambda} (s,x)\Big( 1+
\frac{1}{\lambda \sqrt{2 \pi t}}\Big) e^{\lambda^{2} t/2},
\end{equation}
\begin{equation}\label{osm15b}
    E_0\Big[\psi\big(s \vee (x+S_t)\big) e^{\lambda (s_{\vee}
(x+S_t)-x-X_t)}\Big]  \geq  e^{\lambda^{2} t/2} \frac{2\lambda
}{\sqrt{2 \pi}} \int_{-\infty}^{x-s}e^{-u^2/2}du \int_{s}^\infty
\psi(z)e^{-\lambda z}dz   ,
\end{equation}
\noi if $ t\geq 1$.
\end{lemma}

\begin{prooff} \ {\bf of Lemma \ref{losm2}}.  Let $\Delta$ be the
expectation of $\psi \big(s \vee (x+S_t)  e^{\{\lambda (s \vee
(x+S_t)-x-X_t)\} }$. We split $\Delta$ in two parts $\Delta_1$ and
$\Delta_2$, corresponding respectively to $\{S_t \le s-x\}$ and to
$\{S_t> s-x\}$ :
\begin{eqnarray}
\Delta_1 & = &\psi (s) \, e^{\lambda (s-x)} E_0 [e^{- \lambda
X_{t}} 1_{\{S_{t} < s-x\}}], \nonumber \\
 \Delta_2&=& E_0 \big[\psi (x+S_t)
e^{\lambda (S_t-X_t)} 1_{\{S_{t}
>s-x\}} \big]. \nonumber
\end{eqnarray}

\noi Recall (\cite{KaratShrev}, section 2.8, p 95) and
(\cite{RevYor} section III.3, p 105),  that under $P_0, \,
(S_t,X_t)$ is distributed as :
\begin{equation}\label{osm16}
P_0(S_t \in db, \; X_t \in da)=\frac{2(2b-a) }{ \sqrt{2 \pi
t^{3}}} \,e^{-\frac{(2b-a)^{2} }{ 2t}}\, 1_{\{a<b,\,b>0\}}\,dadb.
\end{equation}
1) Consequently :
$$\Delta_1=\frac{2 \psi(s) e^{\lambda (s-x)} }{ \sqrt{2 \pi t^{3}}}
\int_0^{s-x} db \Big(\int_{-Ê\infty}^b e^{- \lambda a}(2b-a)
e^{-\frac{(2b-a)^{2}}{ 2t}} da\Big).$$
\noi Setting $c=a-2b$ in the $a$-integral, we obtain :
$$\Delta_1=\frac{2 \psi (s) e^{\lambda (s-x)} }{ \sqrt{2 \pi t}} \int_0^{s-x}
e^{-2 \lambda b}\, A(t,b)db, $$

\noi where $\dis A(t,b)=-\frac{1 }{ t} \int_{- \infty}^{-b} e^{-
\lambda c} c e^{-c^{2}/2t}dc$.

\noi  Integrating  by parts, we obtain :
$$
    A(t,b)=e^{\lambda b}\,e^{-b^{2}/2t}+\lambda \int_{-
\infty}^{-b} e^{- \lambda c-c^{2}/2t} dc.
$$
\noi Setting $\dis u= \frac{c }{\sqrt{t}}+ \lambda \sqrt{t}$ we
get :
\begin{equation}\label{osm17}
    A(t,b)=e^{\lambda b}\,e^{-b^{2}/2t}+\lambda \sqrt{t}\, e^{\lambda^{2}
t/2} \int_{-\infty}^{-\frac{b }{ \sqrt{t}}+ \lambda \sqrt{t}} e^{-
u^{2}/2} du.
\end{equation}
\noi Since $ \dis \lim_{t \rightarrow \infty}
\int_{-\infty}^{\frac{-b }{ \sqrt{t}}+\lambda \sqrt{t}}\,
e^{-u^{2}/2}du = \sqrt{2 \pi}$, then  $A(t,b)_{\stackrel
{\sim}{t\rightarrow\infty}} \lambda \sqrt{2 \pi t}e^{\lambda ^2
t/2}$, and :
$$\Delta_{1} \  _{\stackrel{\sim}{t \rightarrow
\infty}} 2\psi(s) \sinh \big(\lambda(s-x)\big)\,e^{\lambda^{2}t
/2}.$$

\smallskip

\noi  2) Mimicking the approach developed in 1), we obtain :
$$\Delta_2=\frac{2 }{ \sqrt{2 \pi t}} \int_{s-x}^{+Ê\infty} \psi
(x+b) \; e^{- \lambda b} A(t,b) \,db.$$

\noi The decomposition (\ref{osm17}) and the finiteness hypothesis
(\ref{osm14}) imply that :
$$\Delta_2 \ _{\stackrel{\sim}{t \rightarrow
\infty}}  \Big(2 \lambda \int_{s-x}^{\infty} \psi (x+b) e^{-
\lambda b} db\Big)\,e^{\lambda^{2} t/2},$$

\noi 3) Suppose $\psi \geq 0$.

\noi a) Applying $\displaystyle \int_{-\infty}^{-\frac{b }{
\sqrt{t}}+ \lambda \sqrt{t}} e^{- u^{2}/2} du \leq \int_\mathbb{R}
e^{- u^{2}/2} du=\sqrt{2\pi}$, and  $\dis \lambda b
-\frac{b^2}{2t} \leq \frac{\lambda ^2 t}{2}$ in(\ref{osm17}) imply
that $A(t,b)\leq e^{\lambda ^2 t /2} (1+\lambda \sqrt{2 \pi t})$
and (\ref{osm15a}).

\noi b) Let $b\in [0,s-x]$ and $ t\geq 1$. Then $\displaystyle
-\frac{b}{\sqrt{t}}+\lambda \sqrt{t}\geq -\frac{b}{\sqrt{t}} \geq
x-s$ and $\displaystyle A(t,b)\geq e^{\lambda^{2} t/2}  \lambda
\sqrt{t}\int_{-\infty}^{x-s}e^{-u^2/2}du$. Using moreover $\Delta
\geq \Delta_2$ we obtain (\ref{osm15b}).

\end{prooff}

\begin{prooff} \ {\bf  of Theorem \ref{tosm2}}

\noi  Let $u>0,\, x \in \mathbb{R}$ and $\Gamma_u \in {\cal F}_u$.
Adapting the proof of Theorem \ref{tosm1} to our new context, we
have :
$$E_x\big[\psi(S_t) e^{\lambda (S_{t}-X_{t})}|{\cal F}_u\big]=
\widetilde{\psi} (S_u,X_u ; t-u),$$

\noi where :
$$\widetilde{\psi} (s,x ; r)=E_0\Big[\psi\big(s \vee (x+S_r)\big)
\,e^{\lambda \big(s \vee (x+S_r)-x-X_r\big)}\Big].$$

\noi Lemma \ref{losm2} gives the rate of increase of
$\widetilde{\psi} (s,x ; r), \; r \to \infty$, if $\psi \geq 0$  :
$$\widetilde{\psi} (s,x ; r) _{\stackrel{ \sim}{r \rightarrow \infty}} 2
\rho_{\lambda} (s,x) \, e^{\lambda^{2} r/2}.$$

\noi In particular, taking  $u=0$, we get :
$$E_x \big[\psi (S_t)\,e^{\lambda (S_t-X_t)}\big]
 _{\stackrel{ \sim}{t \rightarrow \infty}} 2\, \rho_{\lambda}
(x,x) \, e^{\lambda^{2}t/2}.$$
\noi Moreover if $t>u+1$, (\ref{osm15a}) and (\ref{osm15b}) imply
that :
$$\frac{\widetilde{\psi} (S_u,X_u ; t-u)}{E_x \big[\psi (S_t)\,e^{\lambda (S_t-X_t)}\big]}
\leq \frac{k}{\lambda} \rho_\lambda(S_u,X_u)e^{-\lambda u^2/2},$$
where $k$ is a constant depending only on $x,u,\lambda$.

\noi But $\dis
M^{\lambda,\varphi}_u=\frac{1}{\lambda}\rho_\lambda(S_u,X_u)e^{-\lambda
u^2/2}$ and $(M^{\lambda,\varphi}_t)$ is a $P_x$-martingale;
consequently :
$$ \lim_{t \rightarrow \infty} \frac{E_x \big[1_{\Gamma_u} \psi (S_t)
\,e^{\lambda (S_{t}-X_{t})}\big] }{ E_x \big[\psi (S_t)\,
e^{\lambda (S_{t}-X_{t})}\big]}=\frac{1 }{ \rho_{\lambda}
(x,x)}\,E_x\big[ 1_{\Gamma_u} \rho_{\lambda}
(S_u,X_u)\,e^{\lambda^{2} u/2}\big].$$

%

\end{prooff}

\subsection{ Case 2 : the local time at $0$}\label{loc}

\noi From L\'{e}vy's theorem, under $P_0$, $\big((S_t-X_t, S_t) ;
t \ge 0 \big)$ and $\big(|X_t|,\,L_t^0) ; t\ge 0\big)$ have the
same distribution. This implies that $\big(|X_t|\, \varphi
(L_t^0)+1-\Phi(L_t^0)\;;\;t \ge 0\big)$ is a martingale with
respect to the filtration of $(|X_t|)$, hence with respect to
$({\cal F}_t)$, the functions $\varphi$ and $\Phi$ being defined
in Proposition \ref{posm1}. These processes are particular cases
of more general martingales :

\begin{prop}\label{ploc1}
 Let $h^+,h^- : \mathbb{R}_+
\longrightarrow \mathbb{R}_+$, be bounded, Borel functions, and
define :
\begin{equation}\label{loc1}
    H(l)={1 \over 2} \int_0^l \big(h^+ (u)+h^-(u)\big)du, \quad
    l\geq 0.
\end{equation}
\begin{enumerate}
    \item Then :
\begin{equation}\label{loc2}
    M_t^{h^{+}, h^{-}} = 1-H(L_t^0)+X_t^+ h^+(L_t^0)+X_t^- h^-
(L_t^0),
\end{equation}
\noi is a $P_0$-martingale. Moreover $M_0^{h^{+},h^{-}}=1$ $P_0$
a.s., and :
\begin{equation}\label{loc3}
    M_t^{h^{+},h^{-}} =1+\int_0^t \Big(1_{\{X_{s}>0\}} h^+(L^0_s)-
1_{\{X_{s} <0\}} h^-(L^0_s)\Big)d X_s.
\end{equation}

    \item If moreover :
    \begin{equation}\label{loc4}
    \frac{1}{2}\int_0^{\infty} \big(h^+(u)+h^-(u)\big)du=1,
\end{equation}
\noi and $H(l)<1$, for any $l>0$, then $M_t^{h^{+},h^{-}} >0$ and
$M_t^{h^{+},h^{-}} = {\cal E}(J^{h^{+},h^{-}})_t$, where :
\begin{equation}\label{loc5}
    J_t^{h^{+},h^{-}}={1_{\{X_t>0\}} h^{+}(L^0_{t})-1_{\{X_t<0\}}h^{-}
(L^0_{t}) \over 1-H(L^0_{t})+X_t^{+} h^{+} (L^0_t)+X_t^-
h^-(L^0_t)}.
\end{equation}
\end{enumerate}

\end{prop}

\noi  The martingales $(M_t^{h^{+},h^{-}})$ featured in
(\ref{loc2}) and (\ref{loc3})  have already been used, e.g., in
\cite{JY}. Both statements of Propositions \ref{posm1},
\ref{ploc1}  are also found in (\cite{RevYor}, Chapter VI, "first
order calculus") and are particular cases of application of the
balayage formula. More precisely, formula (\ref{loc3}) may be
generalized as follows : if $(h_s^{+} ; s \ge 0)$ and $(h_s^{-} ;
s \ge 0)$ are bounded predictable processes, and if $g_t = {\rm
sup} \{s \le t ; X_s=0\}$, then :
\begin{equation}\label{loc6}
X_t^+ h_{g_{t}}^{+} + X_t^- h_{g_{t}}^{-} - \frac{1}{2} \int_0^t
(h_s^{+}+ h_s^{-})dL_s^0  = \int_0^{t} (h_{g_{s}}^{+}1_{\{X_s
>0\}}-h_{g_{s}}^{-} 1_{\{X_{s} <0\}}) dX_s.
\end{equation}

\begin{rem}\begin{enumerate}\label{rloc1}
    \item Recall that in Proposition \ref{posm2} we have
introduced the family of martingales $(M^{\lambda, \varphi}_t)$;
this should not induce any confusion with the family
$(M^{h^{+},h^{-}}_t)$, since the two parameters indexing the first
(resp. second) family are respectively $(\lambda, \varphi)$ and
 $(h^+,h^-)$, and belong to quite different sets.

    \item In this section we restrict ourselves to $P_0$. Thus $(X_t)$
is a Brownian motion started at $0$. It would also be possible to
work under $P_x$, replacing $(M^{h^{+},h^{-}}_t)$ by
$\big(1-H(L_t^x) +(X_t-x)^+ h^+ (L_t^x)+(X_t -x)^- h^- (L_t^x)
\;;\; t \ge 0\big)$. However, for simplicity we only deal with
$x=0$.
\end{enumerate}

\end{rem}

\noi We now investigate penalizations involving the local time at
$0$ of $X$.

\begin{theo}\label{tloc1}
Let $(h^+,h^-)$ and $(M^{h^{+}, h^{-}}_t)$ be the functions and
the  martingale  defined in Proposition \ref{ploc1}. We suppose
that (\ref{loc4}) holds.
\begin{enumerate}
    \item  Let $s \ge 0$ and $\Gamma_s \in {\cal F}_s$. Then :
\begin{equation}\label{loc7}
\lim_{t \to \infty}\frac{E_0\Big[1_{\Gamma_s} \big(h^{+}
(L_{t}^{0}) 1_{\{X_{t} >0\}}+ h^{-}(L_{t}^{0}) 1_{\{X_{t}
<0\}}\big)\Big]} {E_0\big[h^{+} (L_{t}^{0}) 1_{\{X_{t} <0\}} +
h^{-}(L_{t}^{0}) 1_{\{X_{t} <0\}}\big]} =E_0[1_{\Gamma_s} \;
M_{s}^{h^{+},h^{-}}].
\end{equation}
    \item  Let $Q_0^{h^{+},h^{-}}$ be the probability
measure on $\big(\Omega, {\cal F}_\infty\big)$ satisfying :
$$Q_0^{h^{+},h^{-}} (\Gamma_s)=E_0 [1_{\Gamma_s} \; M_s^{h^{+},h^{-}}],$$

\noi  for any $s \ge 0$ and $\Gamma_s \in {\cal F}_s$.

\noi Then under $Q_0^{h^{+},h^{-}}$, the process :
$$\Big(X_t-\int_0^t \frac{h^+(L_s^0) 1_{\{X_{s} >0\}} - h^{-} (L_s^0)
1_{\{X_{s} <0\}} }{ M_s^{h^{+},h^{-}}} \,ds \;;\; t \ge 0 \Big)$$
\noi is a Brownian motion.
\end{enumerate}
\end{theo}

\begin{rem}\label{rloc2}  As observed in Remark \ref{rosm3}, the $P_0$-martingale $(M_t^{h^{+},h^{-}})$
is not uniformly integrable, in fact, if $\tau_{l}={\rm inf}\{s
>0, \, L_s^0>l\}$, then :
$M_{\tau_{l}}^{h^{+},h^{-}}$ goes a.s. to $0$ as $l \to \infty$.
\end{rem}

\noi The proof of part 1) of Theorem \ref{tloc1} is based on the
following estimate.

\begin{lemma} \label{lloc1} Let $f : \mathbb{R}_+ \longrightarrow \mathbb{R}_+$
be Borel and locally bounded such that $\dis \int_0^{\infty} f(s)
ds < \infty$. Let $a \ge 0$ and $x \in \mathbb{R}$, then
\begin{equation}\label{loc8}
E_x \big[f (a+L_t^0)\,1_{\{X_t >0\}}\big]
_{\stackrel{\sim}{t\rightarrow \infty}} f(a)\sqrt{\frac{2 }{ \pi
t}} \ x^+ + \frac{1 }{ \sqrt{2 \pi t}} \int_a^{\infty} f(s) ds,
\end{equation}
\begin{equation}\label{loc9}
    E_x \big[f (a+L_t^0)\,1_{\{X_t <0\}}\big]
_{\stackrel{\sim}{t\rightarrow \infty}} f(a)\sqrt{\frac{2 }{ \pi
t}} \ x^- + \frac{1 }{ \sqrt{2 \pi t}} \int_a^{\infty} f(s) ds.
\end{equation}
\end{lemma}

\begin{prooff} \ {\bf of Lemma \ref{lloc1}}. Since under $P_x$, $(-X_t)_{t \ge 0}$  is distributed as
$(X_t)$ under $P_{-x}$, (\ref{loc9}) is a direct consequence of
(\ref{loc8}).

\noi To prove (\ref{loc8}), we recall the well-known result :
\begin{equation}\label{loc10}
    P_x (T_a \in dt) = \frac{|x-a| }{ \sqrt{2 \pi t^{3}}}\,\exp
\Big\{-\frac{(x-a)^2 }{ 2t}\Big\} 1_{\{t >0\}} dt.
\end{equation}

\noi On the set $\{T_0 > t\}$, we have, $P_x$ a.s. :
$$f(a+L_t^0) 1_{\{X_t > 0\}} =
\left\{
\begin{array}{cl}
  0 & \mbox{if} \ x<0, \\
  f(a) & \mbox{otherwise} \\
\end{array}
\right.
$$

\noi Then
$$E_x \big[f(a+L_t^0) 1_{\{X_t >0\}}\big]=f(a) P_x (T_0>t)
1_{\{x >0\}} + \Delta_t,$$
\noi where
$$\Delta_t = E_x \big[ f(a+L_t^0)1_{\{X_t>0, \, t\geq T_{0}\}}\big].$$

\noi Using (\ref{osm10}) we have :
$$P_x(T_0 >t)=P_0 (T_{|x|}>t)=P_0(S_t<|x|)=P_0\Big(|X_1|<{|x|
\over \sqrt{t}}\Big)\cdot$$

\noi Therefore :
\begin{equation}\label{loc11}
    P_x (T_0>t) _{\stackrel{\sim}{t \rightarrow + \infty} }\sqrt{\frac{2 }{ \pi t}}
\ |x|.
\end{equation}

\noi To compute $\Delta_t$, we use the strong Markov property at
time $T_0$, and we get :
$$\Delta_t= E_x\big[1_{\{t\geq T_{0}\}}\, g(t-T_{0})\big],$$
\noi with $g(r)=E_0 \big[f(a+L_r^0)\,1_{\{X_r > 0\}}\big]$.

\noi Since under $P_0$, $(-X_u)$ and $(X_u)$ have the same
distribution,
$$g(r) = \frac{1}{ 2}E_0 \big[f(a+L^0_r)\big]= \frac{1}{ 2}E_0 \big[f(a+S_r)\big],$$
\noi the last equality being a consequence of   L\'{e}vy's
theorem.

\noi Applying Lemma \ref{losm1} (with $s=x=a, \, u=r$ and
$\varphi_0 =f$), we obtain :

$$g(r) _{\stackrel{\sim}{r \rightarrow  \infty} } \frac{1}{ \sqrt{2 \pi
r}} \int_a^{\infty} f(s) ds.$$

\noi Consequently, $\displaystyle \Delta_t \  _{\stackrel{\sim}{t
\rightarrow  \infty} } \frac{1 }{ \sqrt{ 2\pi t} }\int_a^{\infty}
f(s) ds$. (\ref{loc8}) follows immediately.
\end{prooff}

\vskip 1cm

\noi The remainder of the proof of Theorem \ref{tloc1} is left to
the reader, since it is similar to the proofs of Theorems
\ref{tosm1} and \ref{tosm2}.


\subsection{ Case 3 : the maximum, the infimum and the local time.}
\label{mil}



 \noi For simplicity we restrict ourselves
to $P_0$. The family of martingales playing a central role in this
section may be new. This family is  related to martingales of the
type $(M_t^{h^{+},h^{-}})$. More precisely, let :
\begin{equation}\label{mil1}
    h^+(l)=\frac{1 }{ a} \, e^{cl}, \, h^-(l)=\frac{1 }{ b}\,e^{cl},
\end{equation}
\noi with $c=\dis\frac {1 }{ 2} \Big(\frac{1 }{ a}+\frac{1 }{
b}\Big),\, a>0$ and $b>0$. Then $H(l)=e^{cl}-1$ and
\begin{equation}\label{mil2}
    M_t^{h^{+},h^{-}}=2-\Big(1- \frac{X_t^+ }{ a}-\frac{X_t^- }{ b}\Big)
e^{cL_{t}^{0}}.
\end{equation}
\noi Since $X_t^+ X_t^-=0$, then :
\begin{equation}\label{mil3}
    M_t^{h^{+},h^{-}}=2-\Big(1-\frac{X_t^+ }{a}\Big) \Big(1-
\frac{X_t^- }{ b}\Big) e^{cL_{t}^{0}}.
\end{equation}
\noi In particular,
\begin{equation}\label{mil4}
    2-M_{t \wedge T_{a} \wedge T_{-b}}^{h^{+},h^{-}}=
\Big(1-\frac{X^{+}_{t} }{ a}\Big)\,\Big(1-\frac{X_{t}^- }{ b}
\Big) e^{cL_{t }^0}\ 1_{\{t\leq  T_{a} \wedge T_{-b}\}},
\end{equation}
\noi is a $P_0$-martingale.

\noi Integrating this identity with respect to some positive
measure  $\nu(da,db)$ we obtain the following result.

\begin{prop}\label{pmil1}
Let $\nu$ be a probability measure on $\mathbb{R}_+ \times
\mathbb{R}_+$, whose support is  included in $[\alpha, + \infty [
\times [\alpha, +\infty[$, for some $\alpha
>0$, and $(M_t^{ \nu})$ be the process :
\begin{eqnarray}
M_t^{\nu} & \dis =  \int_{\mathbb{R}_+\times \mathbb{R}_+}
\Big(1-\frac{X_t^+ }{ a}\Big)\Big(1-\frac{X_t^- }{
b}\Big)\exp\Big\{\frac{1  }{ 2}\Big(\frac{1 }{ a}+\frac{1 }{
b}\Big) L_t^0\Big\} 1_{\{t \le T_{a} \wedge T_{-b}\}}\nu(da,db)
\label{mil5}\\
&\dis =\int_{\mathbb{R}_+\times \mathbb{R}_+} \Big(1-\frac{X_t^+
}{ a}\Big)\Big(1-\frac{X_t^- }{ b}\Big)\exp\Big\{\frac{1  }{
2}\Big(\frac{1 }{ a}+\frac{1 }{ b}\Big) L_t^0\Big\} 1_{\{ S_t \le
a, \, I_t \le b\}} \nu (da,db).\nonumber
\end{eqnarray}
\noi Then, under $P_0$, $(M_t^{\nu}\;;\;t \ge 0)$ is a positive
martingale, $M_0^{\nu}=1$, $M_t^{\nu}={\cal E}(J^{\nu})_t$ with :
\begin{equation}\label{mil6}
    J_t^{\nu} = \frac{1 }{ M_{t}^{\nu}} \int_{\mathbb{R}_+\times \mathbb{R}_+}
\Big(-\frac{1 }{ a} 1_{\{X_t >0\}} +\frac{1
}{b}1_{\{X_{t}<0\}}\Big) \exp\Big\{\frac{1 }{ 2} \Big(\frac{1 }{
a} +\frac{1}{ b} \Big) L_t^0\Big\}  1_{\{ S_t \le a, \, I_t \le
b\}}\nu (da,db).
\end{equation}
\end{prop}

\begin{rem}\label{rmil1}\begin{enumerate}
    \item Obviously if we take for $\nu$ the Dirac measure at
$(a,b)$, for some $a>0,b>0$, then we recover (\ref{mil4}).

    \item We may write $(M_t^{\nu})$ as follows :
    \begin{equation}\label{mil7}
    M_t^{\nu} = F(S_t,I_t,L_t^0)-X_t^+ F^+ (S_t,I_t,L_t^0)-X_t^-
F^-(S_t,I_t,L_t^0),
\end{equation}
with
\begin{equation}\label{mil8}
    F^+(s,i,l)=\int_{\mathbb{R}_{+}^{2}} 1_{\{s \le a ,i \le b\}}\frac{1
    }{
a} \,{\rm exp} \Big\{\frac{1 }{ 2}\Big(\frac{1 }{ a} +\frac{1 }{
b}\Big) l \Big\} \, \nu (da,db),
\end{equation}
 \begin{equation}\label{mil9}
    F^-(s,i,l)=\int_{\mathbb{R}_{+}^{2}} 1_{\{s \le a ,i \le b\}}\frac{1
    }{
b} \,\exp \Big\{\frac{1 }{ 2}\Big(\frac{1 }{ a} +\frac{1 }{
b}\Big) l \Big\} \, \nu (da,db),
\end{equation}
\begin{equation}\label{mil10}
    F(s,i,l)=\int_{\mathbb{R}_{+}^{2}} 1_{\{s \le a ,i \le b\}}
 \exp \Big\{\frac{1 }{ 2}\Big(\frac{1 }{ a} +\frac{1 }{
b}\Big) l \Big\} \, \nu (da,db).
\end{equation}

\noi If $\nu$ is assumed to have a continuous density function,
then it follows that :
\begin{equation}\label{mil11}
    \frac{1 }{ 2} \, (F^+ + F^-)=\frac{\partial F }{\partial l},
\end{equation}
\begin{equation}\label{mil12}
     s \frac{\partial F^{+} }{ \partial s} (s,i,l)=\frac{\partial
F }{
\partial s} (s,i,l) \;;\; i \frac{\partial F^{-}}{ \partial i}
(s,i,l)=\frac{\partial F }{ \partial i} (s,i,l).
\end{equation}

    \item Suppose that $\nu$ is a p.m. with support included
in the diagonal, then $(M_t^{\nu})$ coincides with $M^{\nu_{*}}_t$
where
\begin{equation}\label{mil13}
    M_t^{\nu_{*}}:= \int_{X_{t}^{*}}^{\infty} \Big(1-\frac{|X_t|
    }{a}\Big)\,e^{L^0_{t}/a} \nu_{*}(da),
\end{equation}
\noi $\nu_*$ being, in this case, a p.m. on $[\alpha, +\infty[,\,
\alpha >0$. Consequently, relations (\ref{mil7})- (\ref{mil10})
become :
\begin{equation}\label{mil14}
    M_t^{\nu_{*}}=F(X_t^*, \, L_t^0)-|X_t|\,F^*(X_t^*,\,L_t^0),
\end{equation}
\noi with :
\begin{equation}\label{mil15}
    F(x,l)=\int_x^{\infty}  e^{l/a} \nu_* (da),\quad
F^*(x,l)=\int_x^{\infty}\frac{1 }{ a} e^{l/a} \nu_* (da).
\end{equation}
\noi Moreover $M_t^{\nu_{*}}={\cal E}(J^{\nu_{*}})_t$ with :
\begin{equation}\label{mil17}
    J_t^{\nu_{*}}=-\frac{sgn  (X_t) }{ M_t^{ \nu_*}}
    \int_{X_{t}^{*}}^\infty
\frac{1 }{ a} \, e^{L_{t}^{0}/a} \nu_* (da).
\end{equation}

    \item Sometimes, we shall consider similar martingales involving
only the one-sided maximum $(S_t)$ and $(L^0_t)$; hence,  we shall
use, instead of (\ref{mil4}) :
$$2-M_{t \wedge T_{a}}^{h^{+},0}=\Big(1-\frac{X_{t \wedge T_{a}}^{+}
}{ a}\Big) e^{L_{t \wedge T_{a}}^0/2a},$$

\noi where $h^+(l)=\dis {1 \over a}\, e^{l/2a}$.

\noi Integrating over $[\alpha, \infty[$, with respect to a p.m.
$\nu_+$ on $[\alpha,\infty[$, leads to :
$$M_t^{\nu_{+}}=\int_{S_{t}}^{\infty} \Big(1-\frac{X_{t}^{+} }{ a}
\Big)\, e^{L^0_{t}/2a} \nu_+ (da).$$

\noi As previously $(M_t^{\nu_{+}})$ is a $P_0$-martingale, such
that $M_0^{\nu_{+}} =1$, and it may be written as :
$$M_t^{\nu_{+}} =F(S_t,L^0_t)-X_t^+ \, F^+ (S_t,L^0_t),$$

\noi with $F(s,l)=\dis \int_s^{\infty} e^{l/2a} \nu_+ (da)$ and
$F^+ (s,l)=\dis \int_s^{\infty} {1 \over a}\, e^{l/2a} \nu_+
(da)$.
\end{enumerate}
\end{rem}

\noi We now deal with penalizations involving jointly the maximum,
the minimum and the local time of $X$. Recall (see Proposition
\ref{pmil1}) that $\nu$ is a p.m. on $[\alpha, + \infty[ \times
[\alpha, + \infty[$ for some $\alpha >0$ and $(M_t^{\nu})$ is the
positive martingale defined by (\ref{mil5}).

\begin{theo}\label{tmil1}
 Let $A_{\nu}$ be the function :
 \begin{equation}\label{mil18}
    A_{\nu} (s,i,l)=\int_{\mathbb{R}_{+}^{2}} e^{\frac{1 }{ 2}\big({1
\over a}+\frac{1 }{ b}\big)l} 1_{\{s \le a ,i \le b\}} \nu
(da,db), \ s,i,l \geq 0.
\end{equation}
\begin{enumerate}
    \item Let $u \ge 0$ and $\Gamma_u \in {\cal F}_u$. Then :
\begin{equation}\label{mil19}
    \lim _{t \rightarrow \infty} \; \frac{E_0 \big[1_{\Gamma_u}
A_{\nu} (S_t,I_t,L_t^0)\big] }{ E_0 \big[A_{\nu} (S_t,I_t,
L_t^0)\big]} = E_0 [1_{\Gamma_u} \, M_u^{\nu}\big].
\end{equation}

    \item Let $Q_0^{\nu}$ be the p.m. on $\big(\Omega,
{\cal F}_\infty \big)$ which satisfies :
\begin{equation}\label{mil19a}
    Q_0^{\nu} (\Gamma_u) = E_0 [1_{\Gamma_u} M_u^{\nu}],
\end{equation}
\noi  for any $u \ge 0$ and $\Gamma_u \in {\cal F}_u$.

\noi The process $\Big(X_t - \dis \int_0^t \frac{J_s^{\nu} }{
M_s^{\nu}} ds ; t \ge 0\big)$ is a $Q_0^{\nu}$-Brownian motion,
where $(J_t^{\nu})$ is defined by (\ref{mil6}).
\end{enumerate}
\end{theo}

\noi The asymptotic result (\ref{mil19}) is based on the following
rather striking result.

\begin{lemma} \label{lmil1}
Let $a>0,b>0$ and $c=\dis \big(\frac{1 }{ a}+\frac{1 }{b}\big)$.
Then
\begin{equation}\label{mil20}
    \lim _{t \to \infty} E_0 [e^{cL^0_{t}} 1_{\{t<T_{-b}
\wedge T_{a}\}}]={3 \over 2},
\end{equation}
 \begin{equation}\label{mil21}
   \sup _{t \ge 0} E_0 [e^{cL^0_{t}}
1_{\{t<T_{-b} \wedge T_{a}\}}]<+ \infty.
\end{equation}
\end{lemma}

\noi It would be possible to directly  prove Lemma \ref{lmil1},
but the proof is technical \footnote{ We thank F. Petit who helped
us with such a proof}. Instead, we provide a short proof based on
a disintegration of the p.m. $Q_0^\nu$ (see Theorem \ref{tdir1}
for details). This is why  we have postponed this proof to Section
\ref{dir} (which is devoted to such disintegrations), just after
Theorem \ref{tdir1}.

\begin{prooff} \ {\bf of Theorem \ref{tmil1}}.

\noi Let $u \ge 0$ and $\Gamma_u \in {\cal F}_u$ be fixed, $t>u$
and $\delta (\Gamma_u,t)=E_0\big[1_{\Gamma_u} A_{\nu}
(S_t,I_t,L^0_t)\big]$.

\noi By (\ref{mil18}),  we have :
$$A_{\nu} (S_t,I_t,L^0_t)=\int_{\mathbb{R}_{+}^{2}} e^{\frac{1 }{ 2}\big(
\frac{1 }{ a} + \frac{1 }{ b}\big)L^0_t} 1_{\{T_{a}>t, \,
T_{-b}>t\}} \nu(da,db).$$
\noi Applying the Markov property at time $u$, we obtain :
$$\delta (\Gamma_u, t) = \int_{\mathbb{R}_{+}^{2}}E_0 \big[1_{\Gamma_u}
1_{\{T_{a} \wedge T_{-b}>u\}} e^{cL^0_u} \pi_1 (a,b,X_u,t-u)\big]
\nu (da,db).$$
\noi where $c=\dis \frac{1 }{ 2} \Big(\frac{1 }{ a}+\frac{1 }{
b}\Big)$, and
$$\pi_1 (a,b,x,r) = E_x \big[1_{\{T_{a} \wedge T_{-b}>r\}}
e^{cL_{r}^{0}}\big] \;;\;-b \le x \le a, \, r \ge 0.$$
\noi Let us start with $x>0$.

\noi Since under $P_x$, $\{T_a \wedge T_{-b}>r\}=\{T_a \wedge T_0
>r\} \cup \{T_0<T_a, \,  T_a \wedge T_{-b}>r\}$
we have :
$$
\begin{array}{ccl}
 \pi_1 (a,b,x,r) & = & P_x(T_a \wedge T_0>r)+\pi_2 (a,b,x,r), \\
  \pi_2(a,b,x,r) & = & E_x [1_{\{T_{0}<T_{a},\,T_{a} \wedge T_{-b}
>r\}} e^{cL_{r}^{0}}]. \\
\end{array}
$$

\noi Applying the strong Markov property at time $T_0$, we get :
$$\pi_2 (a,b,x,r)=E_x \big[1_{\{T_{0}<T_{a}, \ T_0<r\}} \pi_3 (a,
b,r- T_0)\big]$$
\noi with :
$$\pi_3(a,b,r')=E_0[e^{cL_{r'}^{0}} 1_{\{r'<T_{-b}Ê\wedge T_{a}
\}}].$$

\noi Lemma \ref{lmil1} implies that $\pi_2$ is bounded and $\pi_2
(a,b,x,r)$ converges to $\dis \frac{3 }{ 2} \, \frac{a-x }{ a}$,
as $r \to + \infty$.

\noi As for $P_x (T_a \wedge T_0>r)$, it is clear that this
probability goes to $0$, as $r \to \infty$.

\noi A similar approach may be developed for $x<0$. Inequality
 (\ref{mil21}) implies that $\pi_1$ is bounded ; consequently, for $\Gamma_u \in
 {\cal F}_u$ :
$$
\lim_{t \to \infty}\delta (\Gamma_u, t) = \frac{3 }{ 2} E_0
\Big[1_{\Gamma_u} \int_{\mathbb{R}_{+}^{2}} e^{cL_{u}^{0}}
1_{\{T_{a} \wedge T_{-b}>u\}} \Big\{\frac{a-X_{u} }{ a}
1_{\{X_{u}\geq 0\}} + \frac{b + X_{u} }{ b} 1_{\{X_{u}<0\}}\Big\}
\nu (da,db)\Big].$$
\noi The result follows because $\dis \frac{a-X_{u}}{ a}1_{\{X_{u}
\geq 0\}}+\frac{b+X_{u} }{b}1_{\{X_{u}<0\}}$ is equal to $\dis
\Big( 1-\frac{X_{u}^{+}}{ a}\Big) \Big(1-\frac{X_{u}^{-}}{
b}\Big)$ and $E_0[M_u^{\nu}]=1$.
\end{prooff}


\subsection{ Case 4 : the down-crossings.}\label{do}


\noi Let $a$ and $b$ be two fixed numbers, $a<b$. To describe
 the down-crossings of $(X_t)$,  from level $b$ to level $a$,
 it is convenient to
introduce the  sequence of stopping times defined inductively as
follows :
\begin{equation}\label{do1}
    \sigma_1 = \inf\{t \ge 0, \; X_t > b\},
\end{equation}
  \begin{equation}\label{do2}
\sigma_2 = \inf \{t \ge \sigma_1, \; X_t < a\},
\end{equation}
\noi and
\begin{equation}\label{do3}
    \sigma_{2n+1}={\rm inf} \{t\ge \sigma_{2n}\;;\; X_t>b\},
\end{equation}
\begin{equation}\label{do4}
    \sigma_{2n+2}=\inf \{t \ge \sigma_{2n+1} \;;\; X_t<a\}.
\end{equation}
\noi Let $D_t$ be the number of down-crossings from level $b$ to
level $a$, up to time $t$ :
\begin{equation}\label{do5}
    D_t = \sum_{n \ge 1} 1_{\{\sigma_{2n} \le t\}}\;;\;t\ge 0.
\end{equation}

\noi We observe that :
\begin{equation}\label{do6}
    \sigma_{2D_{s}}=\sup\{n \ge 1\;;\; \sigma_{2n}\le s\},
\end{equation}

\noi with the convention $\sup\{\emptyset\} = 0$ and $\sigma_0=0$.

\noi The events $\{\sigma_{2D_{s}} +T_b \circ
\theta_{\sigma_{2D_{s}}}
>s\}$ and $\{\sigma_{2D_{s}} + T_b \circ \theta_{\sigma_{2D_{s}}}
\le s\}$ will play a central role below. If $\{\sigma_{2n} \le s
<\sigma_{2n+2}\}$, the first (resp. second) event reduces to $\{
\sigma_{2n} \le s < \sigma_{2n+1}\}$ \big(resp. $\{\sigma_{2n+1}
\le s <\sigma_{2n+2}\}$\big).

\begin{prop} \label{pdo1}
 Let $\big(G(n)\big)_{n\geq 0}$ be a sequence of real numbers, and $(M_t^{\downarrow,G})$ be the
process :
\begin{eqnarray}
  M_t^{\downarrow,G}= & \dis 1_{\{\sigma_{2D_{t}}+T_b \circ
\theta_{\sigma_{2D_{t}}}
>t\}} \biggr(\frac{G(D_{t}) }{ 2}\,\Big(1+\frac{b-X_{t} }{ b-a}\Big)
+\frac{G(1+D_{t}) }{ 2}\Big(\frac{X_{t}-a }{ b-a}\Big)\biggl)\label{do7} \\
   & \dis +1_{\{\sigma_{2D_{t}} +T_b \circ \theta_{\sigma_{2D_{t}}}
\le t\}} \biggr(\frac{G(1+D_{t}) }{ 2}\Big(1+\frac{b-X_{t} }{
b-a}\Big) +\frac{G(D_{t}) }{ 2}\Big(\frac{X_{t}-a }{
b-a}\Big)\biggl). \nonumber
\end{eqnarray}
\begin{enumerate}
    \item Then $(M_t^{\downarrow,G})$ is a continuous
$P_x$-local martingale,
\begin{eqnarray}
 M_t^{\downarrow,G}=&\dis \sum_{n \ge 0}\biggr\{1_{[\sigma_{2n},
\sigma_{2n+1}[}(t) \biggr(\frac{G(n) }{ 2}\Big(1+\frac{b-X_{t} }{
b-a} \Big)+\frac{G(1+n) }{ 2} \Big(\frac{X_{t}-a }{
b-a}\Big)\biggl) \label{do8}
\\
&\dis +1_{[\sigma_{2n+1},\sigma_{2n+2}[} (t)\biggr(\frac{G(1+n)}{
2} \Big(1+\frac{b-X_{t} }{ b-a}\Big)+\frac{G(n)}{
2}\Big(\frac{X_{t}-a }{ b-a}\Big)\biggr\},\nonumber
\end{eqnarray}
\begin{eqnarray}
M_t^{\downarrow,G}=&\dis M_0^{\downarrow,G}+\frac{1 }{ b-a}
\int_0^t \Big(1_{\{\sigma_{2D_{s}}+T_{b} \circ
\theta_{\sigma_{2D_{s}}}
>s\}}-1_{\{\sigma_{2D_{s}}+T_{b} \circ \theta_{\sigma_{2D_{s}}}
 <s\}}\Big)
\label{do9}\\
 &\dis \times \Big(\frac{G(1+D_{s})-G(D_{s}) }{
2}\Big)\, dX_{s},\nonumber
\end{eqnarray}
\noi and under $P_x$ :
\begin{equation}\label{do9a}
    M_0^{\downarrow,G} =
    \left\{
\begin{array}{cl}
\dis \frac{1 }{ 2(b-a)}\big(G(0)(2b-a-x)+G(1)(x-a)\big) & \mbox{
if }
x\le b,\\
    \dis \frac{1 }{ 2(b-a)}\big(G(0)(x-a)+G(1)(2b-a-x)\big)&
    \mbox{otherwise.} \\
\end{array}
 \right.
\end{equation}

    \item Suppose moreover that $\big(G(n)\big)_{n \ge 0}$ is
a decreasing sequence of positive real numbers such that $G(0)=1$
and ${\rm lim}_{n \to \infty} G(n)=0$. Then $M_t^{\downarrow,G}
>0$ and $\dis M_t^{\downarrow,G}= M_0^{\downarrow,G}{\cal E}
(J^{\downarrow, G})_t$ with
\begin{equation} \label{do10}
J^{\downarrow, G}_t=
   \big(G(1+D_t)-G(D_t)\big)\biggr\{\frac{1 }{ G(D_t)(2b-X_t-a)+
G(1+D_t)(X_t-a)} 1_{\{\sigma_{2D_{t}}+T_{b}\circ
\theta_{2D_{t}}>t\}}
\end{equation}
 $$ - \frac{1
}{ G(1+D_{t})(2b-X_{t}-a)+G(D_{t})(X_{t}-a)}
1_{\{\sigma_{2D_{t}}+T_{b}\circ \theta_{2D_{t}}\leq t\}} \biggl\}.
$$
\end{enumerate}
\end{prop}

\begin{prooff}

\noi  1) A priori, $(M_t^{\downarrow,G})$ is a right continuous
process, and may only  jump at times $\sigma_n$. However, it is
easy to check that
$M_{\sigma_{n^{-}}}^{\downarrow,G}=M_{\sigma_{n}}^{\downarrow,G}$,
therefore $(M_t^{\downarrow,G} ; t \ge 0)$ is continuous. (A
posteriori, this is "automatic", as soon as we know that
$(M_t^{\downarrow,G})$ is a local martingale in the Brownian
filtration.)

\noi We have :
\begin{equation}\label{do11}
    M_t^{\downarrow,G}= \frac{G(D_{t})+G(1+D_{t})}{
2}+\Big(\frac{b-X_{t} }{ b-a}\Big) \Big(\frac{G(D_{t})- G(D_{t}+1)
}{ 2}\Big){\rm sign}(\sigma_{2D_{t}}+T_b \circ
\theta_{\sigma_{2D_{t}}} -t).
\end{equation}

\noi Since $t \to D_t$ and $t \to {\rm sgn}(\sigma_{2D_{t}}+ T_b
\circ \theta_{\sigma_{2D_{t}}}-t)$ are piecewise constant
processes, and $(M_t^{\downarrow,G})$ is continuous, applying
Ito's formula leads to (\ref{do9}). Hence $(M_t^{\downarrow,G})$
is a continuous local martingale.

\noi Relation (\ref{do11}) implies that :
$$|M_t^{\downarrow,G}| \le  \sup_{n \ge 0}|G(n)|
\Big(1+ \frac{|b|+X_{t}^{*} }{ b-a}\Big).$$

\noi Consequently, $(M_t^{\downarrow,G}\;;\; t\ge 0)$ is a
martingale.

\smallskip

\noi 2) Suppose that $G(n)>0, \, G(0)=1$, and $n\rightarrow G(n)$
decreases to $0$. If $\sigma_{2D_{t}}+T_b \circ
\theta_{\sigma_{2D_{t}}}
>t$, then $b-X_t \ge 0$ and formula (\ref{do11})  gives :
$$M_t^{\downarrow,G} \ge \frac{G(D_{t})+G(1+D_{t})}{ 2} > 0.$$

\noi If $\sigma_{2D_{t}} + \sigma_1 \circ \theta_{\sigma_{2D_{t}}}
\le t$, we first modify (\ref{do7}) as follows :
$$M_t^{\downarrow,G}=\Big(\frac{X_{t}-a }{ b-a} \Big)\frac {G(D_{t})-G(1+D_{t})
}{ 2}+G(1+D_t).$$

\noi But $X_t \ge a$, consequently $M_t^{\downarrow,G} \ge G(1+
D_t)>0$.

\noi Finally $M_t^{\downarrow,G} >0$.

\noi To end the proof of Proposition \ref{pdo1} we observe that
(\ref{do10}) is a direct consequence of (\ref{do9}) and
(\ref{do7}).
\end{prooff}

\vskip 1cm

\noi We state below our last theorem concerning penalization
results. Here the underlying process is  the number $D_t$ of
down-crossings from $b$ to $a$. Recall that the sequence of
stopping times $( \sigma_{n})_{n \ge 1}$ associated with $(D_t)$
is defined through (\ref{do1})-(\ref{do4}). Let $G \;:\; \Bbb N
\to \mathbb{R}_+$, be a decreasing function, such that :
\begin{equation}\label{do13}
    G(0)=1, \quad G(\infty):= \lim_{n \to \infty}
G(n)=0.
\end{equation}

\noi Let $\big(\Delta G(n)\big)_{n \ge 0}$ be the sequence of
positive numbers :
\begin{equation}\label{do14}
    \Delta G(n)=G(n)-G(n+1), \quad n \ge 0.
\end{equation}
%

\begin{theo} \label{tdo1}  Let $G:\mathbb{N}  \to \mathbb{R}_+$ as above.
\begin{enumerate}
    \item Let $s \ge 0, \; \Gamma_s \in {\cal F}_s$, then :
\begin{equation}\label{do15}
    \lim _{t \to \infty} \frac{E_0 \big[ 1_{\Gamma_s} \Delta
G(D_t)\big] }{ E_0 \big[\Delta G(D_{t})\big]} = \frac{1 }{
M_0^{\downarrow,G}} \, E_0 [1_{\Gamma_s} M_s^{\downarrow,G}],
\end{equation}

\noi  where $(M_t^{\downarrow,G})$ is the positive martingale
defined in Proposition \ref{pdo1}.

    \item  Let $Q_0^{\downarrow,G}$ be the p.m. defined on
$\big(\Omega, {\cal F}_\infty\big)$ :
\begin{equation}\label{do15a}
    Q_0^{\downarrow, G} (\Gamma_s)=\frac{1 }{
M_0^{\downarrow,G}} \, E_0 [1_{\Gamma_s} M_s^{\downarrow,G}],
\end{equation}

\noi for any $s \ge 0$ and $\Gamma_s \in {\cal F}_s$.

\noi Then the process :
\begin{equation}\label{do16}
    X_t-\int_0^t
\big(G(1+D_s)-G(D_s)\big)\biggr\{\frac{1 }{ G(D_s)(2b-X_s-a)+
G(1+D_s)(X_s-a)}
\end{equation}
$$\times 1_{\{\sigma_{2D_{s}}+T_{b}\circ
\theta_{2D_{s}}>s\}}- \frac{1 }{
G(1+D_{s})(2b-X_{s}-a)+G(D_{s})(X_{s}-a)}
1_{\{\sigma_{2D_{s}}+T_{b}\circ \theta_{2D_{s}}\leq s\}}
\biggl\}ds
$$
\noi is a $Q_0^{\downarrow, G}$ Brownian motion.
\end{enumerate}
\end{theo}

\noi The proof of Theorem \ref{tdo1} is based on the following
asymptotic estimate.

\begin{lemma}\label{ldo1}
Let $x \in \mathbb{R}$ and $\big(H(n)\big)_{n \ge 0}$ be a
sequence of positive real numbers satisfying : $\dis \sum_{n \ge
0} H(n) < \infty$. Then :
\begin{equation}\label{do17}
    \lim _{t \to \infty} \sqrt{t} \, E_x \big[H(D_t)\big]
= 2(b-a) \sqrt{\frac{2 }{ \pi}} \Big\{ \sum_{n \ge 1} H(n) + H(0)
\Big(\frac{1 }{ 2} + \frac{|x-b| }{ 2(b-a)}\Big) \Big\}.
\end{equation}
\end{lemma}

\begin{proof} \ On the one hand we observe that the definition
(\ref{do5}) of $D_t$ implies $\{D_t \ge n\}=\{\sigma_{2n} \le
t\}$.

\noi On the other hand, under $P_x,$  $ \sigma_1,
\sigma_2-\sigma_1, \cdots , \sigma_{2n}-\sigma_{2n-1}$ are
independent, $\sigma_1$ (resp. $\sigma_{i+1}-\sigma_{i})$ is
distributed as $T_{|b-x|}$ (resp. $T_{b-a}$) under $P_0$.
Consequently :
$$P_x(D_t \ge n) =P_0 (T_{|b-x|+(2n-1)(b-a)} \le t).$$
\noi Since the probability on the right hand-side equals $P_0
\big(S_t \ge |b-x|+(2n-1)(b-a)\big)$, the scaling property
(\ref{osm10}) implies :
\begin{equation}\label{do18}
    P_x (D_t \geq n)=P_0\Big(n \le \frac {\sqrt{t}|X_1|+b-a-|b-x| }{ 2(b-a)}\Big)
\;;\; n \ge 1.
\end{equation}
\noi $D_t$ being a $\mathbb{N}$-valued r.v., we have :
$$E_x \big[H(D_t)\big]=H(0) + \sum_{n \ge 1} \big(H(n)-H(n-1)\big)
P_x (D_t \ge n)$$

\noi Let $\xi_t:=\dis\frac {\sqrt{t}|X_{1}|+b-a-|b-x| }{2(b-a)}$.
 Using Fubini's theorem and (\ref{do18}) we get
$$
\begin{array}{ccl}
  \dis E_x\big[H(D_t)\big] & = & \dis H(0)-E_0\Big[1_{\{\xi_{t}\ge1\}}
\sum_{n=1}^{[\xi_{t}]} \Big( H(n)-H(n-1)\Big)\Big] \\
   & = & H(0) P_0(\xi_t<1)+E_0 \Big[H\big([\xi_t]\big) 1_{\{\xi_{t} \ge
1\}}\Big].\\
\end{array}
  $$

\noi where $[a]$ denotes the integer part of $a$.

\noi It is easy to compute the last expectation :
$$
E_0 \big[H\big([\xi_t]\big)1_{\{\xi_t \ge 1\}}\big] =
\sqrt{\frac{2 }{ \pi}} \int_0^{\infty} H \Big(\Big[\frac{z
\sqrt{t} +b-a-|b-x] }{ 2(b-a)}\Big]\Big) e^{-z^2/2} 1_{\big\{{ z
\sqrt{t}+b-a-|b-x| \over 2 (b-a)}>1\big\}} dz.$$

\noi Setting $y=\dis \frac{z\sqrt{t}+b-a-|b-x| }{ 2(b-a)}$, we
have :
$$
E_0\big[H\big([\xi_t]\big) 1_{\{\xi_t \ge 1\}}\big] = 2(b-a)
\sqrt{\frac{2 }{ \pi t}} \int_1^{\infty} H\big([y]\big)  \exp
\big\{-\frac{\big(2(b-a)y-b+a+|b-x|\big)^2 }{ 2t}\big\} dy.$$

\noi Since $\dis \sum_{n \ge 1} H(n)< \infty$,

$$ \lim_{t \to \infty} \Big(\sqrt{t}\, E_0\big[H([\xi_t])
1_{\{\xi_{t} \ge 1\}}\big]\Big)=2(b-a) \sqrt{\frac{2 }{ \pi}
}\biggl( \sum_{n \ge 1} H(n)\biggr).$$

\noi It remains to study $P_0 (\xi_t < 1)$. We have successively :
$$
 P_0 (\xi_t < 1) = P_0 \Big(|X_1| < \frac{b-a + |b-x| }{
\sqrt{t}}\Big), $$
\noi hence :

$$  \lim_{t \to \infty}  \sqrt{t} \, P(\xi_t <
1) = \big(b-a+|b-x|\big)\, \sqrt{\frac{2 }{ \pi} }.$$

\noi This ends the proof of (\ref{do17}).
\end{proof}

\begin{prooff} \ {\bf of  Theorem \ref{tdo1}}

\noi Let $s \ge 0, \, \Gamma_s \in {\cal F}_s$ be fixed, $t>s$ and
consider the quantity $E_0\big[1_{\Gamma_s} \Delta G (D_t)\big].$

\noi Let us introduce the  events :
$$\Sigma_1=\big\{\sigma_{2D_{s}}+T_b \circ \theta_{\sigma_{2D_{s}}}
> s\big\}, \; \Sigma_2=\Sigma_1^c=\big\{\sigma_{2D_{s}}+T_b \circ
\theta_{\sigma_{2D_{s}}} \le s\big\} ,$$
\noi and decompose the above expectation accordingly.

 \noi  1) On
$\Sigma_1$ :
\begin{equation}\label{do19}
    D_t=D_{t-s} \circ \theta_s +D_s.
\end{equation}

\noi Applying Lemma \ref{ldo1} with $x=X_s, \, H(n)=\Delta
G(n+D_s)$, after conditioning by ${\cal F}_s$,
 we obtain :
\begin{equation}\label{do20}
    E_0\big[1_{\Gamma _s\cap \Sigma_{1}} \Delta G (D_t)\big] _{\stackrel{\sim}{t \to
    \infty}}
2(b-a) \sqrt{\frac{2 }{ \pi}} E_0 [1_{\Gamma_s \cap \Sigma_{1}}
J_1] \frac{1 }{ \sqrt{t-s}},
\end{equation}

\noi where $J_1=\dis \Delta G (D_s) \, \Big(\frac{1
}{2}+\frac{|X_{s}-b|}{ 2(b-a)}\Big) +\sum_{n \ge 1} \Delta G
(D_s+n)$.

\noi Since on $\Sigma_1, \; X_s \le b$, it is easy to check that :
$$J_1=\frac{G(D_s)}{ 2}\Big(1+\frac{b-X_s}{ b-a}\Big)+
\frac{G(1+D_{s}) }{ 2}\Big(\frac{X_{s}-a }{ b-a}\Big) .$$
\noi As a result :
$$E_0\big[1_{\Gamma_s \cap \Sigma_{1}} \Delta G (D_t)\big] _{\stackrel{\sim}{t \to
    \infty}}\Big\{ 2(b-a) \sqrt{\frac{2 }{ \pi}} E_0 [1_{\Gamma_s \cap
    \Sigma_{1}}M_s^{\downarrow ,G}]\Big\}\frac{1}{\sqrt{t}}.$$

\noi 2)  We decompose $\Sigma_2$ in two disjoint events :
$\Sigma_2 = \Sigma'_2\cup\Sigma_3$ where $\Sigma_2'=\Sigma_2 \cap
\{X_s > b\}$ and $ \Sigma_3= \Sigma_2 \cap \{a <X_s \le b\}$.

\noi On $\Sigma_2 \cap \{X_s > b\}$, (\ref{do19}) holds and as
previously :
$$E_x[1_{\Gamma_s \cap \Sigma_{2} \cap \{X_s > b\}}\Delta G (D_t)] _{\stackrel{\sim}{t \to
    \infty}} \Big(2(b-a) \sqrt{\frac{2 }{ \pi}} E_x [1_{\Gamma _s\cap \Sigma_{2}
\cap \{X_{s} > b\}} J_1]\Big) \frac{1 }{ \sqrt{t-s}} \cdot$$

\noi Since $X_s > b, \; J_1=\dis \frac{G(1+D_{s})}{ 2}  \Big(
1+\frac{b-X_{s} }{ b-a}\Big)+\frac{G(D_s) }{ 2} \Big(\frac{X_s-a
}{ b-a}\Big) $ and
$$E_0\big[1_{\Gamma _s\cap \Sigma_{2} \cap \{X_s > b\}} \Delta G (D_t)\big] _{\stackrel{\sim}{t \to
    \infty}}\Big\{ 2(b-a) \sqrt{\frac{2 }{ \pi}} E_0 [1_{\Gamma _s\cap \Sigma_{2} \cap \{X_s > b\}}
    M_s^{\downarrow ,G}]\Big\}\frac{1}{\sqrt{t}}.$$

\noi 3) We split $\Sigma_3$ in three disjoint subsets :
$\Sigma_3=\Sigma_4 \cup \Sigma_5\cup \Sigma_6$, with
$$
\begin{array}{ccl}
 \Sigma_4 & = & \Sigma_3 \cap \{s+T_a \circ \theta_s<s+T_b \circ \theta_s <t\}, \\
  \Sigma_5 & = & \Sigma_3 \cap \{s+T_b \circ \theta_s <s+T_a \circ \theta_s<t\}, \\
  \Sigma_6 & = & \Sigma_4 \cap \{(s+T_b \circ \theta_s) \wedge (s+T_a \circ
\theta_s) \ge t\}. \\
\end{array}
$$

\noi We set $U_a=s+T_a \circ \theta_s$ and $U_b=s+T_b \circ
\theta_s$.

\smallskip

\noi a)  On $\Sigma_4$, we have :
$$D_t=D_{U_{a}}+D_{t-U_{a}} \circ \theta_{U_{a}}=1+D_s+
D_{t-U_{a}} \circ \theta_{U_{a}}.$$

\noi Applying the strong Markov property at time $U_a$, together
with Lemma \ref{ldo1}, we obtain on $\Sigma_4$ :
$$\lim_{t \to \infty} \sqrt{t}\, E_0 \big[\Delta G
(D_t)|{\cal F}_{U_{a}}\big]=2(b-a) \sqrt{\frac{2 }{\pi}}
G(D_{U_{a}}).$$
\noi But $D_{U_{a}}=1+D_s$; hence, taking the conditional
expectation with respect to ${\cal F}_s$, we get :
$$ \lim_{t \to \infty} \sqrt{t} \, E_0 \big[1_{\Gamma_s
\cap \Sigma_{4}} \Delta G (D_t)\big] = 2(b-a) \sqrt{\frac{2 }{
\pi}}
 E_0\big[1_{\Gamma _s\cap \Sigma_{3}} \frac{b-X_s }{ b-a}
G(1+D_s)\big].$$

\noi b) On $\Sigma_5, \; D_t=D_{U_{b}}+D_{t-U_{b}} \circ
\theta_{U_{b}}=D_s+D_{t-U_{b}} \circ \theta_{U_{b}}$. We proceed
as previously; we obtain  successively  :
$$
 \lim_{t \to \infty} \sqrt{t} E_0 \big[1_{\Sigma_{5}}
\Delta G(D_t) |{\cal F}_{U_{b}}\big]= 2(b-a) \sqrt{\frac{2}{ \pi}}
\frac {G(D_{U_{b}})+G(D_{U_{b}}+1) }{ 2}\,1_{\{U_b<U_a\}},$$

$$
 \lim_{t \to \infty} \sqrt{t} \, E_0 \big[1_{\Gamma_s
\cap\Sigma_{5}} \Delta G(D_t)\big]=2(b-a) \sqrt{\frac{2 }{ \pi}}
 E_0 \Big[1_{\Gamma_s \cap \Sigma_{3}} \frac{X_s-a }{ b-a}\frac{G(D_s)+
G(1+D_{s}) }{ 2}\Big] . $$

\noi c) We claim that $\Sigma_6$ does not contribute to the limit
since
$$
\begin{array}{ccl}
  P_0(\Sigma_6|{\cal F}_s) & = & 1_{\Sigma_{3}} P_{X_s}(T_a \wedge  T_b >t-s) \le
1_{\Sigma_{3}} P_{\frac{a+b} { 2}} (T_a \wedge T_b>t-s) \\
   & \leq & 1_{\Sigma_{3}} \, e^{-\lambda (t-s)}, \\
\end{array}
$$
\noi for some $\lambda > 0$.

\noi Finally :
$$E_0\big[1_{\Gamma_s \cap \Sigma_{2} \cap \{X_s \leq  b\}} \Delta G (D_t)\big] _{\stackrel{\sim}{t \to
    \infty}}\Big\{ 2(b-a) \sqrt{\frac{2 }{ \pi}} E_0 [1_{\Gamma _s\cap \Sigma_{2} \cap \{X_s \leq b\}}
    M_s^{\downarrow ,G}]\Big\}\frac{1}{\sqrt{t}}.$$

\noi 4) Consequently, thanks  to the previous steps we have :
$$E_0\big[1_{\Gamma _s} \Delta G (D_t)\big] _{\stackrel{\sim}{t \to
    \infty}}\Big\{ 2(b-a) \sqrt{\frac{2 }{ \pi}} E_0 [1_{\Gamma_s }
    M_s^{\downarrow ,G}]\Big\}\frac{1}{\sqrt{t}}.$$
This implies (\ref{do15}).

\noi Point 2. of Theorem \ref{tdo1} is a direct consequence of
Proposition \ref{pdo1}.
\end{prooff}

\section{Study of the $Q$-processes. An approach via
enlargements of filtrations} \label{sQ}
%
\setcounter{equation}{0}

 \noi In sections \ref{osm}-\ref{do}  we
have obtained penalization principles involving unilateral
maximum, unilateral maximum and time, local time at $0$, maximum +
minimum + local time at $0$ and finally down-crossings. In each
case, a positive and continuous martingale $(M_t \;;\; t \ge 0)$
appears naturally. This allows to define $Q_x$ on $\big(\Omega,
{\cal F}_\infty\big)$ via : $Q_x(\Gamma_t)=E_x [1_{\Gamma_t}
M_t]$, for any $\Gamma_t \in {\cal F}_t, \; t \ge 0$.

\noi As said in the Introduction, $Q_x$ is a well-defined p.m. on
$(\Omega, {\cal F}_\infty)$.

\noi This leads us to describe the law of $(X_t)_{t \ge 0}$ under
$Q_x, \; x \in \mathbb{R}$.

\noi We are able to handle the two first cases recalled previously
via a general approach, which is developed in subsection
\ref{ssQ}. Unfortunately the other cases cannot be handled in this
way, and we have to study them one by one. Our approach then is
based on enlargements of filtrations. However the schemes of proof
are similar in all cases, consequently we only discuss Case 1 in
details . Concerning the other cases, we  only state the results
and sketch their proofs giving the key points without detailed
arguments.

\subsection{ Some general results and their applications}\label{ssQ}


\noi We consider here  a general setting, where we are given a
filtered probability space $\big(\Omega, {\cal F}, ({\cal F}_t)_{t
\geq 0}, \mathbb{P}\big)$, a strictly positive continuous
martingale $(M_t;t\geq 0)$, with respect to $\big(({\cal F}_t)_{t
\geq 0}, \mathbb{P}\big)$, starting at $1$ at $t=0$, and a second
probability $\mathbb{Q}$ on $(\Omega,{\cal F}_\infty)$ such that :

%
\begin{equation}\label{Q1a}
    \mathbb{Q}(\Gamma_t)=\mathbb{E} [1_{\Gamma_t} M_t], \; \Gamma_t \in {\cal F}_t, \; t \ge 0.
\end{equation}
\noi  We define :
\begin{equation}\label{Q1}
    \underline{M}_t=\inf_{0 \le s \le t} M_s.
\end{equation}
The following discussion  is inspired from \cite{AY3}

\begin{prop} \label{pQ1}\begin{enumerate}
\item The process $\displaystyle (Y_t=\frac{M_t}{\underline{M}_t}
-1,\; t \ge 0)$ is a non-negative, continuous, local
$\mathbb{P}$-submartingale and $\mathbb{P}(Y_0=0)=1$.

    \item Let $(l_t \;;\; t \ge 0)$ be the non-decreasing,
continuous process, such that $l_0=0$ and $(Y_t-l_t)_{t \ge 0}$ is
a continuous local martingale. Then the support of $dl_t$ is
included in $\{t \ge 0 \;;\; Y_t=0\}$ and
\begin{equation}\label{Q2}
    \underline{M}_t= e^{-l_{t}} \;;\; t \ge 0.
\end{equation}
\end{enumerate}
\end{prop}

 \begin{proof}\  Using Ito's formula we have :
$$dY_t = {dM_t \over \underline{M}_t}-{M_t \over \underline{M}_t^2}
\, d \underline{M}_t .$$

\noi Consequently $(Y_t)$ is a local submartingale and :
\begin{equation}\label{Q3}
    dl_t=-\frac{M_t }{ \underline{M}_t^2} \, d\underline{M}_t=-
\frac{1 }{ \underline{M}_t}  d\underline{M}_t .
\end{equation}

\noi This implies that $supp(dl_t) = supp (d\underline{M_t})$.

\noi Integrating (\ref{Q3}) over $[0,t]$ leads to (\ref{Q2}).
\end{proof}

 \noi  \noi To go further, we need
an additional assumption; i.e.,we assume :
\begin{equation}\label{Q5}
M_{\infty}=0 \quad \mathbb{P} \   a.s.
\end{equation}

 It is clear that $M_{\infty}=0$
iff $\underline{M}_{\infty}=0$. From (\ref{Q2}) this condition is
equivalent to $l_{\infty}=+\infty$.

\begin{theo}\label{tQ1} Let $\mathbb{Q}$ be the p.m. defined by
(\ref{Q1a}).
\begin{enumerate}
    \item $\underline{M}_{\infty}$ is a $\mathbb{Q}$-finite r.v. with
uniform distribution on $[0,1]$.
    \item Let $g:={\rm sup}\{t \ge 0, \; M_t=\underline{M}_{\infty}\}$ (with
the convention $\sup \emptyset=0$). Then $\mathbb{Q}(0<g<\infty)
=1$. Let
\begin{equation}\label{Q6}
    Z_t=\mathbb{Q}(g>t|{\cal F}_t) \;;\; t \ge 0.
\end{equation}
    \noi Then
    \begin{enumerate}
        \item $\dis Z_t =\underline{M}_t / M_t$ ,
        \item $(Z_t)$ is a positive, $\mathbb{Q}$-supermartingale with
additive decomposition :
\begin{equation}\label{Q7}
   Z _t=1-\int_0^t \frac{\underline{M}_u
}{ M_u^2}\,d\widetilde{M}_ u +\ln(\underline{M}_t),
\end{equation}
\noi where $\dis \widetilde{M}_t=M_t-\int_0 ^t \frac{d<M>_u}{M_u}$
is the martingale part of $(M_t)$ under $\mathbb{Q}$.
    \end{enumerate}
\end{enumerate}

\end{theo}

\begin{prooff}\ {\bf of Theorem \ref{tQ1}}
i) Let us determine the distribution function of
$\underline{M}_{\infty}$ under $\mathbb{Q}$. Let $t>0$ and
$0<c<1$. We have :
$$\mathbb{Q}(\underline{M}_t <c)=\mathbb{Q} (\sigma(c) < t)=\mathbb{E}_x [1_{\{
\sigma (c)<t\}} M_t],$$

\noi where $\sigma (c)={\rm inf}\{u \ge 0, \; M_u<c\}$.

\noi Applying Doob's optional stopping theorem we obtain :
$$\mathbb{Q}(\underline{M}_t<c)=\mathbb{E} [1_{\{\sigma (c) <t\}} M_{\sigma (c)}
]=c \, \mathbb{P}(\sigma (c) <t).$$

\noi Taking $t \to \infty$, we obtain :
$\mathbb{Q}(\underline{M}_{\infty} <c)=c$.

\smallskip

\noi  ii) Let us compute $E_{\mathbb{Q}}[g>t|{\cal F}_t]$, where
$t>0$ is fixed. Let $\Gamma_t$ be in ${\cal F}_t$. We have :
$$\mathbb{Q}\big(\Gamma_t \cap \{g> t\}\big)=\mathbb{Q}\big(\Gamma_t \cap
\{\sigma'_t < \infty\}\big),$$
\noi where $\sigma'_t={\rm inf}\{s>t\;;\;M_s \leq
\underline{M}_t\}$.

\noi Consequently :
$$\mathbb{Q}\big(\Gamma_t \cap \{g>t\}\big)= \lim_{n \to \infty}
\delta_n,$$
\noi with $\delta_{n}=\mathbb{Q}\big(\Gamma_t \cap \{\sigma'_t \le
t+n\}\big)$.

\noi Using the same technique as in step i), we have successively
:
$$
\begin{array}{ccl}
  \delta_n  & = & \displaystyle \mathbb{E}[1_{\Gamma_t \cap \{\sigma'_t \le t+n\}} M_{t+n}]
  =\mathbb{E}[1_{\Gamma_t \cap
\{\sigma'_t \le t+n\}} M_{\sigma'_t}] \\
   & = & \displaystyle \mathbb{E}[1_{\Gamma_t \cap \{\sigma'_t \le t+n\}} \underline{M}_t]. \\
\end{array}
$$

\noi Letting $n \to \infty$, we get :
\begin{equation}\label{Q8}
    \mathbb{Q}\big(\Gamma_t \cap \{g>t\}\big)=\mathbb{E}[1_{\Gamma_t}\,
     \underline{M}_t]=E_{\mathbb{Q}}
\Big[1_{\Gamma_t} \underline{M}_t / M_t\Big].
\end{equation}
%



\noi This proves that $Z_t=\underline{M}_t/M_t$ (under
$\mathbb{Q}$).

\noi Equivalently to (\ref{Q8}), we have :
$$\mathbb{Q}\big(\Gamma_t\cap \{0<g \leq t\}\big)=
\mathbb{Q}(\Gamma_t)- \mathbb{E}[1_{\Gamma_t}
\,\underline{M}_t],$$
from which it follows that $\mathbb{Q}(0<g <\infty)=1$.
\smallskip

\noi  iii) Since $t \to 1_{\{g>t\}}$ is non-increasing, then
$(Z_t)$ is a $\mathbb{Q}$-supermartingale. Applying Ito's formula
we get :
$$dZ_t={d\underline{M}_t \over M_t}-\frac{\underline{M}_t }{
M_t^2}\,dM_t+\frac{\underline{M}_t }{ M_t^3}d<M>_t.$$
\noi The decomposition (\ref{Q7}) is a direct consequence of
Girsanov's theorem.

\smallskip

\end{prooff}

\begin{rem}\label{rQ1b}  We write :  $\mathbb{P} \longrightarrow \mathbb{Q}$
in case the pair $(\mathbb{P},\mathbb{Q} )$ satisfies the absolute
continuity relationship (\ref{Q1a}), as well as (\ref{Q5}).
We claim that this relation is symmetric, i.e. if $\mathbb{P}
\longrightarrow \mathbb{Q}$, then $\mathbb{Q} \longrightarrow
\mathbb{P}$.

\noi Indeed, let  $N_t=1/M_t, \ t\geq 0$. It is clear that $(N_t)$
is a positive and continuous $\mathbb{Q}$ martingale, starting at
$1$. To prove the claim, we have to check that $\dis
\lim_{t\rightarrow\infty}N_t=0, \ \mathbb{Q}$ a.s.

\noi Let $A>1$ be a real number and $\sigma_A=\inf\{t\geq 0;
M_t\geq A\}$. Recall (\cite{RevYor}, Ex 3.12, Chap. II) that $\dis
\sup_{t\geq 0}M_t$ is distributed as $1/U$ where U is uniformly
distributed on $[0,1]$.

\noi Let $t$ be a fixed real number. We have :
$$\mathbb{Q}(\sigma_A<t)=\mathbb{E}[1_{\{\sigma_A<t\}}M_t]=
\mathbb{E}[1_{\{\sigma_A<t\}}M_{\sigma_A}]=
A\mathbb{P}(\sigma_A<t).$$
Taking $t\rightarrow\infty$, we get :
$$\mathbb{Q}(\sigma_A<\infty)=A\mathbb{P}(\sigma_A<\infty)=A\mathbb{P}(\sup_{t\geq 0}M_t
>A)=A\frac{1}{A}=1.$$

\noi Since $A$ is arbitrary, this implies that $\dis \sup_{t\geq
0}M_t=+\infty$, $\mathbb{Q}$ a.s. Consequently, under
$\mathbb{Q}$, $\dis \lim_{t\rightarrow \infty}N_t=0$.

\noi
\end{rem}

\noi Applying Remark \ref{rQ1b} with Theorem \ref{tQ1} we obtain
the following.

\begin{theo}\label{tQ1b} Let $(N_t;t\geq 0)$ be a positive and
continuous $\mathbb{Q}$ martingale such that under $\mathbb{Q}$ :
$N_0=1$ and $\dis \lim_{t\rightarrow \infty}N_t=0$.

\begin{enumerate}
    \item $\dis\sup_{t\geq 0}N_t\stackrel{(d)}{=}1/U,$  where U is uniformly
distributed on $[0,1]$.
    \item If $\dis g=\sup\{t\geq 0; N_t=\sup_{u\geq 0}N_u\}$, then :
      $\mathbb{Q}(0<g<\infty) =1$. Let $Z_t=\mathbb{Q}(g>t|{\cal F}_t) , t \ge
      0$.  Then
    \begin{enumerate}
        \item $\dis Z_t =\frac{N_t}{\overline{N}_t}$ , where $\dis \overline{N}_t:=\sup_{0\leq u\leq
        t}N_u$.
        \item $(Z_t)$ is a positive, $\mathbb{Q}$-supermartingale, with   Doob-Meyer  decomposition
        : $Z_t=M'_t-\ln(\overline{N}_t)$, where $(M'_t)$ denotes a $\mathbb{Q}$
        martingale.

\end{enumerate}

\end{enumerate}
\end{theo}

\noi We now make the further assumption that $\mathbb{P}=P_x$ is
the Wiener measure on the canonical space ${\cal C} (\mathbb{R}_+,
\mathbb{R})$, and $P_x(X_0=x)=1$. We shall also write $Q_x$ for
$\mathbb{Q}$, in this particular case. We now gather a number of
complements to our previous general results in this particular
instance :
\begin{enumerate}
    \item Since $(M_t)_{t \geq 0}$ is a $P_x$-martingale and $M_0=1$,
then the representation theorem of Brownian martingales
(\cite{RevYor}, section V.3, p.192) implies
 that $(M_t)$ may be written as :
\begin{equation}\label{Q4}
    M_t=1+ \int_0^t m_s\,dX_s,
\end{equation}
\noi where $(m_s)$ is a predictable process, such that for any $t
\ge 0$, $\dis \int_0^t m_s^2 ds<+\infty$ a.s.

    \item  The process $\big(\beta_t=X_t-\dis\int_0^t \frac{m_u }{
M_u}du ; t \ge 0\Big)$ is a $Q_x$-Brownian motion started at $x$.
\end{enumerate}

\begin{rem}\label{rQ1}  For the sake of efficiency, we now use
the technique of progressive enlargement of filtrations (see for
instance \cite{J1} or \cite{Y1}).
Let $({\cal G}_t)$ be the smallest filtration containing $({\cal
F}_t)$ and such that $g$ is a $({\cal G}_t)$- stopping time. Then
:
$$\beta_t=\widetilde{\beta}_t+\int_0^{t \wedge g}
\frac{d<Z,\beta>_u }{ Z_u}-\int_{t \wedge g}^t \frac{d<Z, \beta>_u
}{ 1-Z_u}\;; \; t\ge 0,$$
\noi where $(\widetilde{\beta}_t\;;\;t \ge 0)$ is a $\big( ({\cal
G}_t), \, Q_x\big)$-Brownian motion started at $x$.

\noi Since $Z_t=\underline{M}_t/M_t$, applying directly (\ref{Q7})
leads to :
$$\beta_t=\widetilde{\beta}_t - \int_0^{t \wedge g} \frac{m_u }
{M_u} du + \int_{t \wedge g}^t \frac{\underline{M}_{u} m_u }{M_u
(M_u-\underline{M}_u)}du. $$

\noi Recall that  $X_t=\beta_t+\dis \int_0^t \frac{m_u }{M_u}du$,
consequently :
\begin{equation}\label{Q9}
    X_t=\widetilde{\beta}_t+ \int_{t \wedge g}^t \frac{m_u }{ M_u-
\underline{M}_u}du .
\end{equation}
\noi In particular $X_t=\widetilde{\beta}_t$ for any $t\le g$.
 Let $\sigma (y)={\rm inf}\{t \ge 0, \; M_t < y\}$,
for any $0<y<1$. Then conditionally on $\underline{M}_{\infty}
=y$, the process $(X_t ; 0 \le t \le g)$ is distributed as a
Brownian motion started at $x$ and considered up to its first
hitting time of level  $y$.

\noi Going back to (\ref{Q9}), we note that, for $u>g$, the drift
term equals
$${m_u \over M_u - \underline{M}_{\infty}}\cdot,$$
and, in general, $(X_t)$ is not a diffusion, except in particular
cases, see below Theorems \ref{tQ2} and \ref{tQ3}.

\end{rem}

\noi Our first application of Theorem \ref{tQ1} concerns
$Q_x^{\varphi}$. Recall that
$$Q_x^{\varphi}
(\Gamma_t)=\frac{1}{1-\Phi(x)}E_x[1_{\Gamma_t}\,M_t^{\varphi}],$$
for any $t>0$ and $\Gamma_t \in{\cal F}_t$, where
$M_t^{\varphi}=(S_t-X_t) \,\varphi(S_t)+1-\Phi(S_t)$, and  $
\varphi$ satisfies the conditions given in Proposition \ref{posm1}
and (\ref{osm3}).

\begin{theo}\label{tQ2}\begin{enumerate}
    \item Under $Q_x^{\varphi}$, the r.v. $S_{\infty}$ is
positive, finite and admits $\dis {\varphi (y) \over 1-\Phi(x)}
1_{\{y \ge x\}}$ as a density function.
    \item  Let $g={\rm sup}\{u \ge 0 \;;\;
X_u=S_{\infty}\}$. Then $Q_x^{\varphi}(0<g<\infty)=1$, and under
$Q_x^{\varphi}$ :

\begin{enumerate}
    \item the processes $(X_u ;u \le g)$ and $(X_{g}-X_{u+g} ;u \ge 0)$
are independent,

    \item $(X_g - X_{u+g} ; u\ge 0)$ is distributed as a three
dimensional Bessel process started at $0$,
    \item conditionally on $S_{\infty}=z>x,\;(X_u;u \le g)$ is
distributed as a Brownian motion started at $x$ and stopped at its
first hitting time of $z$.
\end{enumerate}
\end{enumerate}

\end{theo}

\begin{proof} \ a) Let $(a,b)$ be a maximal  interval of excursion of $(X_t)$ below its unilateral
supremum $(S_t)$. Since for any $u \in ]a,b [$, $S_u-X_u>0$ and
$S_u= S_a=S_b$ , then $\dis \inf_{0\le u \le t}
M_u^{\varphi}=1-\Phi(S_t)$. Let $M_t=\dis \frac{1 }{ 1-\Phi(x)}
M_t^{\varphi}$. Then $Q_x^\varphi(M_0=1)=1$ and $M_{\infty}=
\inf_{u \ge 0} M_u=\dis\frac{1-\Phi(S_{\infty})}{ 1-\Phi(x)}$.

\noi From Proposition \ref{posm1}, decomposition (\ref{Q4}) holds
with $m_t=\dis \frac{-\varphi (S_t) }{ 1-\Phi(x)}$. Then
assumption (\ref{osm3}) implies that $P(M_\infty =0)=1$.

 \noi Consequently 1) and 2) c) are   direct consequences  of Theorem
\ref{tQ1}.

\smallskip

\noi  b)  For any $u \ge g$,  we have : $m_u=\dis- \frac{\varphi
(S_{\infty}) }{ 1-\Phi(x)}, \; M_u=\frac{1 }{ 1- \Phi (x)}
\big[(S_{\infty}-X_u)\,\varphi(S_{\infty})+1-\Phi
(S_{\infty})\big]$, $\underline{M}_u=\underline{M}_\infty=\dis
\frac {1-\Phi(S_{\infty}) }{ 1-\Phi (x)}$ and $X_g=S_\infty$.
Setting $R_u=X_g-X_{g+u}$, then the identity (\ref{Q9}) implies
that :
\begin{equation}\label{Q11}
    R_t=\widetilde{\beta}_g-\widetilde{\beta}_{t+g} + \int_0^t
\frac{du }{ R_u} \;;\; t \ge 0.
\end{equation}

\noi Recall that $g$ is a $({\cal G}_t)$ stopping time; hence,
$(\widetilde{\beta}_g-\widetilde{\beta}_{t+g} \;;\; t \ge 0)$ is a
Brownian motion,  independent from ${\cal G}_g$. Points a) and b)
are due to the fact that (\ref{Q11}) has a unique strong solution
whose distribution is the law of a three dimensional Bessel
process started at $0$.

\end{proof}

\noi We now describe the law of $(X_t)$ under $Q^{\lambda
,\varphi}_0$ (we take the starting point to be $0$, for
simplicity).


\begin{theo}\label{tQ2a}
 Let $\varphi, \Phi $ be the functions defined by (\ref{osm6b3}), resp. (\ref{osm6b2}),
and parameterized by the function $\psi$ satisfying
(\ref{osm6b1}).
\begin{enumerate}
    \item Under $Q_0^{\lambda, \varphi}$, the r.v.
$S_{\infty}$ is finite with density function $e^{- \lambda x}
\big(\varphi (x)+ \lambda (1- \Phi(x)\big)=e^{- \lambda x} \psi
(x)$.
    \item Let $g=\sup\{t \ge 0, \; X_t=S_{\infty}\}$.
Then $Q_0^{\lambda, \varphi} (0<g<\infty)=1$ and under
$Q_0^{\lambda, \varphi}$ :

\begin{enumerate}
\item $(S_{\infty}-X_{t+g} ; t \ge 0)$ is independent of $(X_t ; 0
\leq t \le g)$ and is distributed as $(Z_t^{(\lambda)} ; t \ge
0)$, under $P_0$, where
\begin{equation}\label{dir2}
    Z_t^{(\lambda)} =X_t + \lambda \int_0^t {\rm coth}(\lambda
\, Z_u^{(\lambda)}) du.
\end{equation}
    \item Conditionally on $S_{\infty}=x, \; (X_t ; t \le g)$ is
distributed as a Brownian motion with drift $\lambda$ started at
$0$, and stopped when it reaches $x$.
\end{enumerate}
\end{enumerate}

\end{theo}

\noi One proof of Theorem \ref{tQ2a} may be  based on the theory
of enlargements of filtration. This proof  is similar to the proof
of Theorem \ref{tQ1}. Therefore we do not give it. However we will
prove Theorem \ref{tQ2a}, using a   direct approach; see
Proposition \ref{pdir2} in Section \ref{dir}.

 \noi We are now able to deal with
$Q_0^{h^{+},h^{-}}$, where $Q_0^{h^{+},h^{-}}
(\Gamma_t)=E_0[1_{\Gamma_t} \, M_t^{h^{+},h^{-}}]\;;\; \Gamma_t
\in {\cal F}_t$ and $M_t^{h^{+},h^{-}}=X_t^+\,h^+(L_t^0)+
X_t\,h^{-}(L_t^0)+1-H(L_t^0)\;;\;t \ge 0$. We suppose that
$h^{+},h^{-}$ satisfy the conditions given in Proposition
\ref{ploc1}, and (\ref{loc4}). Recall that the function $H$ is
given by (\ref{loc1}).

\begin{theo}\label{tQ3}
\begin{enumerate}
    \item Under $Q_0^{h^{+},h^{-}}$, $L_{\infty}^0$ is a positive,
finite r.v. with density function : $\dis {1 \over 2}(h^+ + h^-)$
.

    \item Let $g={\rm sup}\{u \ge 0\;;\;X_u=0\}$. Then $Q_0^{h^{+},h^{-}} (0<g<\infty)=1$
and under $Q_0^{h^{+},h^{-}}$ :

\begin{enumerate}
    \item the processes $(X_u\;;\;u \le g)$ and $(X_{u+g} \;;\; u \ge 0)$
are independent,
    \item  with probability $\dis \frac{1 }{ 2} \int_0^{\infty} h^+(u)du$
(resp. $\dis \frac{1 }{ 2} \int_0^{\infty}h^- (u) du$), the
process $(X_{u+g} ;u \ge 0)$ (resp. $(-X_{u+g} ; u \ge 0)$) is
distributed as a three dimensional Bessel process, started at $0$.
    \item conditionally on $L_{\infty}^0=l,\,(X_u\;;\;u \le g)$ is
distributed as a Brownian motion started at $0$ and stopped when
its local time at $0$ equals $l$.
\end{enumerate}
\end{enumerate}

\end{theo}

\begin{proof} \  The proofs of 1), 2) a) and c) are similar to the proof of
Theorem \ref {tQ2}, being based on the technique of enlargements
of filtrations. They  are  left to the reader.

\noi However one new  point has to be checked  :
$Q_0^{h^{+},h^{-}} (\Gamma_+)= \dis \frac{1 }{ 2} \int_0^{\infty}
h^+ (u)du$, where $\Gamma_+$ is the set : $\Gamma_+=\{ X_u>0 ;
\forall u>g\}$.

\noi Since either $X_u>0 , \forall u > g$ or $X_u<0 , \forall u
>g$, then  $\dis Q_0^{h^{+},h^{-}}(\Gamma_+) =  \lim_{t \to \infty}
Q_0^{h^{+},h^{-}}(X_t>0)$.

\noi It is easy to compute $Q_0^{h^{+},h^{-}}(X_t>0)$ :
$$\begin{array}{ccl}
  Q_0^{h^{+},h^{-}}(X_t>0) & = & \dis E_0 [1_{\{X_t >0\}} M_t^{h^{+},h^{-}}] \\
   & = & \dis E_0\big[X_t^+ \, h^{+}(L_t^0)\big] +E_0 \big[1_{\{X_t>0\}}
\big(1-H(L_t^0)\big)\big]. \\
\end{array}
$$

\noi The dominated convergence theorem implies that $ \dis \lim_{t
\to \infty}E_0 \big[1_{\{X_t >0\}} \big(1-H (L_t^0)\big)\big]=0$.

\noi Applying Proposition \ref{ploc1}  with $h^-=0$, we get :
$$E_0\big[X_t^+\,h^{+}(L_t^0)\big]=\frac{1 }{ 2}E_0\big[
\int_0^{L_{t}^{0}} h^+ (u)du\big].$$

\noi Therefore  : $Q_0^{h^{+},h^{-}}(\Gamma_+)=\dis \frac{1 }{ 2}
\int_0^{\infty} h^+ (u)du$.
\end{proof}


\subsection{The $Q$-process associated with the bilateral supremum
and local time}\label{ssQ2}


\noi In this sub-section we study the law of $(X_t)$ under
$Q_0^{{\nu}_{*}}$, where $\nu_*$ is a p.m. on $[\alpha ,\infty[$,
for some $\alpha >0$. Recall that  $Q_0^{{\nu}_{*}}$ is the p.m.
on $\big(\Omega , {\cal F}_\infty\big) \;:\; Q_0^{{\nu}_{*}}
(\Gamma_t)=E_0[1_{\Gamma_t} \, M_t^{\nu_{*}}],\, \Gamma_t \in
{\cal F}_t, \; t \ge 0$, and $(M_t^{\nu_{*}})$ is the
$P_0$-martingale defined by (\ref{mil13}).

\begin{theo} \label{tQ4}
\begin{enumerate}
    \item Under $Q_0^{\nu_{*}}, \; L_{\infty}^0$ is
infinite, $X_{\infty}^*$ is a finite r.v., and the distribution of
$X_{\infty}^*$ is $\nu_*$.

    \item For any $t>0$, under $Q_0^{\nu_{*}}, \; |X_t| <
X_{\infty}^*$.

    \item Let $({\cal G}_t)$ be the smallest filtration satisfying the usual conditions
     such that for any $t\geq 0$, ${\cal G}_t$ contains ${\cal F}_t\vee \sigma(X_\infty ^*)$. Then there
     exists $(\widetilde{B}_t)_{t \ge 0}$ a
$\big(Q_0^{\nu_{*}},({\cal G}_t) \big)$-Brownian motion started at
$0$, such that :
\begin{equation}\label{Q12}
    X_t=\widetilde{B}_t-\int_0^t \frac{{\rm sgn}(X_s) }{X_\infty ^* -|X_s|}
ds.
\end{equation}
\noi Moreover :
\begin{equation}\label{Q12a}
    X_\infty ^*-|X_t|=X_\infty ^*+\widehat{B}_t + \int_0^t \frac{ds }{X_\infty ^* -|X_s|}
-L^0_t,
\end{equation}
where $(\widehat{B}_t)_{t \ge 0}$ is the
$\big(Q_0^{\nu_{*}},({\cal G}_t) \big)$-Brownian motion : $\dis
\widehat{B}_t=-\int_0 ^t {\rm sgn}(X_s)d\widetilde{B}_s$.
\end{enumerate}

\end{theo}

\noi Theorem \ref{tdir1} in section \ref{dir} will generalize
Theorem \ref{tQ4} replacing $(X_{\infty}^*, L_{\infty}^0)$ by
$(S_{\infty}, I_{\infty}, L_{\infty}^0)$.

\noi The proof of Theorem
 \ref{tQ4} is divided in two parts, which are separated by Lemma \ref{lQ1}.

\bigskip

\noi {\bf First part of the proof of Theorem \ref{tQ4}}

\noi For simplicity $M$ stands for $M^{\nu_{*}}$, and $Q_0$ for
$Q_0^{\nu_{*}}$. Recall :
\begin{equation}\label{Q13}
    M_t= \int_{X_{t}^{*}}^\infty \Big(1-\frac {|X_t|}{ a}\Big) \, e^{L_{t}^{0}/a}
\nu_* (da)= \int_0^\infty \Big(1-\frac {|X_{t\wedge T_a^*}|}{
a}\Big) \, e^{L_{t\wedge T_a^*}^{0}/a} \nu_* (da), \quad t \ge 0.
\end{equation}

\noi Let us start with the proofs of 1) and 2).

\noi  a) Let $t$ and $c$ be two positive real numbers. We have :
$$Q_0(X_t^* >c)=Q_0(T_c^*<t)=E_0 [1_{\{T_c^* <t\}} M_t]=E_0
[1_{\{T_c^* <t\}}M_{T_{c}^{*}}].$$

\noi By (\ref{Q13}) we have :
\begin{equation}\label{Q14}
    M_{T_{c}^{*}}=\int_c^{\infty} \Big(1-\frac{c }{ a}\Big)
e^{L_{T_{c}^{*}}^0/a} \nu_* (da).
\end{equation}
\noi Then,
$$Q_0(X_{\infty}^* >c)=\int_c^{\infty} \Big(1- {c \over a}\Big)
\,E_0[e^{L_{T_{c}^{*}}^0/a}] \, \nu_* (da).$$
\noi The r.v. $L_{T_{c}^{*}}^0$ is exponentially distributed :
\begin{equation}\label{Q15}
    P_0(L_{T_{c}^{*}}^0>\lambda)=e^{-\lambda/c}.
\end{equation}

\noi Consequently if $a>c$, then $E_0[e^{L_{T_{c}^{*}}^0 /a}]=\dis
\frac{a }{ a-c}$ and $Q_0(X_{\infty}^* >c)=\nu_*
\big(]c,+\infty[\big)$. This proves that the law of $X_{\infty}^*$
under $Q_0$ is $\nu_*$.

\noi  b) By a similar method we get :
$$Q_0(L_{\infty}^0 >l)=E_0[M_{{\tau}_{l}}]=\int_0^{\infty}
e^{l/a} \, P_0(X_{\tau_{l}}^*<a) \nu_*(da).$$

\noi Since $P_0(X_{\tau_{l}}^*<a)=P_0(T_{a}^{*}>\tau_l)=
P_0(L_{T_{a}^{*}}^0 >l)=e^{-l/a}$, then $Q_0(L_{\infty}^0>l)= \dis
\int_0^{\infty} \nu_* (da)=1$, for any $l \ge 0$. This implies
that $Q_0(L_{\infty}^0=\infty)=1$.

\noi  c) Let $g={\rm sup}\{t \ge 0, \; |X_t|=X_{\infty}^*\}$. Let
$t>0$ be fixed and $\sigma$ be the stopping time : $\sigma= {\rm
inf}\{s>t ; |X_s|>X_t^*\}$. We claim that $Q_0(\sigma <
\infty)=1$, this will imply that $g=\infty$, $Q_0$ a.s.

\noi Let $n\geq 1$. Then :
\begin{equation}\label{Q16}
    Q_0(\sigma \le t+n)=E_0[1_{\{\sigma \le t+n\}} M_{t+n}]=
E_0[1_{\{\sigma \le t+n\}} M_{\sigma}].
\end{equation}

\noi Observing that $X_{\sigma}^*=|X_{\sigma}|=X_t^*$ and $P_0
(\sigma<+\infty)=1$ and taking $n \to \infty$, in (\ref{Q16}), we
obtain :
$$Q_0(\sigma < \infty)=E_0[M_\sigma]=\int_0^{\infty} E_0\Big[1_{\{a >X_{t}^{*}\}}
\Big(1-\frac{X_{t}^{*} }{ a}\Big)\,e^{L_\sigma^{0}/a}\Big] \nu_*
(da). $$

\noi On $\{\sigma < T_0 \circ \theta_t\}, \; L_{\sigma}^0= L_t^0$,
hence :
$$E_0 \big[e^{L_{\sigma}^0/a} 1_{\{\sigma < T_{0} \circ
\theta_t\}}|{\cal F}_t\big]=e^{L_{t}^{0}/a} {|X_t| \over
X_t^*}\cdot$$

\noi On $\{\sigma >T_0 \circ \theta_t\}$, we have $L_{\sigma}^0
=L_t^0 + L_{T_z^*}^0 \circ \theta_{{T_{0}} \circ \theta_t}$ with
$z=X_t^*$, consequently :
$$
\begin{array}{ccl}
 \dis E_0 \big[e^{L_{\sigma}^0/a} 1_{\{\sigma>T_{0} \circ
\theta_{t}\}}|{\cal F}_t\big] & = &\dis \Big(1-\frac{|X_t| }{
X_t^*}\Big) \,
e^{L_{t}^0/a} \int_0^{\infty} e^{l/a} e^{-l/z}\frac{dl }{ z} \\
   & = & \dis {X_t^*-|X_t| \over X_t^* \big(1-{X_{t}^{*}\over a }\big)} \,
e^{L_t^0/a}. \\
\end{array}
$$

\noi Finally :
$$
\begin{array}{ccl}
  Q_0(\sigma < \infty) & = & \dis E_0 \Big[\int_{X_{t}^{*}}^{\infty}
\Big(1-\frac{|X_t| }{ a}\Big)\, e^{L_t^0/a} \nu_* (da)\Big] \\
   & = & E_0[M_t]=1. \\
\end{array}
$$

\vskip 1cm

\noi Point 3) of Theorem \ref{tQ4} will be proved via the initial
enlargement of the original filtration with the r.v.
$X_{\infty}^*$. To apply Theorem 1 of \cite{Y1}   we need  to
compute $Q_0\big(X_{\infty}^*
>c|{\cal F}_t\big), \; c>0, t>0$.

\begin{lemma} \label{lQ1}

\noi Let $c>0$, and $t\ge 0$, then :

\begin{eqnarray}
  Q_0\big(X_{\infty}^*\geq c | {\cal F}_t\big) & = & \dis \int_c^{\infty} e^{L_{t
\wedge T_{c}^{*}}^0/a} \Big(1-\frac{|X_{t \wedge T_{c}^{*}}|}{
a}\Big) \,\frac{1 }{ M_{t \wedge T_{c}^{*}}} \nu_* (da) \label{Q17} \\
   & = & \nu_* \big([c,+\infty[\big)+\int_0^t \widetilde{\lambda}_s
dB_s, \label{Q17a}
\end{eqnarray}

\noi with
\begin{equation}\label{Q18}
\widetilde{\lambda}_s=\frac{-{\rm sgn}(X_s)}{ M_s^2} 1_{\{
X_{s}^{*} < c\}} \int_c^{\infty} e^{L_{s}^{0}/a} \Big\{\frac{1}{
a} \int_{X_{s}^{*}}^{\infty} e^{L_s^0/b} \nu_* (db)-
\int_{X_{s}^{*}}^{\infty}\frac{e^{L_{s}^{0}/b} }{ b} \nu_*(db)
\Big\}  \nu_* (da),
\end{equation}
\begin{equation}\label{Q19}
    B_t=X_t-\int_0^t J_s^{\nu_{*}}  ds,
\end{equation}
\noi $(J_s^{\nu_{*}})$ being defined by (\ref{mil17}).

\end{lemma}

\begin{prooff} \ {\bf of Lemma \ref{lQ1}}
1) Let $\Gamma_t \in {\cal F}_t$ and $u>t$. We decompose
$Q_0\big(\Gamma_t \cap \{X^*_u \geq c\}\big)$ as follows :
$$Q_0\big(\Gamma_t \cap
\{X^*_u \geq c\}\big) =Q_0 \big(\Gamma_t \cap \{T^*_c\leq t
\}\big)+Q_0\big(\Gamma_t \cap \{T^*_c > t, \, X_u^* \geq
c\}\big).$$
Obviously : $\Gamma_t \cap \{T^*_c > t, \, X_u^* \geq
c\}=\Gamma_t' \cap \{T_c^* \leq u\}$ where $\Gamma_t'=\Gamma_t
\cap \{T_c^*>t\}\in {\cal F}_{T_{c}^{*} \wedge t}$. Consequently :
$$
Q_0\big(\Gamma_t' \cap \{T_c^* \leq u\}\big) = E_0[1_{\Gamma_t'
\cap \{T_c^*\leq u\}} M_u] =E_0 [1_{\Gamma_t' \cap \{T_c^* \leq
u\}} M_{T_c^*}].$$
\noi Using (\ref{Q13}) to determine $M_{T_c^*}$ and  taking $u \to
\infty$, we get :
$$
Q_0\big(\Gamma_t \cap \{X_{\infty}^*\geq c\}\big)=  Q_0
\big(\Gamma_t \cap \{T_c^* \leq t\}\big)  + \int_c^{\infty}
\Big(1-\frac{c }{a}\Big) E_0 [1_{\Gamma_t'}\,
e^{L_{T_{c}^{*}}^0/a} ] \nu_* (da).$$
\noi To compute $E_0[1_{\Gamma_t'} \, e^{L_{T_{c}^{*}}^0/a}]$ we
proceed as in part c) of the proof of the first part of Theorem
\ref{tQ4}.

\noi  On $\{ T_0 \circ \theta_t > T_c^* \circ \theta_t\} \cap
\Gamma_t'$, we have $L_{T_{c}^{*}}^0=L_t^0$ and :
$$E_0[e^{L_{T_{c}^{*}}^0/a} |{\cal F}_t]=e^{L_t^0/a} {|X_t|
\over c} \cdot$$

\noi On $\{T_0 \circ \theta_t < T_c^* \circ \theta_t\} \cap
\Gamma_t'$, we have $L_{T_{c}^{*}}^0=L_t^0+L_{T_{c}^{*}}^0 \circ
\theta_{T_{0}\circ \theta_t}$.

\smallskip

\noi Using moreover (\ref{Q15}) we deduce that on $\Gamma_t'$,
$$E_0\big[e^{L_{T_{c}^{*}}^0/a} 1_{\{T_0 \circ \theta_t < T_c^*
\circ \theta_t\}} |{\cal F}_t\big]= e^{L_{t}^{0}/a}
\Big(\frac{c-|X_t| }{ c}\Big) \, \frac{1 }{ 1-c/a} .$$

\noi Then :
$$E_0\big[1_{\Gamma_t'}\ e^{L_{T_{c}^{*}}^0/a}\big]=\frac{1 }{ 1-c/a}
E_0\Big[1_{\Gamma_t'}\ e^{L_{t}^{0}/a} \Big(1- \frac{|X_t| }{
a}\Big)\Big],$$
$$
\begin{array}{ccl}
  Q_0 \big(\Gamma_t \cap \{X_{\infty}^* \geq c\}\big) & = & \dis Q_0\big(\Gamma_t \cap
\{T_c^* \leq t\}\big) + \int_c^{\infty} E_0 \Big[1_{\Gamma_t'}\
e^{L_{t}^{0}/a} \Big(1- \frac{|X_t| }{ a}\Big)\Big] \nu_* (da) \\
   & = & \dis Q_0 \big(\Gamma_t \cap \{T_c^* \leq t\}\big) +E_{Q_{0}}
   \Big[\int_c^{\infty} 1_{\Gamma_t \cap \{T_{c}^{*}>t\}}\
e^{L_{t}^{0}/a} \Big(1-\frac{|X_t| }{ a}\Big)\, \frac{1 }{ M_t}
\nu_* (da)\Big].\\
\end{array}
$$

\noi The previous expression may be simplified, using (\ref{Q14})
:
$$Q_0\big(\Gamma_t \cap \{X_{\infty}^* \geq c\}\big)=E_{Q_{0}} \Big[\int_c^{\infty}
1_{\Gamma_t} \, e^{L_{t \wedge T_{c}^{*}}^0/a} \Big(1-\frac{|X_{t
\wedge T_{c}^{*}}|}{ a}\Big)\,\frac{1 }{ M_{t \wedge T_{c}^{*}}}
\nu_* (da)\Big].$$

\noi This leads to (\ref{Q17}).

\smallskip

\noi 2) Let $Z^1$ and $Z^2$ be two  $\big(Q_0,({\cal F}_t)\big)$
continuous semimartingales. We write :
\begin{equation}\label{Q20}
    Z^1 \stackrel{fv}{\equiv} Z^2
\end{equation}

\noi if $Z^1-Z^2$ is a continuous process with finite variation
and $Z_0^1=Z_0^2$.

\noi It will be convenient to use this notion of congruence, which
we shall apply as follows : if $(M_t)$ is a continuous local
martingale, with $M_0=0$, and $M \stackrel{fv}{\equiv} 0$, then
$(M_t)$ is identically $0$.

\noi Due to  Theorem \ref{tmil1}, the process $(B_t)$ defined by
(\ref{Q19}) is a $\big(Q_0,({\cal F}_t)\big)$- Brownian motion
started at $0$. In particular :
\begin{equation}\label{Q21}
    X_t \stackrel{fv}{\equiv} B_t.
\end{equation}

\noi We have successively :
\begin{equation}\label{Q22}
    \Big(1-\frac{|X_t| }{ a}\Big) \, e^{L_{t}^{0}/a} \stackrel{fv}{\equiv} 1-
    \frac{1}{
a} \int_0^t {\rm sgn}(X_s) \, e^{L_{s}^{0}/a} dB_s,
\end{equation}
\begin{eqnarray}
M_t &\stackrel{fv}{\equiv} &\dis \int_0^{\infty}
\Big\{1-\frac{1}{a} \int_0^t 1_{\{X_{s}^{*} <a\}} {\rm sgn}(X_s)
\, e^{L_{s}^{0}/a} dB_s \Big\}
\nu_* (da) \nonumber\\
&\stackrel{fv}{\equiv} &\dis 1-\int_0^t {\rm sgn}(X_s)
\Big(\int_{X_{s}^{*}}^\infty \frac{e^{L_{s}^{0}/a}}{a} \nu_* (da)
\Big) dB_s. \label{Q23}
\end{eqnarray}
\noi Let $Z^1$ and $Z^2$ be two $\big(Q_0,({\cal F}_t)\big)$
semimartingales, the classical rule of stochastic calculus gives :
$$d \frac{Z_t^1 }{ Z_t^2} \stackrel{fv}{\equiv} \frac{1 }{ Z_t^2} dZ_t^1 -
\frac{Z_{t}^{1} }{ (Z_t^2)^2} dZ_t^2.$$

\noi Choosing $Z_t^1 =\dis\Big(1-\frac{|X_t| }{ a}\Big)\,
e^{L_{t}^{0}/a}$ and $Z_t^2=M_t$, we obtain :
\begin{equation}\label{Q24}
\Big(1-\frac{|X_t| }{ a}\Big) \frac{e^{L_{t}^{0}/a}}{ M_t}
\stackrel{fv}{\equiv} 1-\int_0^t {\rm sgn}(X_s) \,e^{L_{s}^{0}/a}
\Big\{ \frac{1 }{ a M_s}-\Big(1-\frac{|X_s|}{ a}\Big) \frac{1}{
M_s^2} \int_{X_{s}^{*}}^{\infty} \frac{e^{L_{s}^{0}/b} }{ b}
\nu_*(db)\Big\} dB_s.
\end{equation}

\noi Since $\big(Q_0(X_{\infty}^*\geq c|{\cal F}_t)\big)_{t \ge
0}$ is a $\big(Q_0,({\cal F}_t)\big)$-martingale, the identity
(\ref{Q17}) implies :
$$
Q_0\big(X_{\infty}^* \geq c |{\cal F}_t\big) = \int_c^{\infty}
\Big[1-\int_0^t {\rm sgn}(X_s) \,e^{L_{s}^{0}/a} 1_{\{X_{s}^{*}
<c\}} $$
 $$\times \Big\{\frac{1}{ aM_s}  -\Big(1-\frac{|X_s|}{ a}\Big)
\frac{1 }{ M_{s}^{2}} \int_{X_{s}^{*}}^{\infty}
\frac{e^{L_{s}^{0}/b} }{ b} \nu_* (db)\Big\}dB_s\Big] \nu_* (da)
.$$

\noi As a result $Q_0\big(X_{\infty}^* \geq c|{\cal F}_t\big)=
\nu_* \big([c,+\infty[\big)+\dis\int_0^t \widetilde{\lambda}_s
dB_s$, where :
$$\widetilde{\lambda}_s=\frac{-{\rm sgn}(X_s) }{ M_{s}^{2}}
1_{\{X_s^*<c\}} \int_c^{\infty} e^{L_{s}^{0}/a} \Big\{\frac{M_s }{
a}-\Big(1-\frac{|X_s|}{ a}\Big) \int_{X_{s}^{*}}^{\infty}
\frac{e^{L_{s}^{0}/b}}{ b} \nu_*(db)\Big\} \nu_* (da).$$

\noi Formula (\ref{Q17a}) is a direct consequence of (\ref{Q13}).
\end{prooff}
\bigskip


\noi {\bf Second part of the proof    of Theorem \ref{tQ4} (i.e.
point 3) }

\noi Obviously Lemma \ref{lQ1}  may be written as follows :
$$Q_0\big(X_{\infty}^* \geq c|{\cal F}_t\big)=\lambda_t (1_{[c,
\infty[}),$$

\noi with the kernel $\lambda_t$ satisfying  :
\begin{equation}\label{Q25}
    \lambda_t \big([c, \infty[\big)=\nu_* \big([c, \infty[\big)+
\int_0^t \lambda_s^0 \big([c,\infty[\big) dB_s,
\end{equation}
 \noi where :
$$
\lambda_s^0(da)=-\frac{{\rm sgn}(X_s)}{ M_s^2} 1_{\{X_{s}^{*}
<a\}} e^{L_{s}^{0}/a} \Big[\frac{1}{ a} \int_{X_{s}^{*}}^{\infty}
e^{L_{s}^{0}/b}\Big(1-\frac{a }{ b}\Big)\nu_* (db)\Big] \nu_*(da).
$$
\noi The relation (\ref{Q17}) directly implies that :
$$\lambda_t (da)=\frac{1 }{ M_t} 1_{\{X_{t}^{*}<a\}} e^{L_{t}^{0}
/a} \Big[1-\frac{|X_t| }{ a}\Big]\nu_* (da).$$
\noi Consequently :
\begin{equation}\label{Q26}
    \lambda_t^0(da)=\Theta (a,t) \lambda_t (da),
\end{equation}
\noi with :
\begin{equation}\label{Q27}
   \Theta(a,t)=-\frac{{\rm sgn}(X_t) }{ M_t(a-|X_t|)} \; \int_{X_{t}^{*}}^{\infty}
e^{L_{t}^{0}/b} \Big(1-\frac{a }{ b}\Big) \nu_*(db).
\end{equation}

\noi The relations (\ref{Q25})- (\ref{Q27}) allow to apply Theorem
1 of \cite{Y1} (see also \cite{Jac}) : there exists
$(\widetilde{B}_t)_{t \ge 0}$ a $\big(Q_0,({\cal
G}_t)\big)$-Brownian motion started at $0$ such that :
$$B_t=\widetilde{B}_t+\int_0^t \Theta (X_{\infty}^{*},s) ds,$$
\noi where $({\cal G}_t)$ denotes the smallest filtration
satisfying the usual conditions
     such that for any $t\geq 0$, the $\sigma$-field ${\cal F}_t\vee \sigma(X_\infty ^*)$
      is included in ${\cal G}_t$.

 \noi  Using moreover (\ref{Q19}),  $(X_t)$ is seen to admit the
 decomposition : $X_t=\widetilde{B}_t +\dis\int_0^t \xi_s
\, ds$, with  :
$$
\xi_s =\Theta (X_{\infty}^*,s)+J_{s}^{\nu_{*}} . $$
\noi Using successively (\ref{Q27}), (\ref{mil17}) and
(\ref{Q13}), we get :
$$
\begin{array}{ccl}
  \xi _s  &= &\dis  \frac{-{\rm sgn}(X_s) }{ M_s \big(X_{\infty}^*-|X_s|\big)} \Big[
\int_{X_{s}^{*}}^{\infty} e^{L_{s}^{0}/b}
\Big(1-\frac{X_{\infty}^* }{ b}\Big) \nu_* (db)
+\big(X_{\infty}^*-|X_s|\big) \int_{X_{s}^{*}}^{\infty}
e^{L_{s}^{0}/b} \frac{1}{ b} \nu_* (db)\Big]\\
   & = & \dis \frac{-{\rm sgn}(X_s )}{ M_s \big(X_{\infty}^*-|X_s|\big)}
\int_{X_{s}^{*}}^{\infty} \Big(1-\frac{|X_s| }{ b}\Big)
e^{L_{s}^{0}/b} \nu_* (db)=\dis \frac{-{\rm sgn}(X_s )}{ X_{\infty}^*-|X_s|}. \\
\end{array}
$$
\noi This proves (\ref{Q12}).

\noi Applying the  Tanaka formula : $|X_t| =\int_0^t {\rm
sgn}(X_s)dX_s+L^0_t$ and the relation (\ref{Q12}) lead to
(\ref{Q12a}).

\hfill $\blacksquare$

\subsection{ The $Q$-process associated with the down-crossings}
\label{ssQ3}

\noi In this last section we are interested in the law of $(X_t)$
under $Q_{0}^{\downarrow,G}$. We have already introduced some
notation concerning down-crossings in section \ref{do}. Let us
briefly recall the main objects involved in this study. $D_t$ is
the number of down-crossings from $b$ to $a$ ($b>a$), up to time
$t$, and $\big(G(n)\big)_{n \ge 0}$ is a decreasing sequence of
positive real numbers, satisfying $G(0)=1$ and $\dis \lim_{n \to
\infty} G(n)=0$. Let $\big(M_{t}^{\downarrow,G}\big)_{t \ge 0}$ be
the positive $P_0$-martingale associated with $\big(G(n) \big)_{n
\ge 0}$ as defined in Proposition \ref{pdo1}. $Q_{0}^{\downarrow,
G}$ will denote the p.m. on the canonical space $(\Omega , {\cal
F}_\infty)$: $Q_{0}^{\downarrow, G} (\Gamma_t)=E_0 [1_{\Gamma_t}
\, M_{t}^{\downarrow, G}], \; \Gamma_t \in {\cal F}_t$.

\begin{theo} \label{tQ5}
\begin{enumerate}
    \item Under $Q_0^{\downarrow, G}$, the r.v.
$D_{\infty}$ is finite and
\begin{equation}\label{Q30}
    Q_0^{\downarrow, G}(D_{\infty} =n)=G(n)-G(n+1) \;;\; n \ge 0.
\end{equation}

    \item Let $g$ and $\overline{g}$ be the two random times :
$$g=\inf\{t \geq 0 ; D_t=D_{\infty}\}\;;\;\overline{g}= \inf\{ t
>g; X_t=b \}.$$
Then : $Q_0^{\downarrow, G}(0<g<\infty)=1$ and
    $Q_0^{\downarrow, G}(\overline{g} < \infty)=1/2$.

    \item Under $Q_0^{\downarrow, G}$ and conditionally on $\{\overline{g} < \infty\}$ :
\begin{enumerate}
    \item $(X_u; 0 \leq u \le g)$, $(X_{g+u}; 0 \le u \le
\overline{g}-g)$ and $(X_{\overline{g} +u} ; u \ge 0)$ are
independent,

    \item $(2b-a-X_{g+u} ; 0 \le u \le \overline{g}-g)$ is distributed
     as a three dimensional Bessel process, started at $2(b-a)$ and stopped at its first hitting time
of level $b-a$,
    \item $(X_{u+\overline{g}}-a ; u \ge 0)$ is distributed
as a three dimensional Bessel process started at $b-a$.
\end{enumerate}

\item Under $Q_0^{\downarrow, G}$ and conditionally on
$\{\overline{g} = \infty\}$ :

\begin{enumerate}
   \item $(X_u; 0 \leq u \le g)$, and $(X_{g+u}; u \geq 0)$  are
independent,
    \item $(2b-a-X_{g+u} ;u \ge
0)$ is distributed  as a three dimensional Bessel process, started
at $2(b-a)$ and conditioned to be greater than $b-a$.
\end{enumerate}

\item Under $Q_0^{\downarrow, G}$ and conditionally to
$\{D_{\infty}=n\}$, $(X_u ; 0 \leq u \le g)$ is distributed
 as a Brownian motion started at $0$, and stopped at the
first time $(\sigma_{2n})$ when the number of down-crossings
equals $n$.
\end{enumerate}

\end{theo}

\begin{proof} \ a) Recall that $(\sigma_n)_{n \ge 0}$ is the sequence of stopping
times defined by (\ref{do1})-(\ref{do4}). Using the definition of
$Q_0^{\downarrow, G}$ and the optional stopping theorem, we obtain
:
$$
\begin{array}{ccl}
  Q_0^{\downarrow, G} (D_t \ge n) & = & Q_0^{\downarrow, G}
(\sigma_{2n} \le t)=E_0 [1_{\{\sigma_{2n} \le t\}}
M_{t}^{\downarrow, G}]  \\
   & = & E_0 [1_{\{\sigma_{2n} \le t\}}
M_{\sigma_{2n}}^{\downarrow,D} ]. \\
\end{array}
  $$

\noi But relation (\ref{do8}) implies that
$M_{\sigma_{2n}}^{\downarrow,G} =G(n)$. Consequently, taking $t
\to \infty$, in the previous identity  leads to :
$Q_0^{\downarrow, G}(D_{\infty} \ge n)=G(n)$. Hence
$Q_0^{\downarrow, G} (D_{\infty}=n)=G(n)-G(n+1)$. Since
$\dis\lim_{n\rightarrow\infty}G(n)=0$, then $Q_0^{\downarrow,
G}(D_{\infty} < \infty)=1$.

\noi b) The proof of 2-5 of Theorem \ref{tQ5} makes use of the
progressive enlargement of filtrations technique. Since we have
already developed this approach in the setting of Theorems
\ref{tQ1} and  \ref{tQ3}, we limit ourselves to state the main
steps without giving details.

\noi We have :
$$Z_t:= Q_0^{\downarrow,G}(g>t|{\cal F}_t)={G(1+D_{t}) \over M_{t}^{\downarrow,G}}
\;;\; \forall t \ge 0.$$
\noi Theorem \ref{tdo1} says that the process $(B_t)$ defined by
relation (\ref{do16}) is a $Q_0^{\downarrow,G}$- Brownian motion,
started at $0$.

\noi It is convenient for the sequel to introduce :
$$
m_s=\frac{G(1+D_s)-G(D_s) }{ 2(b-a)} \big(
1_{\{\sigma_{2D_{s}}+T_{b} \circ \theta_{\sigma_{2Ds}}
>s\}} -1_{\{ \sigma_{2D_{s}}+T_{b} \circ \theta_{\sigma_{
2Ds}}<s\}}\big)  .$$
\noi Consequently :
$$X_t=B_t+\int_0^t \frac{m_{s} }{
M_{s}^{\downarrow,G}} ds,$$
$$M_{t}^{\downarrow, G}=M_{0}^{\downarrow, G}+ \dis \int_0^t m_s \,
dX_s .$$

\noi Applying It\^{o}'s formula, we obtain :
$$
\begin{array}{ccl}
  Z_t & = & \dis  1-\int_0^t \frac{G(1+D_s) }{ (M_{s}^{\downarrow, G})^2}
m_s d  X_s+\int_{[0,t]} \frac {d_sG(1+D_s) }{ M_s^{\downarrow, G}}\\
   & = & \dis  1-\int_0^t \frac{G(1+D_s) }{ (M_{s}^{\downarrow, G})^2}
m_s d  B_s-\int_0^t \frac{G(1+D_s) }{ (M_{s}^{\downarrow, G})^3}
(m_s)^2 d  B_s+ \int_{[0,t]} \frac {d_sG(1+D_s) }{ M_s^{\downarrow, G}}.\\
\end{array}
$$

\noi Let $({\cal G}_t)_{t \ge 0}$ be the smallest filtration,
containing $({\cal F}_t)_{t \ge 0}$, satisfying the usual
conditions and such that $g$ is a $({\cal G}_t)_{t \ge
0}$-stopping time. Then there exists a $\big(Q_0^{\downarrow,
G},\,({\cal G}_t)_{t \ge 0}\big)$-Brownian motion
$(\widehat{B}_t)$ started at $0$, such that :
$$B_t=\widehat{B}_t + \int_0^{t \wedge g}\frac{d<Z, X>_s }{
Z_s} - \int_{t \wedge g}^t \frac{d<Z, X >_ s }{ 1-Z_ s} .
$$

\noi Since :
$$d < Z, \, X>_s=-{G(1+D_s) \over (M_s^{\downarrow,G})^2}\,
m_s \, ds,$$
\noi we have :
$$X_t=\widehat{B}_t+ \int_{t \wedge g}^t {m_s \over M_s^{\downarrow, G}
- G(1+D_s)} ds.$$

\noi It is easy to compute the drift term, via (\ref{do7}) :
$$\frac{m_s }{ M_s^{\downarrow, G}-G(1+D_s)} =
\left \{ \begin{array}{l} \dis   \frac{1 }{ X_s+a-2b} \ \mbox{on}
\ \{\sigma_{2D_{s}}+T_b \circ \theta_{\sigma_{2D_{s}}} >s\},
\\
\dis \frac{1 }{X_s -a} \ \mbox{ otherwise.}
\end{array}
  \right .$$

\noi In particular :
$$X_t = \widehat{B}_t + \int_{t \wedge g}^{t \wedge \overline{g}}
{ds \over X_s +a-2b} + \int_{t \wedge \overline{g}}^t {ds \over
X_s -a} \;;\; t \ge 0.$$

\noi Having obtained  these results,   points 2-5 of Theorem
\ref{tQ5} can now be proved. The details are left to the reader.

\end{proof}

\section{ A direct approach to study the canonical process under
$Q$}\label{dir}

\setcounter{equation}{0}
 \noi To explain the goal of this section let us start with
Case 1. For simplicity we restrict ourselves to $x=0$. Theorem
\ref{tQ2} leads us to consider $Q^{(y)}_0$, the law of $(X_t)_{t
\ge 0}$ conditionally on $S_{\infty}= y$.  Recall that under
$Q^{(y)}_0$,
\begin{itemize}
    \item $(X_t; t \le T_y)$ is a Brownian motion started at
$0$, and considered up to  its first hitting time  of $y$,
    \item $(y-X_{t + T_{y}} \;;\; t \ge 0)$ is a three
dimensional Bessel process started at $0$.
\end{itemize}
 \noi Therefore, Theorem \ref{tQ2} may be summarized as follows :
\begin{equation}\label{dir1}
    Q_0^{\varphi}(\cdot ) = \int_0^{\infty} Q^{(y)}_0(\cdot ) \; \varphi
(y) \, d y.
\end{equation}

\noi This motivated us to prove  (\ref{dir1}) directly without any
enlargement of filtration.

\begin{prop}\label{pdir1}
 Let $\varphi : [0,+\infty[
\longrightarrow \mathbb{R}_+$ as in Proposition \ref{posm1},
satisfying (\ref{osm3}). Then (\ref{dir1}) holds.
\end{prop}

\begin{proof} \  1) Let $0 \le a < b$ and $\varphi (x)= \dis
\frac{1 }{ b-a} 1_{[a,b]} (x)$. Then :
$$\Phi (x)=\int_{0}^x\varphi (t) dt= \frac{x \wedge b-a }{
b-a} 1_{\{x \ge a\}}.$$
\noi Consequently :
$$M_t^{\varphi}=1_{\{S_{t}<a\}}+{b-X_t \over b-a}\,1_{\{a
\le S_t \le b\}} \;;\; t \ge 0.$$

\noi Then using Girsanov's theorem and Proposition \ref{posm1}, it
is easy to check that $(X_t)_{t \ge 0}$ solves :
$$X_t = \beta_t - \int_0^t \frac{1 }{ b-X_u} 1_{\{a \le S_u
\le b\}} du,$$

\noi where $(\beta_t ; t \ge 0)$ is a $Q_0^{\varphi}$-Brownian
motion started at $0$.

\noi Consequently, under $Q_0^{\varphi}$ :
\begin{itemize}
    \item $(X_t \;;\; 0 \le t \le T_a)$ is distributed as a Brownian
motion started at $0$, and stopped at its first hitting of $a$,
    \item $(b-X_{t+T_{a}} \;;\; t \ge 0)$ is distributed as a three
dimensional Bessel process started at $b-a$.
\end{itemize}

\noi 2) To prove (\ref{dir1}), it is convenient to give an
adequate description of the p.m.   $\dis \int_0^{\infty} Q^{(y)}_0
(\cdot) \varphi (y) d y$.
  Let $\xi$ be a r.v.
independent of $(X_t)_{t \ge 0}$ and uniformly distributed on
$[a,b]$. Then under $\dis \int_0^{\infty} Q^{(y)}_0(\cdot)\varphi
(y) d y $:
\begin{itemize}
    \item $ (X_t ; 0 \le t \le T_{\xi})$
is distributed as $(X_t ; 0 \leq t\le T_{\xi})$ under $P_0$,
    \item $(\xi-X_{t+T_{\xi}} ; t \ge 0)$ is distributed as a three
dimensional Bessel process started at $0$.
\end{itemize}

\noi In particular $(X_t ; t \le T_a)$ has the same distribution
 under either  $Q_0^{\varphi}$ or $\dis \int_0^{\infty} Q^{(y)}_0
 (\cdot) \,\varphi (y) d y$.

\noi 3) Let $R_{b-a}$ be a three dimensional Bessel process
started at $b-a$, and $g$ be the unique time $t$ such that $\dis
R_{b-a} (t)= \inf_{u \ge 0} R_{b-a} (u)$. Then :
\begin{itemize}
    \item $\dis \big(R_{b-a} (t+g) -\inf_{u \ge 0} R_{b-a} (u)  ;
t \ge 0\big)$ is independent of $\big(R_{b-a} (t) ; 0\leq t \le
g)$ and is distributed as a three dimensional Bessel process
started at $0$,
    \item $\dis \inf_{u \ge 0}^{} R_{b-a} (u)$ is uniformly
distributed on $[0, b-a]$, and conditionally on $\dis \inf_{u \ge
0} R_{b-a} (u)=x, \; \big(R_{b-a} (t); 0 \leq t \le g \big)$ is
distributed as a Brownian motion started at $b-a$, and stopped
when it reaches $x$.
\end{itemize}

\noi This implies that the law of $(X_{t+T_{a}} ; t \ge0)$ is the
same  under either $Q_0^{\varphi}$ or $\dis \int_0^{\infty}
Q^{(y)}_0 (\cdot)\,\varphi (y) dy$.

\noi 4) Let $\varphi : \mathbb{R}_+ \mapsto \mathbb{R}_+$ be a
Borel function such that $\dis \int _0^\infty \varphi (y) dy =1$.
It is clear that  $Q_0^{\varphi}$ is equal to
$\widetilde{Q}_0^{\varphi}$, with :
$$\widetilde{Q}_0^{\varphi}(\Gamma_t):=E_0\big[1_{\Gamma_t}\big((S_t-X_t)\varphi
(S_t)+\int_{S_t}^\infty \varphi (y) dy \big)\big], t \geq 0,
\Gamma_t \in {\cal F}_t. $$

\noi  But  we have proved that the two p.m.
$\widetilde{Q}_0^{\varphi}$ and $\dis \int_0^{\infty} Q^{(y)}_0
 (\cdot)\,\varphi (y) dy$
coincide when $\varphi (x)= \dis \frac{1 }{ b-a} 1_{[a,b]} (x)$.
Since  $\psi \longrightarrow \widetilde{Q}_0^{\psi}$ and $\psi
\longrightarrow \dis \int_0^{\infty} Q^{(y)}_0
 (\cdot) \psi (y) dy$ are linear (with respect to convex combinations), then
 \begin{equation}\label{dir1b}
    \dis \widetilde{Q} _0^{\varphi}(\cdot ) = \int_0^{\infty}
Q^{(y)}_0(\cdot ) \; \varphi (y) \, d y,
\end{equation}
holds for any elementary function $\varphi$ of the type : $\varphi
= \dis \sum _{i} c_{i} 1_{[a_{i}, b_{i}]}$. Approximating a
continuous function with compact support by a sequence of stepwise
constant functions of the previous type, permits to prove
(\ref{dir1}), $\varphi$ being a positive continuous function with
compact support in $[0,+\infty[$. The Riesz representation theorem
allows to extend (\ref{dir1b})  and (\ref{dir1}) to any Borel and
positive $\varphi$, the details are left to the reader.

 \end{proof}

\noi We have already observed that if we take $\lambda =0$, then
the Kennedy martingale $M^{0,\varphi}$ coincides with
$M^{\varphi}$. This leads us to generalize Proposition
\ref{pdir1}. We give (see Proposition \ref{pdir2} below) a direct
proof of Theorem \ref{tQ2a} via a disintegration of the p.m.
$Q_0^{\lambda, \varphi}$. Formally taking $\lambda \rightarrow 0$
in Proposition \ref{pdir2} permits to recover Proposition
\ref{pdir1}. However the proofs of Propositions \ref{pdir2} and
\ref{pdir1} are different.

\noi Let $\lambda >0, \; x>0$ and $Q^{\lambda, (x)}_0$ be the
unique p.m. on the canonical space such that :

\begin{itemize}
    \item $(X_t ; 0 \leq t \le T_x)$ is distributed as a Brownian motion
with drift $\lambda$, started at $0$, and stopped at its first
hitting time of $x$,
    \item $(x-X_{t+T_{x}}; t \ge 0)$ is distributed as
$(Z_t^{(\lambda)} ; t \ge 0)$ under $P_0$, where :
\begin{equation}\label{dir2a}
    Z_t^{(\lambda)} =X_t + \lambda \int_0^t {\rm coth}(\lambda
\, Z_u^{(\lambda)}) du.
\end{equation}
\end{itemize}

We observe that :
\begin{equation}\label{dir2b}
    Q_0^{(x)}=\lim_{\lambda \rightarrow 0}Q^{\lambda, (x)}_0.
\end{equation}
%

\begin{prop}\label{pdir2}
 Let $\varphi, \Phi $ be the functions defined by (\ref{osm6b3}), resp. (\ref{osm6b2}),
and parametrized by the function $\psi$ satisfying $\dis
\int_0^\infty \psi (z)e^{-\lambda z} dz=1$. Then :
\begin{equation}\label{dir3}
  Q^{\lambda, \varphi}_0(\cdot)= \int_0^{\infty}Q^{\lambda, (y)}_0(\cdot) e^{- \lambda y} \psi (y) dy.
\end{equation}

\end{prop}

\begin{proof} \ Let $\Gamma \in {\cal F}_\infty$. Using the
definition of $Q_0^{\lambda , \varphi}$ and (\ref{osm6h}), we have
:
\begin{equation}\label{dir4a}
 Q_0^{\lambda , \varphi}(\Gamma)=  E_{P^{(-\lambda)}_x}\Big[\frac{\psi (S_\infty)}
{2\lambda} e^{\lambda S_\infty} 1_\Gamma\big].
\end{equation}
Conditioning on $S_\infty$ and using (\ref{osm6i}) leads to :
$$
Q_0^{\lambda , \varphi}(\Gamma)=\int_0^\infty e^{-\lambda
y}\psi(y) P^{(-\lambda)}_x (\Gamma | S_\infty = y) dy.
$$

\noi The equality :  $Q^{\lambda, (y)}_0(\Gamma)= P^{(-\lambda)}_x
(\Gamma | S_\infty = y)$ follows directly from the well-known
theorem of Williams \cite{Wi}.


\end{proof}

\noi We now investigate the law of $(X_t)$ under $Q_0^{\nu}$,
where $\nu$ is a p.m. on $[\alpha ,\infty[\times [\alpha
,\infty[$, for some $\alpha
>0$. The p.m.  $Q_0^\nu$  has been already introduced in Theorem
\ref{tmil1} and $(M_t^\nu)$ is the $P_0$-martingale defined in
Proposition \ref{pmil1}. Note that the penalization result has
been proved (see Theorem \ref{tmil1}) for the triplet
$(S_t,I_t,L_t^0)$ but we have only described the law of $(X_t)$
under $Q_0^{\nu_*}$; this p.m. being associated with the two
dimensional process $(X^*_t=S_t \vee I_t,L^0_t)$. Hence Theorem
\ref{tdir1} below, generalizes Theorem \ref{tQ4}. Moreover, its
proof hinges on a disintegration of $Q^{\nu}_0$ and does not use
enlargement of filtration .

\noi Let $Q^{s,i}_0$ be the p.m. on $(\Omega ,{\cal F}_\infty)$
defined as the law of $(Y_t^{s,i})$ under $P_0$, where $s,i \geq
0$ and $(Y_t^{s,i})$ is the solution of the following SDE :
\begin{equation}\label{dir5}
    Y_t=X_t-\int_0^t\frac{1}{s-Y_u}1_{\{Y_u>0\}}du +
    \int_0^t\frac{1}{i+Y_u}1_{\{Y_u<0\}}du .
\end{equation}
\noi It will cause no confusion to keep  the same letter $Q$ to
designate  the p.m used in both Proposition \ref{pdir2} above  and
 Theorem \ref{tdir1} below, since the
first p.m.'s is always indexed by $(\lambda ,(y))$ and the second
p.m. by $(s,i)$.

\noi Coming back to (\ref{dir5}), F. Knight (\cite{Kn}) already
considered the process $(Y_t^{s,s})$  and  called it the "Brownian
taboo process". Intuitively it is a Brownian motion conditioned on
never reaching $\pm s$, the taboo levels.

\noi Likewise, from (\ref{dir5}), it can be proved that :
\begin{itemize}
    \item the process $(Y_t^{s,i})$ takes its values in $]-i,s[$,
    \item $\dis \sup_{t\geq 0}Y_t^{s,i} =s$ and $\dis \inf_{t\geq 0}Y_t^{s,i}
    =-i$.
\end{itemize}
\noi These properties may be also proved via the classification of
boundary points $s,-i$ of the diffusion process $(Y_t^{s,i})$ (see
for instance (\cite{RoWi}, section V 50-51).

\begin{theo} \label{tdir1}
\begin{enumerate}
    \item Under $Q_0^{\nu}, \; L_{\infty}^0$ is
infinite, $S_{\infty}$, $I_{\infty}$ are  finite r.v., and the
distribution of $(S_{\infty}, I_{\infty})$ is $\nu$.

    \item We have :
    \begin{equation}\label{dir6}
Q_0^{\nu}(\cdot) =\int
_{([0,+\infty[)^2}Q^{s,i}_0(\cdot)\nu(ds,di).
\end{equation}
\end{enumerate}
\end{theo}

\begin{rem} \label{rdir1} The probabilistic interpretation of
the disintegration property (\ref{dir6}) is the following :
conditionally on $S_{\infty}=s , I_{\infty}=i$, the law of $(X_t)$
under $Q_0^{\nu}$ coincides with the law of $(Y_t^{s,i})$ under
$P_0$.

\end{rem}
%

\begin{prooff} \ {\bf of Theorem \ref{tdir1}} 1) To determine the law of $(S_{\infty},
I_{\infty})$ and prove that $L^0 _\infty \equiv   \infty $ a.s.,
we can proceed as in the proof of Theorem \ref{tQ4}. The details
are left to the reader.

\noi 2) Let us prove (\ref{dir6}). Let $s,i \geq 0$ fixed, and $b$
be the function :
$$
b(x)=-\frac{1}{s -x}1_{\{0<x<s\}}+\frac{1}{i+x}1_{\{-i<x<0\}} .$$
\noi The Girsanov theorem implies that :
$$
Q^{s,i}_0\  _{|{\cal F}_t}=\Theta_t 1_{\{t<T_s\wedge T_{-i}\}}P_0
\ _{|{\cal F}_t},$$
where
$$ \Theta_t =\exp\Big\{\int_0 ^t b(X_u)dX_u -\frac{1}{2}\int_0 ^t
b(X_u)^2du \Big\}.$$
Applying the It\^{o}-Tanaka formula we obtain :
$$
\begin{array}{ccl}
  \ln(s-X_t^+) & = & \dis \ln(s) -\int_0^t\frac{1}{s-X_u}1_{\{X_u>0\}}dX_u
-\frac{1}{2}\int_0^t\frac{1}{(s-X_u)^2}1_{\{X_u>0\}}du
-\frac{1}{2s}L^0_t,
  \\
  \ln(i-X_t^-) & = &  \dis \ln(i) +\int_0^t\frac{1}{i+X_u}1_{\{X_u<0\}}dX_u
-\frac{1}{2}\int_0^t\frac{1}{(i+X_u)^2}1_{\{X_u<0\}}du
-\frac{1}{2i}L^0_t,\\
\end{array}
$$
where $t<T_s \wedge T_{-i}$.

\noi It follows that :
$$ \Theta_t =\Big(1-\frac{X_t^+}{s}\Big)\Big(1-\frac{X_t^-}{i}\Big)
\exp\Big\{ \frac{1}{2}\big(\frac{1}{s}
+\frac{1}{i}\big)L^0_t\Big\}.$$
Summarizing previous calculations, we get :
$$
Q^{s,i}_0\  _{|{\cal F}_t}=
\Big(1-\frac{X_t^+}{s}\Big)\Big(1-\frac{X_t^-}{i}\Big) \exp\Big\{
\frac{1}{2}\big(\frac{1}{s} +\frac{1}{i}\big)L^0_t\Big\} 1_{\{ t <
T_s \wedge T_{-i}\}}\ P_0 \
 _{|{\cal F}_t}.$$
\noi  Integrating with respect to $\nu(ds,di)$ implies directly
(\ref{dir6}), since the martingale $(M^\nu _t)$ is defined by
(\ref{mil5}).
\end{prooff}

\noi Our Theorem \ref{tdir1} provides a short  proof of Lemma
\ref{lmil1}. It is actually possible to demonstrate directly Lemma
\ref{lmil1}, however tedious calculations are necessary.

\noi Let $a,b>0, \nu (ds, di):=\delta_a (ds)\otimes \delta_b (di)$
and $M_t=M_t^\nu$.

\noi Applying the definition (\ref{mil5}) of $(M_t)$, we have :
$$M_t =\frac{1}{ab}(a-X_t^+)(b-X_t^-)e^{cL_t^0}1_{\{t\leq T_a
\wedge T_{-b}\}}; \ t \geq 0,$$
\noi with $\dis c=\frac{1}{2}\big(\frac{1}{a}+\frac{1}{b}\big)$.

\noi Due to the definition of $Q_0^{a,b}$, we get :
$$
\begin{array}{ccl}
 \dis  E_{Q_0^{a,b}}\Big[\frac{1}{(a-X_t^+)(b-X_t^-)}\Big] & = &\dis
  E_0 \Big[\frac{1}{(a-X_t^+)(b-X_t^-)}M_t\Big]\\
   & = & \dis  \frac{1}{ab} E_0\big[e^{cL_t^0}1_{\{t\leq T_a
\wedge T_{-b}\}}\big].\\
\end{array}
$$
\noi Consequently :
$$ \lim_{t\rightarrow \infty}E_0\big[e^{cL_t^0}1_{\{t\leq T_a
\wedge T_{-b}\}}\big]=ab \lim_{t\rightarrow
\infty}E_{Q_0^\nu}\Big[\frac{1}{(a-X_t^+)(b-X_t^-)}\Big].$$
\noi But under $Q_0^\nu$, $(X_t)$ is a recurrent diffusion. It is
easy to compute its invariant density function $g$ since this
function solves :
$$
\left\{
    \begin{array}{cl}
   \dis   \frac{1}{2} g"(x)+\frac{g'(x)}{a-x}+\frac{g(x)}{(a-x)^2}=0& \mbox{ if } \ x>0, \\
   \dis   \frac{1}{2} g"(x)-\frac{g'(x)}{b+x}+\frac{g(x)}{(b+x)^2}=0& \mbox{ if } \ x<0,\\
    \end{array}
    \right.
$$
\noi Finally :
$$g(x)=\frac{3}{(a+b)a^2b^2}\Big[
b^2(a-x)^21_{\{0<x<a\}}+a^2(b+x)^21_{\{-b<x<0\}} \Big].$$
\noi We observe that :
$$\int_\mathbb{R}\frac{1}{(a-x^+)(b-x^-)}g(x)dx =\frac{3}{2ab}.$$
In particular the  integral above is finite and
$$\lim_{t\rightarrow
\infty}E_{Q_0^\nu}\Big[\frac{1}{(a-X_t^+)(b-X_t^-)}\Big]=\frac{3}{2}.$$
This finishes the proof of Lemma \ref{lmil1}, since :
\begin{equation}\label{dir7}
    \lim_{t\rightarrow \infty}E_0\big[e^{cL_t^0}1_{\{t\leq T_a
\wedge T_{-b}\}}\big]=\frac{3}{2}.
\end{equation}

\begin{rem}\label{rdir2} \begin{enumerate}
    \item It is easy to deduce from the previous analysis that
    \begin{equation}\label{dir8}
\lim_{t\rightarrow \infty}E_0\big[e^{\alpha L_t^0}1_{\{t\leq T_a
\wedge T_{-b}\}}\big]= \left \{
\begin{array}{cl}
  \infty & \mbox{ if }\  \alpha > c,\\
  0 & \mbox{ if } \ \alpha < c. \\
\end{array}
 \right.
\end{equation}
    \item Lemma \ref{lmil1} may be generalized as follows. Let
    $(Y_t)$ be Walsh's Brownian motion with parameters $(p_i)_{1 \leq i\leq n}$, where
    $0<p_i<1$ and $\dis \sum_{i=1}^n p_i=1$. We recall (see Walsh's original paper \cite{Wal}, and also
    \cite{BPY1} and \cite{BPY2} for detailed constructions), that this process takes its
    values in a union of half-lines $I_1,\cdots, I_n$ such that
    $\dis \cap_{i=1}^n I_i=\{0\}$. Heuristically, this process :
    \begin{itemize}
        \item moves as a one-dimensional Brownian motion inside
        each $I_i$,
        \item when it reaches $0$, it chooses at random, with
        probability $p_i$, to evolve in  $I_i$ .
    \end{itemize}
    \noi Now, the statement of Lemma \ref{lmil1} may be extended in the following manner :
\begin{equation}\label{dir9}
    \lim_{t\rightarrow \infty}E\Big[e^{\alpha L^0_t}\prod _{i=1}^n
    1_{\{S^i_t \leq a_i\}}\Big]= \left \{
\begin{array}{cl}
  0 & \mbox{ if }\  \alpha < c,\\
\dis 3/2 & \mbox{ if }\  \alpha = c,\\
\infty & \mbox{ if } \ \alpha > c, \\
\end{array}
 \right.
\end{equation}
with $\dis c:=\sum_{i=1}^n \frac{p_i}{a_i}$.
\end{enumerate}
\end{rem}

\section{Further developments} \label{fur}

\setcounter{equation}{0}

\noi In this section, we sketch a number of results which shall
appear in \cite{RVY4}.

\noi For simplicity, we shall only discuss here  Case 1.
\begin{enumerate}
    \item To prove our main limit result, rather
    than considering $E[\varphi(S_t)|{\cal F}_s]$, as $t\rightarrow \infty$,
 we study $P(\Gamma_s|S_t=y]$, which is shown to converge  as $t\rightarrow
\infty$ towards  $Q^{(y)}_0(\Gamma_s)$.
    \item Also, we study the speed of convergence of
    $Q_{0,t}^\varphi(\Gamma_s)$, where $Q_{0,t}^\varphi$ denotes
    the p.m. $Q_{0,t}^F$ defined by (\ref{int6}) with
    $F_t=\varphi(S_t)$. More precisely, suppose that $\varphi$ satisfies
    moreover $\dis \int_0^\infty \varphi (t)t^2 dt<\infty$, we then prove that there exists
    a $P_0$-martingale $(N^{\varphi}_t)$ such that :
\begin{equation}\label{fur1}
    \frac{E\big[1_{\Gamma _s}\varphi(S_t)\big]}{E[\varphi(S_t)]}-
    E\big[1_{\Gamma _s}M^{\varphi}_s\big] _{\stackrel{
    \sim}{t\rightarrow\infty}}\frac{E\big[1_{\Gamma
    _s}N^{\varphi}_s\big]}{t} , \quad \forall \Gamma_s\in {\cal F}_s.
\end{equation}
We obtain more generally   an expansion of the left side of
(\ref{fur1})  in powers of $(1/t)$, to any order.

    \item Finally we consider ( \cite{RVY4}, \cite{RVY3}) the same problems with
    $F_t=f(X_t,S_t)$ where $f :
    \mathbb{R}\times\mathbb{R}_+\mapsto \mathbb{R}_+$. We suppose
    :
    \begin{equation}\label{fur2}
    \overline{f}:=\int_\mathbb{R}dx\int_{x_+}^\infty (2y-x)f(x,y)dy
    <\infty .
\end{equation}
Denote $\widehat{f}=1/\overline{f}$.

\noi  Then we establish that $Q _{0,t}(\Gamma _s)$ converges to
$\dis \widehat{f} E_0\big[1_{\Gamma _s}M^{\varphi}_s\big]$,
$t\rightarrow\infty$, for any  $ \Gamma_s\in {\cal F}_s$ where
$\varphi$ is given by :
$$
\varphi(y) =\widehat{f} \int_\mathbb{R}dx\int_{y\vee x_+}^\infty
f(x,z)dz+ \int_{-\infty}^yf(x,y)(y-x)dx.
$$

\end{enumerate}

\def\refname{References}
\bibliographystyle{plain}
\bibliography{refmaxminRVYII}
\end{document}